\documentclass[10pt,oneside,a4paper,reqno]{amsart}
\usepackage{amsmath,amssymb,amsthm}
\usepackage{bm}
\usepackage{esint}
\usepackage[normalem]{ulem}
\usepackage{comment}

\usepackage[utf8]{inputenc}
\usepackage[T1]{fontenc}
\usepackage[english]{babel}
\usepackage[dvipsnames,svgnames]{xcolor}

\usepackage{newtxtext}

\usepackage[
    a4paper,
    textwidth=15cm,
    hcentering,
    top=1in,
    bottom=1.2in,
    headheight=15pt,
    headsep=0.4in,
    footskip=1cm
]{geometry}
\makeatletter
\def\@tocline#1#2#3#4#5#6#7{\relax
  \ifnum #1>\c@tocdepth 
  \else
    \par \addpenalty\@secpenalty\addvspace{#2}
    \begingroup \hyphenpenalty\@M
    \@ifempty{#4}{
      \@tempdima\csname r@tocindent\number#1\endcsname\relax
    }{
      \@tempdima#4\relax
    }
    \parindent\z@ \leftskip#3\relax \advance\leftskip\@tempdima\relax
    \rightskip\@pnumwidth plus4em \parfillskip-\@pnumwidth
    #5\leavevmode\hskip-\@tempdima
    \ifcase #1
     \or\or \hskip 1em \or \hskip 2em \else \hskip 3em \fi
    #6\nobreak\relax
    \dotfill\hbox to\@pnumwidth{\@tocpagenum{#7}}\par 
    \nobreak
    \endgroup
  \fi}
\def\l@subsection{\@tocline{2}{0pt}{3pc}{5pc}{}}
\makeatother

\definecolor{darkorange}{HTML}{85423e}

\usepackage{enumerate}

\numberwithin{equation}{section}

\theoremstyle{plain}
\newtheorem{theorem}{Theorem}[section]
\newtheorem{proposition}[theorem]{Proposition}
\newtheorem{lemma}[theorem]{Lemma}
\newtheorem{corollary}[theorem]{Corollary}

\theoremstyle{definition}
\newtheorem{definition}[theorem]{Definition}

\theoremstyle{remark}
\newtheorem{remark}[theorem]{Remark}

\newcommand{\R}{\mathbb R}
\newcommand{\C}{\mathbb C}
\newcommand{\Z}{\mathbb Z}
\newcommand{\eps}{\varepsilon}
\newcommand{\Sch}{\mathcal S}
\newcommand{\abs}[1]{\lvert#1\rvert}

\newcommand{\LLL}{\mathcal L_\eps}
\newcommand{\E}{\mathcal E}
\newcommand{\F}{\mathcal F}

\newcommand{\norm}[1]{\left\|#1\right\|}
\newcommand{\ip}[2]{\langle #1,#2\rangle}

\DeclareMathOperator{\Rea}{Re}

\DeclareMathOperator{\dist}{dist}

\usepackage{tikz}
\usetikzlibrary{arrows.meta,patterns,calc}
\usepackage{graphicx}
\usepackage{flafter}
\usepackage{hyperref}
\hypersetup{
    pdftitle={The Triangular Abrikosov Constant and the Thomas--Fermi Gamma-Limit of Rotating Bose--Einstein Condensates},
    pdfauthor={Giovanni Di Fratta and Valeriy V. Slastikov},
    colorlinks=true,
    citecolor=darkorange,
    linkcolor=darkorange,
    urlcolor=darkorange,
    breaklinks=true
}

\title[The triangular Abrikosov constant]{The Triangular Abrikosov Constant and the Thomas--Fermi $\Gamma$-Limit of Rotating Bose--Einstein Condensates}

\author{Giovanni Di Fratta}
\address{\textsc{Giovanni Di Fratta:}  Dipartimento di Matematica e Applicazioni ``R. Caccioppoli'', Universit\`a degli Studi di Napoli ``Federico II'', Via Cintia, 80126 Napoli, Italy.}
\email{giovanni.difratta@unina.it}

\author{Valeriy V. Slastikov}
\address{\textsc{Valeriy V. Slastikov:} School of Mathematics, University of Bristol, Bristol BS8 1TW, United Kingdom.}
\email{valeriy.slastikov@bristol.ac.uk}



\date{}
\begin{document}

\begin{abstract}
We resolve the Abrikosov coefficient conjecture
for rapidly rotating two-dimensional Bose--Einstein
condensates within the variational framework of
$\Gamma$-convergence.
For a harmonic trap, we establish the
$\Gamma$-convergence of the rescaled
Gross--Pitaevskii energy functionals in the
fast-rotation regime, through weak $L^2$
convergence of the rescaled densities, and
identify the effective interaction coefficient
$\beta_*$ with $b_\triangle$, the Abrikosov ratio
associated with the triangular vortex lattice.
The limiting Thomas--Fermi functional determines
the first correction to the minimum energy and
the macroscopic density of almost minimizers.
For the auxiliary rescaled model and its
lowest-Landau-level restriction, we also compute
a strictly positive limiting $L^2$ distance
between minimizing densities and the
Thomas--Fermi profile.

The identification of the effective interaction coefficient $\beta_*$
with the triangular Abrikosov ratio $b_\triangle$ relies on the sharp
Gaussian-minorant theorem established in a recent OpenAI
preprint~\cite[Theorem~1.2]{Atomic}. A Fourier representation 
of the quartic energy
and a weighted Poisson identity yield a lower
bound for arbitrary magnetic-periodic
lowest-Landau-level states.
The sign conditions on the minorant and its
Fourier transform give the bound, while the
lattice contact conditions identify its
coefficient with $b_\triangle$.
Cutoff and magnetic projection transfer this
estimate to local holomorphic states without
boundary conditions or periodicity assumptions.
Consequently, the free-cell coefficients
converge to $b_\triangle$ with error
$O(L^{-1/5})$ on squares of side length $L$.
\end{abstract}

\maketitle

\begingroup
\small
\tableofcontents
\endgroup
\section{Introduction}\label{sec:intro}
We resolve the longstanding conjecture that the
triangular Abrikosov ratio $b_\triangle$ is the
optimal effective interaction coefficient for
rapidly rotating two-dimensional Bose--Einstein
condensates. More precisely, we prove that states
with triangular-periodic density attain the least
limiting interaction cost among arbitrary,
possibly nonperiodic competitors.
We combine this identification with a
$\Gamma$-convergence theorem describing the
limiting energy through the rescaled density.
The coefficient conjecture will be stated precisely
in Subsection~\ref{sec:intro-cell}; see
\cite{ABN06,AS07} and
\cite[Section~7.1]{Hedenmalm2020}
for earlier formulations and related problems.

We consider a condensate confined by the quadratic
potential $V(x)=\abs{x}^2$, commonly called a
harmonic trap. The parameter $\Omega\in(0,1)$
denotes the ratio of the rotation frequency to
the natural oscillation frequency of the trap.
In the rotating frame, the Gross--Pitaevskii
energy is
\begin{equation}\label{eq:original}
  E_\Omega^{\mathrm{GP}}(\Psi)
  =\int_{\R^2}\abs{\nabla\Psi-i\Omega x^\perp\Psi}^2\,\mathrm dx
   +(1-\Omega^2)\int_{\R^2}\abs{x}^2\abs{\Psi}^2\,\mathrm dx
   +\int_{\R^2}\abs{\Psi}^4\,\mathrm dx,
\end{equation}
where $x=(x_1,x_2)\in\R^2$ and
$x^\perp=(-x_2,x_1)$.
We minimize this energy over complex-valued
wave functions of mass one:
\begin{equation}\label{eq:mass-constraint}
  \Psi\in\Sigma
  :=\{u\in H^1(\R^2;\C):
       \abs{x}u\in L^2(\R^2;\C)\},
  \qquad
  \norm{\Psi}_{L^2(\R^2)}=1.
\end{equation}
The associated density $\abs{\Psi}^2$ belongs to
$L^1(\R^2;\R)\cap L^2(\R^2;\R)$, is nonnegative,
and has integral one.
The quartic term is the integral of its square;
it penalizes concentration and models repulsive
interactions, whose coefficient we fix equal
to one.

To express rotation in \eqref{eq:original},
we use the magnetic derivative
$\nabla-i\mathcal A$, with
$\mathcal A(x)=\Omega x^\perp$.
In two dimensions, the corresponding magnetic
field is the scalar curl
\begin{equation}\label{eq:intro-effective-field}
  \operatorname{curl}\mathcal A
  :=\partial_{x_1}\mathcal A_2
    -\partial_{x_2}\mathcal A_1
  =2\Omega,
  \qquad x\in\R^2.
\end{equation}
Here $\mathcal A$ and its curl represent rotation
in the energy and are therefore called an
\emph{effective vector potential} and an
\emph{effective magnetic field}.
Writing the kinetic and rotational terms in this
form leaves the effective confining potential
$(1-\Omega^2)\abs{x}^2$, obtained by subtracting
the centrifugal contribution from the harmonic
potential; see \cite[Section~V.A]{Fet09} and
\cite[Sections~2.1--2.2]{Cooper2008}.
We take \eqref{eq:original} as our starting point.
For many-body derivations of Gross--Pitaevskii
theory, including rotating gases, and for
dimensional reduction, we refer to
\cite{LSY00,LiebSeiringerYngvason2001,LiebSeiringer2006,AB08}. For an account of the rigorous derivation of
Gross--Pitaevskii theory, Bose--Einstein
condensation, and the subsequent developments
in Bogoliubov theory, see the survey
\cite{Schlein2022}.

As $\Omega\uparrow1$, the effective confinement
weakens and the condensate spreads over a growing
region. We will rescale space to describe its
macroscopic density, while retaining the
interaction cost of the finer density variations
imposed by rotation.
The triangular configuration gives an explicit
upper bound for this cost.
Our matching lower bound will use the sharp
Gaussian-minorant theorem proved by OpenAI in
\cite[Theorem~1.2]{Atomic}.

Together with $\Gamma$-convergence, the coefficient
identification will determine the limiting density
functional, the first correction to the minimum
energy, and the macroscopic density profile;
see Subsection~\ref{sec:intro-main}.
The density conclusion will also apply to
almost minimizers, meaning families whose energy
excess is negligible compared with that correction,
as specified in \eqref{eq:intro-almost-minimizers}.
The classification of their microscopic vortex
arrangements remains a separate delicate question.

\subsection{The lowest Landau level and vortices}\label{sec:intro-model}

To understand which wave functions are favored by the
first term in \eqref{eq:original}, we first consider the
magnetic kinetic energy alone. Its minimum under the
mass constraint is $2\Omega$ and it is achieved on 
the \emph{lowest Landau level} (LLL) --- the linear space
of functions in $L^2(\R^2;\C)$ satisfying
\eqref{eq:intro-lll-equation} below.
Its elements belonging to $\Sigma$ and having mass one
are precisely the minimizers of the magnetic kinetic
energy.
We will see that this space defines the local
constraint in our limiting variational problem.

We identify $\R^2$ with $\C$ by $z=x_1+ix_2$ and write
$2\partial_{\bar z}=\partial_{x_1}+i\partial_{x_2}$.
For every $\Psi\in\Sigma$, magnetic factorization gives
\begin{equation}\label{eq:intro-magnetic-factorization}
  \norm{(\nabla-i\Omega x^\perp)\Psi}_{L^2(\R^2)}^2
  =2\Omega\norm{\Psi}_{L^2(\R^2)}^2
   +\norm{(2\partial_{\bar z}+\Omega z)\Psi}_{L^2(\R^2)}^2.
\end{equation}
Thus, when $\norm{\Psi}_{L^2(\R^2)}=1$, the magnetic
kinetic energy is at least $2\Omega$, with equality
if and only if
\begin{equation}\label{eq:intro-lll-equation}
  (2\partial_{\bar z}+\Omega z)\Psi=0
  \qquad\text{in }\mathcal D'(\R^2).
\end{equation}
The normalized Gaussian attains equality, so this lower
bound is the minimum.

To describe the solutions of \eqref{eq:intro-lll-equation},
we set $f=e^{\Omega\abs{z}^2/2}\Psi$.
For $\Psi\in L^2(\R^2;\C)$, we have
$f\in L^2_{\mathrm{loc}}(\C;\C)$, and the product rule gives
\begin{equation}\label{eq:intro-holomorphic-factor}
  2\partial_{\bar z}f
  =e^{\Omega\abs{z}^2/2}
   (2\partial_{\bar z}+\Omega z)\Psi
  \qquad\text{in }\mathcal D'(\R^2).
\end{equation}
Hence \eqref{eq:intro-lll-equation} holds precisely when
$f$ has an entire holomorphic representative. The LLL
therefore consists of the functions
\begin{equation}\label{eq:intro-lll-representation}
  \Psi(z)=e^{-\Omega\abs{z}^2/2}f(z),
  \qquad f\in\mathcal O(\C),
  \qquad
  \int_{\C}e^{-\Omega\abs{z}^2}\abs{f(z)}^2\,\mathrm dx
  <\infty,
\end{equation}
where $\mathcal O(\C)$ is the set of entire functions.
In Subsection~\ref{sec:magnetic}, we will state the
factorization on its full magnetic energy domain;
see also \cite{Bargmann1961,ABN06} and
\cite[Section~1.1]{GGT2019}.

The holomorphic representation also describes the zeros
of LLL states. If $f\not\equiv0$ vanishes at $a\in\C$,
then, in a neighborhood of $a$,
\begin{equation}\label{eq:intro-vortex-factorization}
  f(z)=(z-a)^m g(z),
  \qquad m\in\mathbb N\setminus\{0\},
\end{equation}
where $g$ is holomorphic and $g(a)\neq0$.
The Gaussian factor is strictly positive, so $\Psi$
has the same zeros as $f$. On every sufficiently small
positively oriented circle centered at $a$, the map
$\Psi/\abs{\Psi}$ has degree $m$: its phase makes
$m$ complete turns around the zero. Such a zero is
called a \emph{quantized vortex}; see
\cite[Chapter~4]{Ahlfors1979}, \cite{BBH94}, and
\cite[Section~II]{Fet09}.
In the free-cell problem introduced below, we will
prescribe neither the number nor the positions of
these zeros.

\subsection{Free cells and the coefficient conjecture}
\label{sec:intro-cell}

The quartic term in \eqref{eq:original} measures the
concentration of the density. To determine its contribution
to the limiting energy, we will study a local minimization
problem under the LLL constraint
\eqref{eq:intro-lll-equation}.

For $L>0$, let $Q_L=(-L/2,L/2)^2$ be an open
subset of $\C$.
A dilation reduces the coefficient $\Omega$ in
\eqref{eq:intro-lll-equation} to one.
For $u\in L^2(Q_L;\C)$, the resulting local constraint
$(2\partial_{\bar z}+z)u=0$ in $\mathcal D'(Q_L)$
is equivalent to
$u(z)=e^{-\abs{z}^2/2}f(z)$ almost everywhere on $Q_L$,
with $f\in\mathcal O(Q_L)\cap L^2(Q_L;\C)$.
We call $(Q_L,u)$ a \emph{holomorphic cell} if,
in addition,
$L^{-2}\int_{Q_L}\abs{u}^2\,\mathrm dx=1$.
Thus the density has spatial mean equal to one.

We impose no boundary conditions and require no holomorphic
extension beyond $Q_L$; in particular, we prescribe neither
periodicity nor an arrangement of zeros.
We define
\begin{equation}\label{eq:cell}
  c_L=\inf\left\{
    \frac1{L^2}\int_{Q_L}\abs{u}^4\,\mathrm dx:
    (Q_L,u)\text{ is a holomorphic cell}
  \right\},
\end{equation}
with the quartic integral interpreted as $+\infty$
when it diverges.
The normalization makes $c_L$ the infimum of the spatial
mean of the squared density. By Cauchy--Schwarz, every
admissible value is at least one, with equality only
for a constant density. We will prove the strict
inequality $c_L>1$ in Step~2 of the proof of
Lemma~\ref{lem:cell}: local holomorphic compactness
will rule out even a sequence of admissible densities
whose normalized quartic integrals tend to one.
In Remark~\ref{cert:rem:quotient}, we will interpret
the corresponding quotient as a measure of the
deviation of the density from its spatial mean.

To compare this unrestricted problem with periodic
configurations, we introduce the number
\begin{equation}\label{eq:intro-btri}
  \begin{aligned}
    b_\triangle
      &=\sum_{p\in\Lambda_\triangle}
          e^{-\pi\abs{p}^2}\simeq1.1596,\\
    \Lambda_\triangle
      &=t^{-1/2}\{(m+n/2,tn):(m,n)\in\Z^2\},
      \qquad t=\sqrt3/2.
  \end{aligned}
\end{equation}
Here $\Lambda_\triangle$ is the triangular lattice
normalized so that its fundamental cells have area one.
The number $b_\triangle$ is called the
\emph{triangular Abrikosov ratio}.
In Lemma~\ref{lem:periodic}, we will recall the
theta-function construction of periodic LLL states;
see \cite{ABN06}.
For the triangular configuration, this construction
gives a nonzero state
$u(z)=e^{-\abs{z}^2/2}f(z)$, with
$f\in\mathcal O(\C)$, whose density $\abs{u}^2$
is periodic with respect to a triangular lattice.
If $C\subset\R^2$ is a fundamental cell
of that lattice, then
$\langle\abs{u}^4\rangle/
\langle\abs{u}^2\rangle^2=b_\triangle$,
where
$\langle\abs{u}^p\rangle
:=\abs{C}^{-1}
\int_{C}\abs{u(x)}^p\,\mathrm dx$
for $p\in\{2,4\}$.

In Lemma~\ref{lem:cell}, we will prove
\begin{equation}\label{eq:intro-cell-limit}
  1<c_L\leqslant b_\triangle\quad(L>0),
  \qquad
  \beta_*:=\lim_{L\to\infty}c_L
          =\sup_{L>0}c_L\in(1,b_\triangle].
\end{equation}
The coefficient $\beta_*$ will therefore describe the
optimal quartic interaction on large squares, allowing
all admissible holomorphic functions. The triangular
construction will give the upper bound. The remaining
question is whether arbitrary, possibly nonperiodic,
states can have a smaller limiting energy.

The \emph{coefficient conjecture} asserts that they cannot:
\begin{equation}\label{eq:intro-coefficient-conjecture}
  \beta_*=b_\triangle.
\end{equation}
The fine spatial variations allowed by the local LLL
constraint will contribute the term
$\beta_*\int_{\R^2}\rho^2\,\mathrm dx$ to the limiting
functional \eqref{eq:ThFen}.
Here $\rho \sim |\Psi|^2$ will denote the limiting rescaled density,
introduced in Subsection~\ref{sec:intro-main}.
Thus $\beta_*$ will measure the least interaction
cost compatible with a prescribed macroscopic density.

We prove that this cost is attained asymptotically
by a profile with
triangular-periodic density. More precisely, let
$u_\triangle(z)=e^{-\abs{z}^2/2}f_\triangle(z)$,
with $f_\triangle\in\mathcal O(\C)$, be the nonzero
profile that we will construct
in Lemma~\ref{lem:periodic}.

Although $u_\triangle\notin L^2(\R^2;\C)$, its
restriction to every $Q_L$, after normalization
to mean density one, is a free-cell competitor.
Periodicity of $\abs{u_\triangle}^2$ and its quartic
ratio give
\begin{equation}\label{eq:intro-triangular-optimality}
  \lim_{L\to\infty}
  \frac{L^2\int_{Q_L}\abs{u_\triangle}^4\,\mathrm dx}
       {\left(\int_{Q_L}\abs{u_\triangle}^2\,\mathrm dx\right)^2}
  =b_\triangle.
\end{equation}

Together with $c_L\to\beta_*$, our identification
$\beta_*=b_\triangle$ will show that these normalized
restrictions form an asymptotically minimizing family
for the free-cell problem.
In this sense, the triangular configuration achieves
the optimal limiting interaction cost without any
periodicity assumption on the competitors.

The $\Gamma$-convergence theorem will connect this local
optimization with the macroscopic problem: after minimizing
over the unresolved spatial variations, the remaining energy
will be the Thomas--Fermi functional $\F$, defined in \eqref{eq:ThFen}.
Minimizing $\F$ will then determine the limiting density
profile and the first correction to the minimum
Gross--Pitaevskii energy.
The coefficient identity
will establish the optimality of
the triangular configuration at the level of this limiting
interaction cost; it will not imply
that every minimizing sequence has triangular order.

The prediction that this factor multiplies the
density-squared term in the limiting energy appears in
\cite{Ho01,WBP04,ABD05}; see also
\cite[Section~2.2, equations~(43)--(45)]{Cooper2008}
and \cite[Section~V.A]{Fet09}.
For its mathematical formulation and the distinction
between a periodic trial state and a lower bound for
arbitrary admissible states, see \cite{AB06,ABN06}.
In Theorem~\ref{thm:abrikosov}, we will prove
\eqref{eq:intro-coefficient-conjecture}, together with
the finite-cell estimate \eqref{eq:cellrate}.

\subsection{The variational limit and its consequences}
\label{sec:intro-main}

As $\Omega\uparrow1$, the condensate spreads over a
growing region. We will therefore rescale space before
studying the limiting density.
In this subsection, we denote the spatial variable
in \eqref{eq:original} by $y\in\R^2$ and the rescaled
 variable by $x\in\R^2$. We set
$\eps=\sqrt{1-\Omega}$ and $x=\sqrt\eps\,y$.
For a family $\Psi_\Omega\in\Sigma$ with
$\norm{\Psi_\Omega}_{L^2(\R^2)}=1$, we define
\begin{equation}\label{eq:intro-rescaled-density}
  \rho_\eps(x)
  =\eps^{-1}\abs{\Psi_\Omega(x/\sqrt\eps)}^2,
  \qquad\text{for almost every }x\in\R^2.
\end{equation}
The prefactor preserves mass:
$\int_{\R^2}\rho_\eps\,\mathrm dx=1$.
A fixed region in the $x$ variable corresponds
to a region enlarged by the factor
$\eps^{-1/2}=(1-\Omega)^{-1/4}$ in the $y$ variable.

Our main result identifies the limiting energy
as $\eps\downarrow0$, equivalently $\Omega\uparrow1$.
After subtracting $2\Omega$ and dividing the energy
by $\eps$, we prove $\Gamma$-convergence,
through weak convergence of the rescaled densities
$\rho_\eps$ in $L^2(\R^2;\R)$, to the Thomas--Fermi
functional
\begin{equation}\label{eq:ThFen}
  \F(\rho)
  =2\int_{\R^2}\abs{x}^2\rho\,\mathrm dx
   +\beta_*\int_{\R^2}\rho^2\,\mathrm dx.
\end{equation}
Here $\rho\in L^1(\R^2;\R)\cap L^2(\R^2;\R)$ is
nonnegative almost everywhere on $\R^2$, has
integral one, and satisfies
$\int_{\R^2}\abs{x}^2\rho\,\mathrm dx<\infty$.
We set $\F(\rho)=+\infty$ for all other
$\rho\in L^2(\R^2;\R)$.
Both integrals in \eqref{eq:ThFen} are taken in
the rescaled variable $x$.

We will introduce the rescaled energies in
Subsection~\ref{sec:intro-scaling} and state their
$\Gamma$-convergence through weak $L^2$ convergence
of densities in Theorem~\ref{thm:gamma}.
The proof will follow the three standard steps:
compactness, proved in Section~\ref{sec:compactness};
the $\Gamma$-liminf inequality treated in Section~\ref{sec:lower};
and the construction of recovery sequences in
Section~\ref{sec:recovery}.
These arguments will use the intrinsic coefficient
$\beta_*$ from the free-cell problem.
We will then identify $\beta_*=b_\triangle$ in
Theorem~\ref{thm:abrikosov}, proved in
Section~\ref{sec:certificate}.

Compactness and $\Gamma$-convergence will give
convergence of the minimum energies and ensure
that every weak $L^2$ cluster point of the densities
of almost minimizers minimizes $\F$.
Here an \emph{almost-minimizing family} means a
family of admissible wave functions whose rescaled
energies exceed the corresponding infima by a
quantity tending to zero.
This minimization principle is usually called
the fundamental theorem of $\Gamma$-convergence;
see \cite[Chapter~7]{DalMaso1993} and
\cite[Chapter~1]{Braides2002}.
Since $\F$ has a unique minimizer, the whole
family of rescaled densities will converge weakly
to it in $L^2(\R^2;\R)$.
We will state the minimization conclusions in
Theorem~\ref{thm:min} and prove them, together with
the assertion for almost minimizers, in
Subsection~\ref{sec:min}.

The variational limit will also determine the minimum
energy and the limiting density in the original
variables. We describe these consequences here;
their precise statement, including the condition
on almost minimizers, will be given in
Corollary~\ref{thm:intro} in Section~\ref{sec:main}.

As $\Omega\uparrow1$, the minimum energy will satisfy
\begin{equation}
  \inf_{\substack{\Psi\in\Sigma\\
                  \norm{\Psi}_{L^2(\R^2)}=1}}
      E_\Omega^{\mathrm{GP}}(\Psi)
  =2\Omega+\frac43\sqrt{\frac{2b_\triangle}{\pi}}
      \sqrt{1-\Omega}
   +o\bigl(\sqrt{1-\Omega}\bigr).
\end{equation}
Thus the first correction to the magnetic kinetic
minimum $2\Omega$ will have order $\sqrt{1-\Omega}$.
For mass-one families whose energy excess above the
infimum is $o(\sqrt{1-\Omega})$, the rescaled densities
\eqref{eq:intro-rescaled-density} will converge weakly
in $L^2(\R^2;\R)$ to
\begin{equation}
  \rho_*(x)=\frac{(s-\abs{x}^2)_+}{b_\triangle},
  \qquad s=\sqrt{\frac{2b_\triangle}{\pi}},
  \qquad x\in\R^2,
\end{equation}
where $t_+=\max\{t,0\}$ for $t\in\R$.
This profile is radially symmetric, decreases
quadratically inside the disk of radius $\sqrt{s}$,
and vanishes outside it. It will describe the limiting
mass distribution at the spatial scale
$(1-\Omega)^{-1/4}$ in the original coordinate $y$.

The correction coefficient in
\eqref{eq:intro-energy} will be the minimum of
$\F$ after identifying $\beta_*=b_\triangle$,
and \eqref{eq:intro-TF-profile} will be its unique
minimizer.
We will justify the comparison with the original
energy in Proposition~\ref{prop:original}, proved
in Subsection~\ref{sec:original}.
The deduction of Corollary~\ref{thm:intro} will
be explained
after its statement in Section~\ref{sec:main}.

The weak convergence will describe the mass
distribution on the rescaled spatial scale,
while allowing finer variations of the density.
In Theorem~\ref{thm:min}, we will prove that,
for minimizers, the rescaled densities remain
at a strictly positive $L^2$ distance from
$\rho_*$. Thus these variations will persist
in the quartic energy, even though the limiting
mass profile is uniquely determined.

\subsection{Previous results}\label{sec:intro-sota}
The literature on rotating condensates and vortex
lattices is extensive, and we will not attempt
a comprehensive survey.
We discuss only works closely related to our aims:
the variational fast-rotation limit, the LLL
reduction, and the identification of the effective
interaction coefficient.

The study of periodic vortex lattices began with
Abrikosov's work on the Ginzburg--Landau equations
\cite{Abr57}, followed by the triangular-lattice
comparison in \cite{KRA64}.
Montgomery's theorem \cite{Mon88} establishes
that the triangular lattice minimizes Gaussian
theta sums among planar lattices of fixed covolume.
Since this minimization is restricted to lattices,
it does not yield the lower bound for arbitrary
holomorphic states required in the free-cell
problem \eqref{eq:cell}.

For rapidly rotating condensates, the LLL
description was developed in \cite{Ho01}.
The works \cite{WBP04,ABD05} explained how
distortions of the vortex lattice can produce
a Thomas--Fermi density profile after averaging
over the vortex scale; see also the reviews
\cite{Cooper2008,Fet09}.
The mathematical analysis in \cite{AB06,ABN06}
uses the Gaussian-weighted holomorphic
representation to construct such states and
estimate their energies. These constructions
give upper bounds, while leaving open the
matching triangular lower bound for unrestricted
holomorphic competitors.
The reduction from the three-dimensional
Gross--Pitaevskii model to a two-dimensional
problem, followed by restriction to the LLL,
is justified in \cite{AB08}.

Nguyen and Rougerie \cite{NR22,NR25}
subsequently established minimum-energy
and density asymptotics, with an implicitly
defined thermodynamic coefficient. Its identification with
$b_\triangle$ remained open; see
\cite[Section~1, after equation~(1.14)]{NRv3}.
In our formulation, the coefficient is defined
by $\beta_*:=\sup_{L>0}c_L$, with $c_L$ given
by the free-cell problem \eqref{eq:cell}.
We will prove in Lemma~\ref{lem:cell} that
this supremum is finite and equals
$\lim_{L\to\infty}c_L$.

Our main result will describe the limiting
energy for every admissible density.
Theorem~\ref{thm:gamma} will establish compactness,
the lower bound, and recovery within the LLL.
Theorem~\ref{thm:min} will then give convergence
of the minima and minimizing densities, together
with separate limits of the energy contributions.
These arguments will use only the intrinsic
coefficient $\beta_*$ and impose no periodicity
on the free-cell competitors.

Petrache and Serfaty
\cite[Section~1.2]{PS2020} asked whether the
Cohn--Kumar conjecture could help determine
the minimum of the quartic theta-function
energy arising in high-field superconductivity
and rotating superfluids. They pointed out
that this energy has no evident
pair-interaction representation, so its
relation to universal optimality for point
configurations requires a separate argument.

Our identification $\beta_*=b_\triangle$,
proved in Theorem~\ref{thm:abrikosov}, uses
the sharp Gaussian-minorant theorem
established by OpenAI
in \cite[Theorem~1.2]{Atomic}.
On magnetic tori, we express the quartic
LLL energy as a Gaussian sum with
nonnegative weights. A weighted Poisson
identity then turns the sign conditions
on the minorant and its Fourier transform
into a lower bound. The lattice contact
conditions identify its value with
$b_\triangle$, and cutoff and projection
transfer the bound to free cells.
The variational arguments of
Sections~\ref{sec:proofs}--\ref{sec:assembly}
use only the intrinsic coefficient
$\beta_*$ and are independent of this
identification.

For a model of rotating superfluids,
Serfaty \cite{Serfaty2001} established
vortex formation above a critical angular
velocity, estimates for the minimum energy,
and a limiting free-boundary problem.
Related localization and LLL projection
arguments appear in superconductivity
problems \cite{AS07}.
For the analysis of Ginzburg--Landau vortices,
see \cite{SS07}; for the derivation of a
renormalized interaction energy governing
vortex configurations, see \cite{SS12}.
The coefficient identity considered here
does not classify microscopic minimizing
configurations or establish their stability.
For dynamical questions, we refer to
\cite{GGT2019} on the cubic LLL equation and
to \cite{GST25} on the linear stability of
the triangular Abrikosov lattice and the
instability of rectangular lattices.

Strongly correlated rotating Bose gases
lead to a different variational problem.
Rougerie, Serfaty, and Yngvason
\cite{RSY2014} use the plasma analogy
to estimate the densities and energies
of Laughlin-type trial states.
These are correlated many-body wave
functions, whereas the functional
\eqref{eq:original} describes the
condensate by a single wave function.

\smallskip

\noindent\textbf{Organization of the paper.}
The rest of the paper is organized as follows.
Section~\ref{sec:main} introduces the
functional-analytic setting and rescaling, and
states the main results.
Section~\ref{sec:proofs} provides the analytical
preliminaries.
We prove compactness, the lower bound, and recovery
in Sections~\ref{sec:compactness},
\ref{sec:lower}, and~\ref{sec:recovery}, respectively,
using only the intrinsic definition of $\beta_*$.
Section~\ref{sec:assembly} treats minimization and
the return to the original variables.
Finally, Section~\ref{sec:certificate} identifies
$\beta_*=b_\triangle$ and compares the free-cell
and thermodynamic formulations.

\section{Main results}\label{sec:main}
We now state the variational limit in the variables that keep the
condensate at a fixed macroscopic scale. We first specify the
notation and the domains of the rescaled energies. Theorems~\ref{thm:gamma}
and~\ref{thm:min} then give the variational limit and its consequences
for minimization. Proposition~\ref{prop:original} relates them to
the original energy \eqref{eq:original}.

\subsection{Notation}\label{sec:framework}

Throughout the paper, $U\subseteq\R^2$ denotes an open
set. Wave functions are complex-valued, while densities
are real-valued. We specify the target spaces in
statements and omit them from norms when they are clear
from the context. For vector-valued functions, the
$L^p$ norm is computed using the pointwise Euclidean
length. For example, if $F=(F_1,F_2)\in L^p(U;\C^2)$
and $1\leqslant p<\infty$, then
\begin{equation}
  \norm{F}_{L^p(U)}
  =\left(\int_U
    \bigl(\abs{F_1(x)}^2+\abs{F_2(x)}^2\bigr)^{p/2}
    \,\mathrm dx\right)^{1/p}.
\end{equation}
We identify $x=(x_1,x_2)\in\R^2$ with
$z=x_1+ix_2\in\C$ and write $\mathcal O(U)$ for
the space of holomorphic functions on $U$, viewed
as a subset of $\C$. Derivatives are understood
in the distributional sense unless greater
regularity is stated.

Following the convention commonly used in physics,
we take the complex $L^2$ inner product to be
conjugate-linear in its first argument:
\begin{equation}
  \ip{f}{g}_{L^2(U)}
  =\int_U\bar f g\,\mathrm dx,
  \qquad f,g\in L^2(U;\C).
\end{equation}
We omit the subscript when the underlying Hilbert
space is clear.
We write $\mathcal D(U)=C_c^\infty(U;\C)$ and
denote by $\mathcal D'(U)$ its continuous
complex-linear dual, with $\mathcal D(U)$ equipped
with the usual inductive limit topology.
For $f\in L^1_{\mathrm{loc}}(U;\C)$, the associated
distribution $T_f\in\mathcal D'(U)$ is defined by
\begin{equation}
  \ip{\varphi}{T_f}_{\mathcal D,\mathcal D'}
  =\int_U\varphi f\,\mathrm dx,
  \qquad \varphi\in\mathcal D(U).
\end{equation}
Thus the distributional pairing is complex-linear
in both arguments and involves no complex conjugation.

If $h\in L^1_{\mathrm{loc}}(\R^2;\C)$ is periodic
with respect to a planar lattice and $C\subset\R^2$
is a fundamental cell, we write
\begin{equation}
  \langle h\rangle
  =\frac1{\abs{C}}\int_C h\,\mathrm dx.
\end{equation}
Periodicity makes this average independent of the
chosen fundamental cell.

For $t\in\R$, we set $t_+=\max\{t,0\}$.
For a Lebesgue-measurable set $B\subseteq\R^2$,
we denote its indicator by $\mathbf1_B$ and its
Lebesgue measure by $\abs{B}$. For $R>0$, we write
$B_R=\{x\in\R^2:\abs{x}<R\}$.
The symbols $\rightharpoonup$ and $\to$ denote
weak and strong convergence, respectively, in
the stated space. Constants denoted by $C$ may
change between estimates; their relevant parameter
dependence will be specified.

We conclude by recording two useful representations
of quantities introduced in Section~\ref{sec:intro}.
The triangular Abrikosov ratio \eqref{eq:intro-btri}
can also be written as
\begin{equation}\label{eq:btri}
  b_\triangle=\sum_{k,l\in\Z}
    \exp\left(-\frac{\pi}{t}\abs{l-k\tau}^2\right),
  \qquad \tau=\frac12+i\frac{\sqrt3}{2},
  \qquad t=\Im\tau=\frac{\sqrt3}{2}.
\end{equation}
Indeed, the substitution $(m,n)=(l,-k)$ identifies
this sum with \eqref{eq:intro-btri}.

The normalization in the free-cell problem
\eqref{eq:cell} can be removed by writing the energy
as a quotient. Recalling that $Q_L=(-L/2,L/2)^2$,
we have
\begin{equation}\label{eq:cellquotient}
  c_L=\inf_{\substack{f\in\mathcal O(Q_L)\\
           0<\int_{Q_L}\abs{f}^2e^{-\abs{z}^2}\,\mathrm dx<\infty}}
  \frac{L^2\int_{Q_L}\abs{f}^4e^{-2\abs{z}^2}\,\mathrm dx}
       {\left(\int_{Q_L}\abs{f}^2e^{-\abs{z}^2}\,\mathrm dx\right)^2},
  \qquad L>0.
\end{equation}
We interpret the quotient as $+\infty$ when its
numerator is infinite. It is invariant under
$f\mapsto\lambda f$ for every $\lambda\in\C\setminus\{0\}$.
In particular, replacing $f$ by
\begin{equation}
  \frac{L f}
  {\left(\int_{Q_L}\abs{f}^2e^{-\abs{z}^2}
    \,\mathrm dx\right)^{1/2}}
\end{equation}
gives a Gaussian-weighted function of mean density
one, as required in \eqref{eq:cell}, without changing
the quotient.

\subsection{Rescaling and the lowest Landau level}\label{sec:intro-scaling}
We introduce the spatial and energy rescalings needed
for the $\Gamma$-convergence result. The spatial
rescaling will keep the density on a fixed
scale, while the energy rescaling will isolate the
correction above the magnetic minimum $2\Omega$.
We will then introduce the lowest Landau level,
by which we mean the linear space of wave functions
associated with the lowest magnetic energy;
its precise definition is given in
\eqref{eq:lllintro}.
The gap to the next energy level will allow us to
approximate bounded-energy states by their orthogonal
projections onto this space.

We choose the spatial scale so that confinement
and interaction contribute at the same order.
Let $\rho\in L^1(\R^2;\R)\cap L^2(\R^2;\R)$ be
nonnegative, with mass one and finite second moment.
For $R>0$, the dilation
$\rho_R(x)=R^{-2}\rho(x/R)$ preserves mass and satisfies
\begin{equation}
  \int_{\R^2}\abs{x}^2\rho_R(x)\,\mathrm dx
  =R^2\int_{\R^2}\abs{x}^2\rho(x)\,\mathrm dx,
  \qquad
  \int_{\R^2}\rho_R^2\,\mathrm dx
  =R^{-2}\int_{\R^2}\rho^2\,\mathrm dx.
\end{equation}
Thus increasing $R$ raises the confinement energy
and lowers the interaction energy.
Their scaling factors in \eqref{eq:original} are
$(1-\Omega^2)R^2$ and $R^{-2}$.
Setting $\eps=\sqrt{1-\Omega}$ and choosing
$R=\eps^{-1/2}$ makes these factors equal to
$(1+\Omega)\eps$ and $\eps$, respectively.
Both have order $\eps$, and their ratio tends to $2$.
This motivates dividing the energy above the
magnetic minimum $2\Omega$ by $\eps$, so that
both contributions remain present in the limit.

We implement this rescaling at the level of wave
functions. Writing $y\in\R^2$ for the original
coordinate and $x\in\R^2$ for the rescaled coordinate,
we set
\begin{equation}\label{eq:intro-rescaling}
  \eps=\sqrt{1-\Omega},
  \qquad x=\sqrt\eps\,y,
  \qquad
  \Psi(y)=\sqrt\eps\,\psi(\sqrt\eps\,y),
  \qquad y\in\R^2.
\end{equation}
For $\Omega\in(0,1)$, we have $\eps\in(0,1)$ and
$\Omega=1-\eps^2$.
The map $\psi\mapsto\Psi$ is a bijection of $\Sigma$
onto itself and preserves the $L^2$ norm.
For $\psi\in\Sigma$ with
$\norm{\psi}_{L^2(\R^2)}=1$, changing variables in
\eqref{eq:original} gives the exact rescaled energy
\begin{equation}\label{eq:exactscaled}
  \begin{aligned}
    \widehat\E_\eps(\psi)
    &:=\frac{E_\Omega^{\mathrm{GP}}(\Psi)-2\Omega}{\eps}\\
    &=\norm{(\nabla-i\Omega x^\perp/\eps)\psi}_{L^2(\R^2)}^2
      -\frac{2\Omega}{\eps}\\
    &\quad +(1+\Omega)\int_{\R^2}
          \abs{x}^2\abs{\psi}^2\,\mathrm dx
      +\int_{\R^2}\abs{\psi}^4\,\mathrm dx.
  \end{aligned}
\end{equation}
Here we used $(1-\Omega^2)/\eps^2=1+\Omega$.
We extend $\widehat\E_\eps$ to $L^2(\R^2;\C)$
by setting it equal to $+\infty$ outside the
unit sphere of $\Sigma$.
All integrals in \eqref{eq:exactscaled} are taken
in the rescaled variable $x$.

For the variational proofs, we will replace
\eqref{eq:exactscaled} by an energy in which
$\Omega$ is replaced by its limiting value $1$.
We first define this simpler energy and then state
the estimate that justifies the replacement.
For every $\eps>0$, let
$D_\eps=\nabla-i x^\perp/\eps$ and
$A_\eps=2\partial_{\bar z}+z/\eps$, with
$z=x_1+ix_2$. We define the magnetic energy domain by
\begin{equation}\label{eq:magnetic-energy-domain}
  \mathcal H_\eps
  :=\{v\in L^2(\R^2;\C):
       D_\eps v\in L^2(\R^2;\C^2)\}.
\end{equation}
The factorization \eqref{eq:factor}, proved on this
domain in Subsection~\ref{sec:magnetic}, gives
\begin{equation}\label{eq:Keps}
  K_\eps(v)
  :=\norm{D_\eps v}_{L^2(\R^2)}^2
    -\frac2\eps\norm{v}_{L^2(\R^2)}^2
  =\norm{A_\eps v}_{L^2(\R^2)}^2
  \geqslant0,
  \qquad v\in\mathcal H_\eps.
\end{equation}
Thus $K_\eps$ measures the magnetic kinetic energy
above its minimum per unit mass.
We then define $\E_\eps:L^2(\R^2;\C)\to[0,+\infty]$ by
\begin{equation}\label{eq:modelintro}
  \E_\eps(v)
  =K_\eps(v)
   +2\int_{\R^2}\abs{x}^2\abs{v}^2\,\mathrm dx
   +\int_{\R^2}\abs{v}^4\,\mathrm dx
\end{equation}
for $v\in\Sigma$ with
$\norm{v}_{L^2(\R^2)}=1$, and by $+\infty$ otherwise.
Indeed, the identity
$\nabla v=D_\eps v+i\eps^{-1}x^\perp v$
shows that
\begin{equation}\label{eq:energy-domain}
  \mathcal H_\eps\cap
  \{v\in L^2(\R^2;\C):
       \abs{x}v\in L^2(\R^2;\C)\}
  =\Sigma.
\end{equation}
The quartic integral is finite on $\Sigma$ by the
classical Sobolev embedding in dimension two.

Proposition~\ref{prop:original} will give the following
comparison: for every $M>0$, there exist
$C_M>0$ and $\eps_M>0$ such that
\begin{equation}
  \abs{\widehat\E_\eps(v)-\E_\eps(v)}
  \leqslant C_M\sqrt\eps
\end{equation}
whenever $0<\eps<\eps_M$, $v\in\Sigma$,
$\norm{v}_{L^2(\R^2)}=1$, and at least one of
$\widehat\E_\eps(v)$ and $\E_\eps(v)$ is at most $M$.
Thus the two energies differ by a quantity tending
uniformly to zero on bounded-energy families.
This will allow us to establish the variational
limit for $\E_\eps$ and then transfer it to
$\widehat\E_\eps$.

To analyze $\E_\eps$, we first identify the functions
for which its kinetic contribution vanishes.
By \eqref{eq:Keps}, this amounts to solving the
first-order equation $A_\eps v=0$.
Its $L^2$ solutions form the closed subspace called
the \emph{lowest Landau level}:
\begin{equation}\label{eq:lllintro}
  \begin{aligned}
    \mathcal L_\eps
    &:=\{v\in L^2(\R^2;\C):
          A_\eps v=0\text{ in }\mathcal D'(\R^2)\}\\
    &=\{v\in L^2(\R^2;\C):
          v(z)=e^{-\abs{z}^2/(2\eps)}f(z)
          \text{ a.e. on }\R^2,\quad
          f\in\mathcal O(\C)\}.
  \end{aligned}
\end{equation}
The second description in the previous relation can be obtained as follows. For
$v\in L^2(\R^2;\C)$ we set
$f=e^{\abs{z}^2/(2\eps)}v
\in L^2_{\mathrm{loc}}(\R^2;\C)$.
The product rule gives
\begin{equation}
  2\partial_{\bar z}f
  =e^{\abs{z}^2/(2\eps)}A_\eps v
  \qquad\text{in }\mathcal D'(\R^2).
\end{equation}
By the distributional Cauchy--Riemann theorem,
$A_\eps v=0$ is equivalent to the existence of an
entire holomorphic representative of $f$.
This proves the representation in \eqref{eq:lllintro}.
In Subsection~\ref{sec:magnetic}, we will also recall
that every $v\in\mathcal L_\eps$ belongs to
$\mathcal H_\eps$. Hence \eqref{eq:Keps} gives
$K_\eps(v)=0$ for every such $v$.
Consequently, for a mass-one function
$v\in\Sigma\cap\mathcal L_\eps$, the energy
\eqref{eq:modelintro} reduces to
\begin{equation}
  \E_\eps(v)
  =2\int_{\R^2}\abs{x}^2\abs{v}^2\,\mathrm dx
   +\int_{\R^2}\abs{v}^4\,\mathrm dx.
\end{equation}
We call its minimization over this class the
\emph{LLL-restricted problem}.
It retains the confinement and interaction energies
while requiring the wave function to have the
Gaussian-weighted holomorphic form
\eqref{eq:lllintro}.

To compare this restricted problem with the full
energy $\E_\eps$, we must control the component
of a wave function outside $\mathcal L_\eps$.
Let $P_\eps:L^2(\R^2;\C)\to\mathcal L_\eps$
be the orthogonal projection.
The magnetic Laplacian associated with $D_\eps$
has lowest eigenvalue $2/\eps$, with eigenspace
$\mathcal L_\eps$, and next eigenvalue $6/\eps$;
see \cite[Section~4.1, equation~(24)]{Thiang2025},
with magnetic field $b=2/\eps$.
Thus the component $(I-P_\eps)v$ has kinetic
energy at least $4/\eps$ above the lowest level.
Since $K_\eps$ measures this kinetic excess,
the spectral decomposition gives
\begin{equation}\label{eq:gapintro}
  \norm{(I-P_\eps)v}_{L^2(\R^2)}^2
  \leqslant\frac\eps4K_\eps(v),
  \qquad v\in\mathcal H_\eps.
\end{equation}
In particular, if $\E_\eps(v_\eps)$ remains bounded
as $\eps\downarrow0$, then
$\norm{v_\eps-P_\eps v_\eps}_{L^2(\R^2)}
=O(\sqrt\eps)$.
This gives an $L^2$ approximation by LLL states.
To compare the quartic interaction terms, we will
also need to control the same remainder in $L^4$;
Lemma~\ref{lem:projection} will provide that estimate.

\subsection{Convergence through densities}

The formal scaling argument above suggests that the
limiting energy should depend only on the density.
Indeed, the confinement and interaction terms in
\eqref{eq:modelintro} depend on $\abs{v}^2$, while
the magnetic kinetic term favors states in the LLL.
A nonzero LLL state cannot have a positive \emph{constant}
density on an open set. The free-cell limit
\eqref{eq:intro-cell-limit} quantifies the interaction
cost of the density variations imposed by this
constraint. We will prove that, for a prescribed
limiting density $\rho$ with finite limiting energy,
the least possible limiting quartic interaction is
$\beta_*\int_{\R^2}\rho^2\,\mathrm dx$.
Since $\beta_*>1$, this exceeds the value
$\int_{\R^2}\rho^2\,\mathrm dx$ obtained by ignoring
the fine spatial variations.

This motivates describing convergence of wave functions
through weak convergence of their densities, without
imposing convergence of their phases.
Indeed, both the energy and the
density are unchanged under
$v_\eps\mapsto e^{i\theta_\eps}v_\eps$, where
$\theta_\eps\in\R$ may depend on $\eps$.
A description of the microscopic phases would
require a separate analysis.
For $v_\eps\in L^2(\R^2;\C)\cap L^4(\R^2;\C)$,
$\eps\in(0,1)$, and $\rho\in L^2(\R^2;\R)$,
we write
\begin{equation}\label{eq:density-convergence}
  v_\eps\xrightarrow{\mathrm{dens}}\rho
  \quad\Longleftrightarrow\quad
  \abs{v_\eps}^2\rightharpoonup\rho
  \quad\text{in }L^2(\R^2;\R)
  \quad\text{as }\eps\downarrow0.
\end{equation}
Every such limit $\rho$ is nonnegative almost
everywhere on $\R^2$, since the cone of nonnegative
functions is weakly closed in $L^2(\R^2;\R)$.

We will formulate the variational limit in terms
of $\Gamma$-convergence. For the reader's convenience,
we recall its definition for the convergence
\eqref{eq:density-convergence}.

\begin{definition}[$\Gamma$-convergence through densities]
\label{def:density-convergence}
We say that $\E_\eps$ $\Gamma$-converges through
densities to
$\mathcal G:L^2(\R^2;\R)\to[0,+\infty]$
as $\eps\downarrow0$ if the following conditions hold:
\begin{enumerate}[(i)]
  \item \emph{Lower bound.}
  For every $\rho\in L^2(\R^2;\R)$ and every
  family $v_\eps\xrightarrow{\mathrm{dens}}\rho$,
  \begin{equation}\label{eq:density-Gamma-values}
    \liminf_{\eps\downarrow0}\E_\eps(v_\eps)
    \geqslant\mathcal G(\rho).
  \end{equation}

  \item \emph{Recovery.}
  For every $\rho\in L^2(\R^2;\R)$ with
  $\mathcal G(\rho)<\infty$, there exists a family
  $v_\eps\xrightarrow{\mathrm{dens}}\rho$ such that
  \begin{equation}\label{eq:density-recovery}
    \lim_{\eps\downarrow0}\E_\eps(v_\eps)
    =\mathcal G(\rho).
  \end{equation}
  We call such a family a \emph{recovery family}.
\end{enumerate}
All wave functions in this definition belong to
$L^2(\R^2;\C)\cap L^4(\R^2;\C)$.
\end{definition}

Compactness is not a part of this definition.
We will prove, separately, that bounded-energy
sequences have subsequences whose densities
converge weakly in $L^2(\R^2;\R)$.
Together with $\Gamma$-convergence, this property
will yield convergence of the minimum energies
and identify the limiting densities of families
whose energy approaches the minimum.
More precisely, we call a family
$\psi_\eps\in\Sigma$, $\eps\in(0,1)$, with
$\norm{\psi_\eps}_{L^2(\R^2)}=1$,
\emph{almost minimizing} for $\E_\eps$ if
\begin{equation}\label{eq:almost-minimizing-family}
  \E_\eps(\psi_\eps)
  -\inf_{\substack{v\in\Sigma\\
                   \norm{v}_{L^2(\R^2)}=1}}
      \E_\eps(v)
  \longrightarrow0
  \qquad\text{as }\eps\downarrow0.
\end{equation}
Thus the energy excess tends to zero on the
rescaled energy scale. Exact minimizers satisfy
this condition, but attainment of the infimum
is not required to choose such a family.
We will prove that almost-minimizing families
have the same limiting density as exact minimizers.
This notion is useful when comparing two energies.
A minimizer of $\widehat\E_\eps$ need not minimize
$\E_\eps$ exactly. However,
Proposition~\ref{prop:original} will show that the
two energies differ by $o(1)$ on bounded energy
sublevels, and that their infima also differ by
$o(1)$. Consequently, a family minimizing
$\widehat\E_\eps$ will be almost minimizing for
$\E_\eps$. The density convergence proved for
almost minimizers of $\E_\eps$ will therefore
apply to this family as well.
For the general minimization principle underlying
this argument, see \cite[Chapter~7]{DalMaso1993}
and \cite[Chapter~1]{Braides2002}.

\subsection{The variational limit and minimization}\label{sec:results}

We now state our main result, which determines the
least limiting energy compatible with a prescribed
density $\rho$.
For mass-one wave functions $v_\eps\in\Sigma$,
write $\rho_\eps=\abs{v_\eps}^2$.
The quartic interaction in \eqref{eq:modelintro}
is then
\begin{equation}
  \int_{\R^2}\abs{v_\eps}^4\,\mathrm dx
  =\int_{\R^2}\rho_\eps^2\,\mathrm dx
  =\norm{\rho_\eps}_{L^2(\R^2)}^2.
\end{equation}
Consequently, weak convergence
$\rho_\eps\rightharpoonup\rho$ in $L^2(\R^2;\R)$
does not by itself imply convergence of the
interaction energy to $\int_{\R^2}\rho^2\,\mathrm dx$.

The magnetic kinetic term imposes an additional
restriction on bounded-energy families: they must
approach the LLL in the sense described by
\eqref{eq:gapintro}. We will prove that the resulting
local holomorphic constraint forces the stronger
interaction bound
\begin{equation}
  \liminf_{\eps\downarrow0}
  \int_{\R^2}\abs{v_\eps}^4\,\mathrm dx
  \geqslant
  \beta_*\int_{\R^2}\rho^2\,\mathrm dx,
\end{equation}
whenever $\sup_{\eps\in(0,1)}\E_\eps(v_\eps)<\infty$
and $\abs{v_\eps}^2\rightharpoonup\rho$ in
$L^2(\R^2;\R)$.

The theorem below makes this description precise.
The compactness assertion will then allow us to
apply this variational convergence result to minima and
almost minimizers.

\begin{theorem}[Variational limit through densities]\label{thm:gamma}
The free-cell constants \eqref{eq:cell} satisfy
\begin{equation}\label{eq:celllimit}
  \beta_*:=\lim_{L\to\infty}c_L
           =\sup_{L>0}c_L=b_\triangle.
\end{equation}
Define $\F:L^2(\R^2;\R)\to[0,+\infty]$ as follows.
For every nonnegative
$\rho\in L^1(\R^2;\R)\cap L^2(\R^2;\R)$ with
$\int_{\R^2}\rho\,\mathrm dx=1$ and
$\int_{\R^2}\abs{x}^2\rho\,\mathrm dx<\infty$, set
\begin{equation}\label{eq:limit}
  \F(\rho)
  =2\int_{\R^2}\abs{x}^2\rho\,\mathrm dx
   +\beta_*\int_{\R^2}\rho^2\,\mathrm dx.
\end{equation}
For all other $\rho\in L^2(\R^2;\R)$, set
$\F(\rho)=+\infty$.
Then the following statements hold.
\begin{enumerate}[(1)]
\item \emph{Lower bound.}
For every sequence $\eps_n\downarrow0$, every
$\rho\in L^2(\R^2;\R)$, and every
$v_n\in L^2(\R^2;\C)\cap L^4(\R^2;\C)$ such that
$\abs{v_n}^2\rightharpoonup\rho$ in $L^2(\R^2;\R)$,
\begin{equation}\label{eq:directliminf}
  \liminf_{n\to\infty}\E_{\eps_n}(v_n)\geqslant\F(\rho).
\end{equation}
\item \emph{LLL recovery.}
For every $\rho\in L^2(\R^2;\R)$ with $\F(\rho)<\infty$,
there is a family $\psi_\eps\in\Sigma\cap\mathcal L_\eps$,
$\eps\in(0,1)$, with $\norm{\psi_\eps}_{L^2(\R^2)}=1$, such that
\begin{equation}\label{eq:directrecovery}
  \abs{\psi_\eps}^2\rightharpoonup\rho
     \quad\text{in }L^2(\R^2;\R),
  \qquad \E_\eps(\psi_\eps)\longrightarrow\F(\rho)
     \quad\text{as }\eps\downarrow0.
\end{equation}
\item \emph{Density compactness.}
If $\eps_n\downarrow0$ and $v_n\in L^2(\R^2;\C)$ satisfy
$\sup_n\E_{\eps_n}(v_n)<\infty$, then, after passing to a
subsequence, $\abs{v_n}^2\rightharpoonup\rho$ in
$L^2(\R^2;\R)$ for some nonnegative $\rho$ of mass one
and finite second moment.
\end{enumerate}
Consequently, $\E_\eps$ and its LLL restriction have the same
sequential $\Gamma$-limit $\F$ in the sense of
Definition~\ref{def:density-convergence}.
\end{theorem}

We prove the theorem using the three canonical steps of a 
$\Gamma$-convergence result.
First, we establish compactness, assertion~(3), in
Section~\ref{sec:compactness}.
We then prove the lower bound, assertion~(1), in
Subsection~\ref{sec:liminfproof}, using the averaged
cell inequality and the projection estimates.
Finally, we construct the recovery families required
by assertion~(2) in Subsection~\ref{sec:recoveryproof}:
we approximate the target density by piecewise constant
densities, assemble localized cell profiles, and project
them onto the LLL.

These arguments use only the intrinsic characterization
of $\beta_*$ in Lemma~\ref{lem:cell}.
We identify $\beta_*=b_\triangle$ separately in
Theorem~\ref{thm:abrikosov}, which also gives the
finite-cell estimate \eqref{eq:cellrate}.
This identification uses the Gaussian minorant
provided by \cite[Theorem~1.2]{Atomic}.

We now apply Theorem~\ref{thm:gamma} to the minimization
of the rescaled energy $\E_\eps$.
We consider two problems: minimization over all
mass-one functions in $\Sigma$, and minimization
over those that also belong to $\mathcal L_\eps$.
The latter is a smaller admissible class, so its
minimum cannot be lower. Nevertheless, the LLL
recovery families in Theorem~\ref{thm:gamma}
will show that the two minima have the same limit.
Compactness will then identify the limiting
densities of minimizers of either problem.

\begin{theorem}[Minima and limiting densities]
\label{thm:min}
For every $\eps>0$, the energy $\E_\eps$ attains
its minimum both on
$\{v\in\Sigma:\norm{v}_{L^2(\R^2)}=1\}$
and on
$\{v\in\Sigma\cap\mathcal L_\eps:
\norm{v}_{L^2(\R^2)}=1\}$.
Moreover,
\begin{equation}\label{eq:minvalue}
  \begin{aligned}
    \lim_{\eps\downarrow0}
    \min_{\substack{v\in\Sigma\\
                    \norm{v}_{L^2(\R^2)}=1}}
      \E_\eps(v)
    &=
    \lim_{\eps\downarrow0}
    \min_{\substack{v\in\Sigma\cap\mathcal L_\eps\\
                    \norm{v}_{L^2(\R^2)}=1}}
      \E_\eps(v)\\
    &=e_*:=\min_{\rho\in L^2(\R^2;\R)}\F(\rho)
      =\frac43\sqrt{\frac{2\beta_*}{\pi}}.
  \end{aligned}
\end{equation}
The unique minimizer of $\F$ is
\begin{equation}\label{eq:TF}
  \rho_*(x)=\frac{(s-\abs{x}^2)_+}{\beta_*},
  \qquad s=\sqrt{\frac{2\beta_*}{\pi}},
  \qquad x\in\R^2.
\end{equation}
For each $\eps>0$, let $\psi_\eps$ be a mass-one
minimizer of $\E_\eps$ on $\Sigma$, or a mass-one
minimizer of $\E_\eps$ on $\Sigma\cap\mathcal L_\eps$.
Then
\begin{equation}
  \E_\eps(\psi_\eps)\longrightarrow e_*,
  \qquad
  \rho_\eps:=\abs{\psi_\eps}^2
  \rightharpoonup\rho_*
  \quad\text{in }L^2(\R^2;\R)
  \quad\text{as }\eps\downarrow0.
\end{equation}
The three contributions to $\E_\eps(\psi_\eps)$
also converge separately:
\begin{align}
  K_\eps(\psi_\eps)
    &\longrightarrow0,
      \label{eq:partK}\\
  2\int_{\R^2}\abs{x}^2\rho_\eps\,\mathrm dx
    &\longrightarrow
      2\int_{\R^2}\abs{x}^2\rho_*\,\mathrm dx,
      \label{eq:partV}\\
  \int_{\R^2}\rho_\eps^2\,\mathrm dx
    &\longrightarrow
      \beta_*\int_{\R^2}\rho_*^2\,\mathrm dx.
      \label{eq:partI}
\end{align}
In particular, the densities do not converge
strongly in $L^2(\R^2;\R)$:
\begin{equation}\label{eq:defect}
  \lim_{\eps\downarrow0}
  \norm{\rho_\eps-\rho_*}_{L^2(\R^2)}^2
  =(\beta_*-1)\norm{\rho_*}_{L^2(\R^2)}^2>0.
\end{equation}
\end{theorem}

We prove Theorem~\ref{thm:min} in
Subsection~\ref{sec:min}.
The convergence of the minimum energies and
minimizing densities follows from compactness,
the lower bound, and LLL recovery in
Theorem~\ref{thm:gamma}.
To obtain the separate limits
\eqref{eq:partK}--\eqref{eq:partI}, we also use
the interaction lower bound in
Proposition~\ref{prop:lower}.
Together with convergence of the total energy,
it forces each contribution to attain its
corresponding limiting value.
Finally, weak convergence and
\eqref{eq:partI} give \eqref{eq:defect} by
expanding the squared $L^2$ distance.
Thus the additional factor $\beta_*$ occurs
in the quartic interaction, while the magnetic
kinetic excess tends to zero.

The density conclusion also holds for almost
minimizers. More precisely, suppose that
$\psi_\eps\in\Sigma$,
$\norm{\psi_\eps}_{L^2(\R^2)}=1$, and
\begin{equation}
  \E_\eps(\psi_\eps)
  -\min_{\substack{v\in\Sigma\\
                   \norm{v}_{L^2(\R^2)}=1}}
      \E_\eps(v)
  \longrightarrow0
  \qquad\text{as }\eps\downarrow0.
\end{equation}
Then $\E_\eps(\psi_\eps)\to e_*$.
By compactness and the lower bound in
Theorem~\ref{thm:gamma}, every weak $L^2$ cluster
point of $\abs{\psi_\eps}^2$ minimizes $\F$.
The uniqueness in Theorem~\ref{thm:min} therefore
gives
$\abs{\psi_\eps}^2\rightharpoonup\rho_*$
in $L^2(\R^2;\R)$ for the whole family.
The same conclusion holds for almost minimizers
of the LLL-restricted problem, since its minimum
has the same limit \eqref{eq:minvalue}.

These results concern the simpler rescaled
energy $\E_\eps$. To obtain conclusions for the
original energy $E_\Omega^{\mathrm{GP}}$, we must
first compare $\E_\eps$ with the exact rescaled
energy $\widehat\E_\eps$ defined in
\eqref{eq:exactscaled}.
The following estimate shows that their difference
is negligible on every bounded energy sublevel.

\begin{proposition}[Comparison with the exact rescaled energy]
\label{prop:original}
For every $M>0$, there exist $C_M>0$ and
$\eps_M\in(0,1)$ such that
\begin{equation}\label{eq:comparison}
  \abs{\widehat\E_\eps(\psi)-\E_\eps(\psi)}
  \leqslant C_M\sqrt\eps
\end{equation}
whenever $0<\eps<\eps_M$, $\psi\in\Sigma$,
$\norm{\psi}_{L^2(\R^2)}=1$, and
$\min\{\E_\eps(\psi),\widehat\E_\eps(\psi)\}
\leqslant M$.
Consequently, $\widehat\E_\eps$
$\Gamma$-converges through densities to the
same functional $\F$ as $\E_\eps$.
\end{proposition}

We prove Proposition~\ref{prop:original} in
Subsection~\ref{sec:original}.
The estimate applies to bounded sublevels of
either energy, so it allows us to compare their
infima in both directions.
Combining it with \eqref{eq:minvalue} gives
\begin{equation}
  \lim_{\eps\downarrow0}
  \inf_{\substack{\psi\in\Sigma\\
                  \norm{\psi}_{L^2(\R^2)}=1}}
    \widehat\E_\eps(\psi)
  =e_*.
\end{equation}
It also shows that every almost-minimizing family
for $\widehat\E_\eps$ is almost minimizing for
$\E_\eps$, and hence has densities converging
weakly to $\rho_*$.
Returning to the original variables through
\eqref{eq:intro-rescaling} and using
$\beta_*=b_\triangle$ from
Theorem~\ref{thm:abrikosov}, we obtain the
following corollary.

\begin{corollary}[Fast-rotation energy and density asymptotics]
\label{thm:intro}
Let $E_\Omega^{\mathrm{GP}}$ and $\Sigma$ be defined
by \eqref{eq:original}--\eqref{eq:mass-constraint}.
Then the following statements hold.

\begin{enumerate}[(i)]
\item \emph{Minimum energy.}
As $\Omega\uparrow1$,
\begin{equation}\label{eq:intro-energy}
  \inf_{\substack{\Psi\in\Sigma\\
                  \norm{\Psi}_{L^2(\R^2)}=1}}
    E_\Omega^{\mathrm{GP}}(\Psi)
  =2\Omega+\eps e_*+o(\eps),
\end{equation}
where we recall that
\begin{equation}
  \eps=\sqrt{1-\Omega},
  \qquad
  e_*=\min_{\rho\in L^2(\R^2;\R)}\F(\rho)
     =\frac43\sqrt{\frac{2b_\triangle}{\pi}}.
\end{equation}
Thus the remainder in \eqref{eq:intro-energy},
divided by $\eps$, tends to zero.
\smallskip

\item \emph{Densities of almost minimizers.}
Let $(\Psi_\Omega)_{\Omega\in(0,1)}$ be a family
in $\Sigma$ with
$\norm{\Psi_\Omega}_{L^2(\R^2)}=1$.
Assume that its energy excess above the infimum
is $o(\eps)$, that is,
\begin{equation}\label{eq:intro-almost-minimizers}
    E_\Omega^{\mathrm{GP}}(\Psi_\Omega)=
    \displaystyle\inf_{\substack{\Psi\in\Sigma\\
                       \norm{\Psi}_{L^2(\R^2)}=1}}
       E_\Omega^{\mathrm{GP}}(\Psi)+o(\eps).
\end{equation}
Define
\begin{equation}
  \rho_\eps(x)
  =\eps^{-1}\abs{\Psi_\Omega(x/\sqrt\eps)}^2
  \qquad\text{for almost every }x\in\R^2.
\end{equation}
Then $\rho_\eps\rightharpoonup\rho_*$
in $L^2(\R^2;\R)$ as $\eps\downarrow0$, where
\begin{equation}\label{eq:intro-TF-profile}
  \rho_*(x)=\frac{(s-\abs{x}^2)_+}{b_\triangle},
  \qquad s=\sqrt{\frac{2b_\triangle}{\pi}},
  \qquad x\in\R^2.
\end{equation}
Consequently, for every bounded
Lebesgue-measurable set $B\subset\R^2$,
\begin{equation}\label{eq:intro-mass-regions}
  \int_{\eps^{-1/2}B}
    \abs{\Psi_\Omega(y)}^2\,\mathrm dy
  =\int_B\rho_\eps(x)\,\mathrm dx
  \longrightarrow
    \int_B\rho_*(x)\,\mathrm dx.
\end{equation}
\end{enumerate}
\end{corollary}

Assertion~(i) identifies $\eps e_*$ as the first
correction to the minimum magnetic kinetic energy
$2\Omega$. Assertion~(ii) describes the mass
distribution of families whose energy excess is
negligible compared with this correction.
It applies, in particular, to exact minimizers.

We prove Corollary~\ref{thm:intro} in
Subsection~\ref{sec:original}, after establishing
the comparison with the exact rescaled energy.

\section{Analytical preliminaries}\label{sec:proofs}

We establish the estimates needed for the compactness,
lower-bound, and recovery assertions of
Theorem~\ref{thm:gamma}.
The magnetic factorization \eqref{eq:factor} and
the projection estimates will allow us to approximate
bounded-energy wave functions by LLL states.
Lemma~\ref{lem:translation} will let us translate
and rescale local profiles, while
Lemmas~\ref{lem:periodic} and~\ref{lem:cell}
will establish the bounds and large-square limit
defining the intrinsic coefficient $\beta_*$.
These results will be used in
Sections~\ref{sec:compactness},
\ref{sec:lower}, and~\ref{sec:recovery}.
They do not require the identification
$\beta_*=b_\triangle$.

\subsection{Magnetic identities and projection errors}
\label{sec:magnetic}

To use LLL projection in the variational argument,
we need to specify its relation to the magnetic
energy domain and estimate the resulting errors.
We first recall well-known operator identities, then state
the projection estimates in the form required below.

Fix $\eps>0$ and recall
$D_\eps=\nabla-i x^\perp/\eps$ and
$A_\eps=2\partial_{\bar z}+z/\eps$.
The magnetic energy domain
\eqref{eq:magnetic-energy-domain} can equivalently
be described using either operator:
\begin{equation}
  \begin{aligned}
    \mathcal H_\eps
    &=\{v\in L^2(\R^2;\C):
          D_\eps v\in L^2(\R^2;\C^2)\}\\
    &=\{v\in L^2(\R^2;\C):
          A_\eps v\in L^2(\R^2;\C)\}.
  \end{aligned}
\end{equation}
Here the derivatives are distributional.
On this domain, magnetic factorization gives
\begin{equation}\label{eq:factor}
  \norm{D_\eps v}_{L^2(\R^2)}^2
  =\norm{A_\eps v}_{L^2(\R^2)}^2
   +\frac2\eps\norm{v}_{L^2(\R^2)}^2,
  \qquad v\in\mathcal H_\eps.
\end{equation}
The identity is first established for smooth,
compactly supported functions and extended to
the stated domain by approximation; see
\cite{Friedrichs1944,GGT2019}.

The equation $A_\eps u=0$ also has a local
description that does not require global
square integrability.
For every open set $U\subseteq\R^2$ and every
$u\in L^2_{\mathrm{loc}}(U;\C)$, the product rule
and the distributional Cauchy--Riemann theorem give
\begin{equation}\label{eq:local-lll-representation}
  \begin{aligned}
    A_\eps u=0\text{ in }\mathcal D'(U)
    &\quad\Longleftrightarrow\quad
    \partial_{\bar z}
      \bigl(e^{\abs{z}^2/(2\eps)}u\bigr)
      =0\text{ in }\mathcal D'(U)\\
    &\quad\Longleftrightarrow\quad
    u=e^{-\abs{z}^2/(2\eps)}f
    \text{ a.e. on }U
    \text{ for some }f\in\mathcal O(U).
  \end{aligned}
\end{equation}
Taking $U=\R^2$ and imposing global square
integrability recovers the lowest Landau level
defined in \eqref{eq:lllintro}:
\begin{equation}
  \mathcal L_\eps
  =\{v\in L^2(\R^2;\C):
      v(z)=e^{-\abs{z}^2/(2\eps)}f(z)
      \text{ a.e. on }\R^2,\quad
      f\in\mathcal O(\C)\}.
\end{equation}
By \eqref{eq:factor}, this is precisely the set
of $v\in\mathcal H_\eps$ for which
$K_\eps(v)=\norm{A_\eps v}_{L^2(\R^2)}^2=0$.
We denote its $L^2$-orthogonal projection by
$P_\eps:L^2(\R^2;\C)\to\mathcal L_\eps$.

We next estimate the remainder $(I-P_\eps)v$.
Its $L^2$ norm controls the change in mass,
while its $L^4$ norm is needed to compare the
quartic interaction energies before and after
projection. We will use these estimates in
Propositions~\ref{prop:lower}
and~\ref{prop:steprecover}.

The estimates are the whole-plane counterparts
of Lemma~\ref{cert:lem:projection}.
Here the magnetic field is $2/\eps$ and the
spectral gap is $4/\eps$; on the magnetic tori
used in Subsection~\ref{cert:sec:transfer},
they are $2\pi$ and $4\pi$.
In the torus case, the constants must also be
independent of the side length $S\geqslant1$.
Both proofs combine the spectral gap with
magnetic factorization and an interpolation
inequality applied to the modulus.
On $\R^2$ we use the Gagliardo--Nirenberg
inequality; on the torus, its periodic version
\eqref{cert:eq:ladyzhenskaya} contains an
additional $L^2$ term.
We give only a proof sketch here, since the
argument is standard; the proof of
Lemma~\ref{cert:lem:projection} gives the
corresponding periodic argument in detail.

\begin{lemma}[Projection errors]\label{lem:projection}
Let $\eps>0$, $v\in\mathcal H_\eps$, and
$r=(I-P_\eps)v$. Then
\begin{align}
  \norm{r}_{L^2(\R^2)}
  &\leqslant\frac{\sqrt\eps}{2}
    \norm{A_\eps v}_{L^2(\R^2)},
    \label{eq:r2}\\
  \norm{D_\eps r}_{L^2(\R^2)}
  &\leqslant\sqrt{\frac32}
    \norm{A_\eps v}_{L^2(\R^2)},
    \label{eq:rD}\\
  \norm{r}_{L^4(\R^2)}
  &\leqslant C\eps^{1/4}
    \norm{A_\eps v}_{L^2(\R^2)}.
    \label{eq:r4}
\end{align}
The constant $C>0$ is independent of $v$ and $\eps$.
\end{lemma}
\begin{proof}[Proof sketch]
All norms below are over $\R^2$. Since
$P_\eps v\in\mathcal L_\eps\subset\mathcal H_\eps$, we have
$r\in\mathcal H_\eps$, $r\perp\mathcal L_\eps$, and
$A_\eps r=A_\eps v$.
The gap \eqref{eq:gapintro} gives \eqref{eq:r2}.
Substituting this bound into
$\norm{D_\eps r}_{L^2}^2
=\norm{A_\eps v}_{L^2}^2+(2/\eps)\norm{r}_{L^2}^2$
gives \eqref{eq:rD}.
The diamagnetic inequality gives $\abs{r}\in H^1(\R^2;\R)$;
combined with the planar Gagliardo--Nirenberg inequality
\cite[Section~8.5]{LiebLoss2001}, it yields
\begin{equation}
  \norm{r}_{L^4}^4
  \leqslant C\norm{r}_{L^2}^2\norm{D_\eps r}_{L^2}^2
  \leqslant C\eps\norm{A_\eps v}_{L^2}^4.
\end{equation}
Taking fourth roots proves \eqref{eq:r4}.
\end{proof}
For the recovery construction we must also control the confinement
term $2\int_{\R^2}\abs{x}^2\abs{v}^2\,\mathrm dx$ after projection.
The next lemma bounds the weighted remainder for states supported
in a fixed ball; we use it in Step~2 of the proof of
Proposition~\ref{prop:steprecover}. It requires the Gaussian decay
of the projection kernel in addition to the spectral gap.
Section~\ref{cert:sec:transfer} only compares mass and quartic
interaction on a torus, so it does not require this weighted estimate.
We again sketch the standard kernel argument.
\begin{lemma}[Weighted projection error]\label{lem:kernel}
Let $R>0$ and $\eps>0$. If $v\in\mathcal H_\eps$ is supported in $B_R$, then
\begin{equation}\label{eq:weightedr}
 \norm{\abs{x}(I-P_\eps)v}_{L^2(\R^2)}
 \leqslant C_R\sqrt\eps
 \bigl(\norm{v}_{L^2(\R^2)}+\norm{A_\eps v}_{L^2(\R^2)}\bigr),
\end{equation}
where $C_R$ is independent of $v$ and $\eps$.
\end{lemma}

\begin{proof}[Proof sketch]
Set $r=(I-P_\eps)v$. On $B_{2R}$, we have
$\norm{\abs{x}r}_{L^2(B_{2R})}\leqslant2R\norm{r}_{L^2(\R^2)}$.
Outside $B_{2R}$, we have $r=-P_\eps v$ and use the Bargmann
kernel formula \cite[equation~(1.1)]{GGT2019}, after dilation:
\begin{equation}\label{eq:projection-kernel-modulus}
  \abs{P_\eps(x,y)}
  =\frac1{\pi\eps}e^{-\abs{x-y}^2/(2\eps)},
  \qquad x,y\in\R^2.
\end{equation}
Here $P_\eps(x,y)$ denotes the integral kernel of $P_\eps$.
For $x\notin B_{2R}$ and $y\in B_R$, we have
$\abs{x}\leqslant2\abs{x-y}$.
The kernel $x\mapsto2(\pi\eps)^{-1}\abs{x}e^{-\abs{x}^2/(2\eps)}$
has $L^1(\R^2)$ norm $C\sqrt\eps$, so Young's inequality gives
$\norm{\abs{x}r}_{L^2(\R^2\setminus B_{2R})}
\leqslant C\sqrt\eps\norm{v}_{L^2(\R^2)}$.
Combining the two regions and applying \eqref{eq:r2}
proves \eqref{eq:weightedr}.
\end{proof}
In the following lemma we establish invariance of the local LLL equation under magnetic translation and scaling. It is required in Lemma~\ref{lem:cell} to compare translated cells, in Lemma~\ref{lem:average} to apply the cell inequality at scale $\sqrt\eps$, and in Proposition~\ref{prop:steprecover} to assemble profiles and compute their magnetic excess.

\begin{lemma}[Magnetic translation and scaling]
\label{lem:translation}
Let $U\subseteq\R^2$ be open,
$\bm a=(a_1,a_2)\in\R^2$, and $\eps>0$.
We write $\bm a^\perp=(-a_2,a_1)$ for the
counterclockwise rotation of $\bm a$ by $\pi/2$.
For $h\in L^2_{\mathrm{loc}}(U;\C)$, define
\begin{equation}
  h_{\bm a,\eps}(x)
  =e^{i\bm a^\perp\cdot x/\eps}
    h\bigl((x-\bm a)/\sqrt\eps\bigr),
  \qquad x\in\bm a+\sqrt\eps\,U.
\end{equation}
Then $\abs{h_{\bm a,\eps}(x)}=\abs{h\bigl((x-\bm a)/\sqrt\eps\bigr)}$ for almost every $x\in\bm a+\sqrt\eps\,U$.
Moreover,
\begin{equation}\label{eq:translation}
  \begin{aligned}
    D_\eps h_{\bm a,\eps}(x)
    &=\eps^{-1/2}e^{i\bm a^\perp\cdot x/\eps}
      (D_1h)\bigl((x-\bm a)/\sqrt\eps\bigr),\\
    A_\eps h_{\bm a,\eps}(x)
    &=\eps^{-1/2}e^{i\bm a^\perp\cdot x/\eps}
      (A_1h)\bigl((x-\bm a)/\sqrt\eps\bigr)
  \end{aligned}
\end{equation}
in distributions on $\bm a+\sqrt\eps\,U$.
Here $D_1=\nabla-i x^\perp$ and
$A_1=2\partial_{\bar z}+z$ are the operators
with parameter equal to one.

Consequently, $A_\eps h_{\bm a,\eps}=0$ in
$\mathcal D'(\bm a+\sqrt\eps\,U)$ if and only if
$A_1h=0$ in $\mathcal D'(U)$, or equivalently,
by \eqref{eq:local-lll-representation}, if
$h=e^{-\abs{z}^2/2}f$ almost everywhere on $U$
for some $f\in\mathcal O(U)$.

For every $p\in[1,\infty)$, we also have
\begin{equation}
  \int_{\bm a+\sqrt\eps\,U}
    \abs{h_{\bm a,\eps}(x)}^p\,\mathrm dx
  =\eps\int_U\abs{h(y)}^p\,\mathrm dy.
\end{equation}
If $h$ vanishes almost everywhere outside a
compact subset of $U$, extend $h$ and
$h_{\bm a,\eps}$ by zero to $\R^2$.
The identities \eqref{eq:translation} then hold
in $\mathcal D'(\R^2)$.
If, in addition, $A_1h\in L^2(U;\C)$, then
\begin{equation}
  \norm{A_\eps h_{\bm a,\eps}}_{L^2(\R^2)}
  =\norm{A_1h}_{L^2(\R^2)}.
\end{equation}
\end{lemma}

\begin{proof}
The phase has modulus one, which gives the
pointwise identity.
To compute the derivatives, set
$y=(x-\bm a)/\sqrt\eps\in U$.
The product rule gives
\begin{equation}
  \begin{aligned}
    D_\eps h_{\bm a,\eps}(x)
    &=e^{i\bm a^\perp\cdot x/\eps}
      \left(
        \eps^{-1/2}\nabla h(y)
        +\frac{i(\bm a^\perp-x^\perp)}{\eps}h(y)
      \right)\\
    &=\eps^{-1/2}e^{i\bm a^\perp\cdot x/\eps}
      \bigl(\nabla h(y)-iy^\perp h(y)\bigr),
  \end{aligned}
\end{equation}
where we used $x=\bm a+\sqrt\eps\,y$.
This proves the identity for $D_\eps$.
Taking the first component plus $i$ times the
second proves the identity for $A_\eps$.
These calculations hold in distributions because
multiplication by the smooth phase and an invertible
affine change of variables are valid distributional
operations. The phase compensation is the usual
magnetic translation; see \cite{Zak1964}.

The substitution $x=\bm a+\sqrt\eps\,y$
has Jacobian $\eps$, giving the integral identity.
When $A_1h\in L^2(U;\C)$, the same substitution
in \eqref{eq:translation} yields
\begin{equation}
  \int_{\bm a+\sqrt\eps\,U}
    \abs{A_\eps h_{\bm a,\eps}}^2\,\mathrm dx
  =\int_U\abs{A_1h}^2\,\mathrm dy:
\end{equation}
the Jacobian $\eps$ cancels the squared
derivative factor $\eps^{-1}$.

Finally, if $h$ vanishes outside a compact subset
of $U$, neither zero extension introduces a
boundary distribution. The differential identities
therefore hold on $\R^2$, and the last integral
identity gives the asserted equality of norms.
\end{proof}

\subsection{Holomorphic cells and their large-square limit}\label{sec:cells}
In this section we state and prove results related to the intrinsic coefficient $\beta_*$. Its supremum characterization gives the lower bound in Lemma~\ref{lem:average}, while its limit characterization gives the recovery cost in Proposition~\ref{prop:steprecover}. The strict inequality $\beta_*>1$ implies the density defect in Theorem~\ref{thm:min}. This section is independent of Section~\ref{sec:certificate}, where we identify $\beta_*$ with $b_\triangle$.
\begin{lemma}[Intrinsic cell coefficient]\label{lem:cell}
For every $L>0$, the constants in \eqref{eq:cell} satisfy $1<c_L\leqslant b_\triangle$. The large-square limit exists and satisfies
\begin{equation}\label{eq:intrinsiccell}
 \beta_*:=\sup_{L>0}c_L=\lim_{L\to\infty}c_L,
 \qquad 1<\beta_*\leqslant b_\triangle.
\end{equation}
\end{lemma}

The following periodic profiles provide the comparison functions used later
to prove $c_L\leqslant b_\triangle$.

\begin{lemma}[Periodic Gaussian-weighted profiles]
\label{lem:periodic}
Let $\tau=\tau_1+i\tau_2$ with $\tau_2>0$, and set
$a=\sqrt{\pi/\tau_2}$. There exists a nonzero entire function $f_\tau$
such that
\begin{equation}
  u_\tau(z)=e^{-\abs{z}^2/2}f_\tau(z),
  \qquad
  (2\partial_{\bar z}+z)u_\tau=0
  \quad\text{on }\C,
\end{equation}
and the density $q_\tau=\abs{u_\tau}^2$ is periodic with respect to
the lattice generated by $a$ and $a\tau$. Every fundamental cell $C$
has area $\pi$, and
\begin{equation}\label{eq:btau}
  \frac{\pi\int_C\abs{u_\tau}^4}
       {\left(\int_C\abs{u_\tau}^2\right)^2}
  =\frac{\langle q_\tau^2\rangle}{\langle q_\tau\rangle^2}
  =b(\tau)
  :=\sum_{k,l\in\Z}
    e^{-\pi\abs{l-k\tau}^2/\tau_2},
\end{equation}
where $\langle h\rangle=\pi^{-1}\int_C h$.

In particular, for $\tau=\frac12+i\frac{\sqrt3}{2}$, the lattice
is triangular and
\begin{equation}\label{eq:periodic-triangular-ratio}
  \frac{\pi\int_C\abs{u_\tau}^4}
       {\left(\int_C\abs{u_\tau}^2\right)^2}
  =b_\triangle.
\end{equation}
Thus $b_\triangle$ is attained as the normalized quartic quotient
over a fundamental cell of a solution of the LLL equation with
triangular-periodic density.
\end{lemma}

\begin{proof}
Apply \cite[Proposition~2.2]{GST25} (see also \cite[Theorem~1.4]{ABN06}) with $N=1$ and $\gamma=a$.
Part~(i) gives the required Gaussian-weighted entire function,
and parts~(ii)--(iii) give
\begin{equation}
  \frac{\pi\int_C\abs{u_\tau}^4}
       {\left(\int_C\abs{u_\tau}^2\right)^2}
  =\sum_{k,l\in\Z}e^{-a^2\abs{l-k\tau}^2}
  =b(\tau),
\end{equation}
since $a^2=\pi/\tau_2$. The equality with the quotient of averages
follows from $\abs{C}=\pi$. The triangular case is precisely
the definition \eqref{eq:btri}.
\end{proof}

\subsubsection{Proof of Lemma~\ref{lem:cell}:
the intrinsic cell coefficient}

\begin{proof}[Proof of Lemma~\ref{lem:cell}]
We first obtain a uniform upper bound from periodic
profiles. We then prove $c_L>1$ for each fixed $L>0$
by holomorphic compactness. Finally, a tiling argument
will show that $c_L$ converges to its supremum.

For every admissible density $q$ on $Q_L$,
Cauchy--Schwarz gives
\begin{equation}
  \left(\int_{Q_L}q\,\mathrm dx\right)^2
  \leqslant L^2\int_{Q_L}q^2\,\mathrm dx.
\end{equation}
Hence $c_L\geqslant1$ for every $L>0$.

\smallskip
\noindent\emph{Step 1. The periodic upper bound.}
Fix $L>0$ and $\tau\in\C$ with $\Im\tau>0$.
Let $f_\tau$ and
$q_\tau=\abs{f_\tau}^2e^{-\abs{z}^2}$
be supplied by Lemma~\ref{lem:periodic}.
Set
$u=e^{-\abs{z}^2/2}f_\tau/\langle q_\tau\rangle^{1/2}$
and $q=\abs{u}^2$.
Then $\langle q\rangle=1$ and
$\langle q^2\rangle=b(\tau)$ by \eqref{eq:btau}.
For each $\bm a\in\R^2$,
Lemma~\ref{lem:translation} shows that
$u_{\bm a}(x)=e^{-i\bm a^\perp\cdot x}u(x+\bm a)$
satisfies the local LLL equation on $Q_L$ and
$\abs{u_{\bm a}(x)}^2=q(x+\bm a)$.
Its mass is positive because $f_\tau$ is
holomorphic and not identically zero.
Thus \eqref{eq:cellquotient} gives
\begin{equation}\label{eq:cell-periodic-local}
  \frac{c_L}{L^2}
  \left(\int_{\bm a+Q_L}q\,\mathrm dx\right)^2
  \leqslant
  \int_{\bm a+Q_L}q^2\,\mathrm dx,
  \qquad \bm a\in\R^2.
\end{equation}
Let $C$ be a fundamental cell of the periodicity
lattice. For $p\in\{1,2\}$, periodicity and
Fubini's theorem yield
\begin{equation}\label{eq:cell-periodic-average}
    \frac1{\abs{C}}\int_C
      \int_{\bm a+Q_L}q^p(x)\,\mathrm dx\,
      \mathrm d\bm a
    =\int_{Q_L}\frac1{\abs{C}}\int_C
      q^p(x+\bm a)\,\mathrm d\bm a\,\mathrm dx
    =L^2\langle q^p\rangle.
\end{equation}
In particular, Jensen's inequality gives
\begin{equation}
  \frac1{\abs{C}}\int_C
    \left(\int_{\bm a+Q_L}q\,\mathrm dx\right)^2
    \mathrm d\bm a
  \geqslant
  \left(
    \frac1{\abs{C}}\int_C
      \int_{\bm a+Q_L}q\,\mathrm dx\,\mathrm d\bm a
  \right)^2
  =L^4.
\end{equation}
Averaging \eqref{eq:cell-periodic-local} and using
\eqref{eq:cell-periodic-average}, we obtain
\begin{equation}
  c_LL^2
  \leqslant
  \frac{c_L}{L^2\abs{C}}\int_C
    \left(\int_{\bm a+Q_L}q\,\mathrm dx\right)^2
    \mathrm d\bm a
  \leqslant L^2b(\tau).
\end{equation}
Therefore $c_L\leqslant b(\tau)$.
Choosing $\tau=\frac12+i\frac{\sqrt3}{2}$
gives $c_L\leqslant b_\triangle$ for every $L>0$.

\smallskip
\noindent\emph{Step 2. The strict inequality on a fixed cell.}
Fix $L>0$ and suppose, for contradiction, that $c_L=1$.
By the definition of the infimum in \eqref{eq:cell},
for every $n\in\mathbb N$ we can choose
$f_n\in\mathcal O(Q_L)\cap L^2(Q_L;\C)$ such that,
with $q_n=\abs{f_n}^2e^{-\abs{z}^2}$, $ \norm{q_n}_{L^1(Q_L)}=L^2$ and $L^2\leqslant\norm{q_n}_{L^2(Q_L)}^2\leqslant L^2(1+1/n)$.
Expanding the square gives
\begin{equation}\label{eq:cell-density-constant-limit}
  \norm{q_n-1}_{L^2(Q_L)}^2
  =\int_{Q_L}q_n^2\,\mathrm dx-L^2
  \longrightarrow0.
\end{equation}
Since $\abs{z}^2\leqslant L^2/2$ on $Q_L$, we also have
\begin{equation}
  \int_{Q_L}\abs{f_n}^2\,\mathrm dx
  \leqslant e^{L^2/2}\int_{Q_L}q_n\,\mathrm dx
  =e^{L^2/2}L^2.
\end{equation}
For a nonempty compact set $K\subset Q_L$, set
$d:=\frac12\dist(K,\partial Q_L)>0$.
The submean inequality for $\abs{f_n}^2$ gives
\begin{equation}
  \sup_{z\in K}\abs{f_n(z)}^2
  \leqslant
  \frac1{\pi d^2}\int_{Q_L}\abs{f_n}^2\,\mathrm dx
  \leqslant\frac{e^{L^2/2}L^2}{\pi d^2};
\end{equation}
see \cite[Chapter~4]{Ahlfors1979}.
Thus $(f_n)$ is locally uniformly bounded.
By Montel's theorem, a subsequence converges
locally uniformly on $Q_L$ to some
$f\in\mathcal O(Q_L)$.
It follows that
$q_n\to\abs{f}^2e^{-\abs{z}^2}$
locally uniformly.
Together with \eqref{eq:cell-density-constant-limit},
this gives $\abs{f(z)}^2e^{-\abs{z}^2}=1$
for every $z\in Q_L$.
Hence $f$ has no zeros and $\log\abs{f}^2$
is harmonic. But $\log\abs{f(z)}^2=\abs{z}^2$
implies
$0=\Delta\log\abs{f}^2=\Delta\abs{z}^2=4$,
a contradiction.
Therefore $c_L>1$.

\smallskip
\noindent\emph{Step 3. The large-square limit.}
Fix $\ell>0$.
For each $L\geqslant\max\{1,\ell\}$, the definition
of the infimum in \eqref{eq:cell} allows us to choose
an admissible profile whose normalized energy
differs from $c_L$ by at most $L^{-1}$.
Thus there exists
$f_L\in\mathcal O(Q_L)\cap L^2(Q_L;\C)$ such that,
with $q_L=\abs{f_L}^2e^{-\abs{z}^2}$,
\begin{equation}\label{eq:cell-near-minimizer}
  \int_{Q_L}q_L\,\mathrm dx=L^2,
  \qquad
  \int_{Q_L}q_L^2\,\mathrm dx
  \leqslant(c_L+L^{-1})L^2
  \leqslant(b_\triangle+1)L^2.
\end{equation}
The last inequality uses $c_L\leqslant b_\triangle$
from Step~1 and $L\geqslant1$.

Set $n=\lfloor L/\ell\rfloor$.
Partition $Q_{n\ell}$, up to a set of measure zero,
into $n^2$ disjoint open squares $S_j$ of side
$\ell$, indexed by $j\in\{1,\ldots,n^2\}$.
The remainder
$B_L=Q_L\setminus\bigcup_{j=1}^{n^2}S_j$
satisfies
$\abs{B_L}=L^2-n^2\ell^2
=(L-n\ell)(L+n\ell)\leqslant2\ell L$.
By Cauchy--Schwarz and
\eqref{eq:cell-near-minimizer},
\begin{equation}\label{eq:remainder}
  \int_{B_L}q_L\,\mathrm dx
  \leqslant
    \abs{B_L}^{1/2}\norm{q_L}_{L^2(Q_L)}
  \leqslant C L^2\sqrt{\frac\ell L},
\end{equation}
where $C=\sqrt{2(b_\triangle+1)}$.
For $j\in\{1,\ldots,n^2\}$, write
$M_j=\int_{S_j}q_L\,\mathrm dx$.
The cell inequality \eqref{eq:cellquotient},
transported to $S_j$ by
Lemma~\ref{lem:translation}, gives
\begin{equation}
  \int_{S_j}q_L^2\,\mathrm dx
  \geqslant\frac{c_\ell}{\ell^2}M_j^2.
\end{equation}
Moreover,
$\sum_{j=1}^{n^2}M_j
=L^2-\int_{B_L}q_L\,\mathrm dx$.
For sufficiently large $L$, so that
$C\sqrt{\ell/L}<1$, we obtain
\begin{equation}
  \begin{aligned}
    (c_L+L^{-1})L^2
    &\overset{\eqref{eq:cell-near-minimizer}}{\geqslant}
      \int_{Q_L}q_L^2\,\mathrm dx
    \geqslant
      \frac{c_\ell}{\ell^2}
      \sum_{j=1}^{n^2}M_j^2
    \geqslant
      \frac{c_\ell}{n^2\ell^2}
      \left(\sum_{j=1}^{n^2}M_j\right)^2\\
    &\overset{\eqref{eq:remainder}}{\geqslant}
      \frac{c_\ell L^4}{n^2\ell^2}
      \left(1-C\sqrt{\frac\ell L}\right)^2
    \geqslant
      c_\ell L^2
      \left(1-C\sqrt{\frac\ell L}\right)^2.
  \end{aligned}
\end{equation}
The third inequality is Cauchy--Schwarz, and
the last uses $n^2\ell^2\leqslant L^2$.
Dividing by $L^2$ and letting $L\to\infty$ gives
$\liminf_{L\to\infty}c_L\geqslant c_\ell$.

Since $\ell>0$ was arbitrary,
$\liminf_{L\to\infty}c_L
\geqslant\sup_{\ell>0}c_\ell$.
The reverse inequality for the upper limit follows
from $c_L\leqslant\sup_{\ell>0}c_\ell$.
Thus the limit exists and
$\beta_*:=\lim_{L\to\infty}c_L
=\sup_{L>0}c_L\in(1,b_\triangle]$,
where the final bounds follow from Steps~1 and~2.
\end{proof}

\section{Compactness: proof of assertion (3)}\label{sec:compactness}
In this section we prove the compactness claim. The quartic term bounds the densities in $L^2$, while confinement prevents mass from escaping to infinity. Together they give assertion~(3) of Theorem~\ref{thm:gamma}. 

\begin{proof}[Proof of assertion (3) of Theorem~\ref{thm:gamma}]
Let $\eps_n\downarrow0$ and
$v_n\in L^2(\R^2;\C)$, $n\in\mathbb N$, satisfy the uniform bound
\begin{equation}
  M:=\sup_{n\in\mathbb N}\E_{\eps_n}(v_n)<\infty.
\end{equation}
Recall that $ \Sigma=\{v\in H^1(\R^2;\C):\abs{x}v\in L^2(\R^2;\C)\}$ and by definition (see \eqref{eq:modelintro}), $\E_\eps(v)=+\infty$
outside the unit sphere of $\Sigma$.
The energy bound therefore implies
$v_n\in\Sigma$ and
$\norm{v_n}_{L^2(\R^2)}=1$ for every
$n\in\mathbb N$.

Set $\rho_n=\abs{v_n}^2$.
Since $v_n\in H^1(\R^2;\C)\subset L^4(\R^2;\C)$,
each $\rho_n$ belongs to
$L^1(\R^2;\R)\cap L^2(\R^2;\R)$.
Moreover,
\begin{equation}
  \rho_n\geqslant0
  \quad\text{a.e. on }\R^2,
  \qquad
  \int_{\R^2}\rho_n\,\mathrm dx=1,
  \qquad n\in\mathbb N.
\end{equation}

\par\medskip\noindent\textit{Step 1. Uniform bounds and weak compactness.}
By \eqref{eq:factor}, $K_{\eps_n}(v_n)=\norm{A_{\eps_n}v_n}_{L^2(\mathbb R^2)}^2\geqslant0$; hence the decomposition \eqref{eq:modelintro} of the energy gives, for every $n$,
\begin{equation}\label{eq:compactbounds}
 \norm{\rho_n}_{L^2(\mathbb R^2)}^2
 +2\int_{\mathbb R^2}|x|^2\rho_n\,\mathrm dx
 \leqslant\E_{\eps_n}(v_n)\leqslant M.
\end{equation}
Weak compactness of bounded sets in $L^2(\mathbb R^2;\mathbb R)$ yields
indices $n_k\uparrow\infty$ and $\rho\in L^2(\mathbb R^2;\mathbb R)$ such that
\begin{equation}
 \rho_{n_k}\rightharpoonup\rho\quad\text{in }L^2(\mathbb R^2;\mathbb R).
\end{equation}
The weak limit $\rho$ is also nonnegative.
Indeed, for every $\varphi\in L^2(\R^2;\R)$
with $\varphi\geqslant0$ almost everywhere,
weak convergence and $\rho_{n_k}\geqslant0$ give
\begin{equation}
  \int_{\R^2}\rho\varphi\,\mathrm dx
  =\lim_{k\to\infty}
    \int_{\R^2}\rho_{n_k}\varphi\,\mathrm dx
  \geqslant0.
\end{equation}
We may take $\varphi=\rho_-:=\max\{-\rho,0\}$,
since $\rho_-\in L^2(\R^2;\R)$.
Using $\rho\rho_-=-(\rho_-)^2$ almost everywhere,
we obtain
\begin{equation}
  0\leqslant\int_{\R^2}\rho\rho_-\,\mathrm dx
  =-\norm{\rho_-}_{L^2(\R^2)}^2.
\end{equation}
Hence $\rho_-=0$ almost everywhere, or equivalently
$\rho\geqslant0$ almost everywhere on $\R^2$.

\smallskip
\noindent\emph{Step 2. Preservation of mass.}
We show that no mass is lost in the weak limit,
that is, $\norm{\rho}_{L^1(\R^2)}=1$.
We cannot test weak $L^2$ convergence against the
constant function $1$, since $1\notin L^2(\R^2)$, but
the workaround is a standard truncation argument:
we test on bounded balls and control the mass
outside them by the second-moment bound.
Precisely, for $R>0$, let $B_R=\{x\in\R^2:\abs{x}<R\}$.
Since $\rho_n\geqslant0$ and $\abs{x}^2\geqslant R^2$
on $\R^2\setminus B_R$, we have
\begin{equation}
  \int_{\R^2\setminus B_R}\rho_n\,\mathrm dx
  \leqslant\frac1{R^2}
    \int_{\R^2}\abs{x}^2\rho_n\,\mathrm dx
  \overset{\eqref{eq:compactbounds}}{\leqslant}
    \frac{M}{2R^2},
  \qquad n\in\mathbb N.
\end{equation}
Thus the mass outside $B_R$ tends to zero
uniformly in $n$ as $R\to\infty$.
Using $\norm{\rho_n}_{L^1(\R^2)}=1$ and
testing the weak convergence
$\rho_{n_k}\rightharpoonup\rho$ against
$\mathbf1_{B_R}\in L^2(\R^2;\R)$, we obtain
\begin{equation}
  1-\frac{M}{2R^2}
  \leqslant
  \lim_{k\to\infty}\int_{B_R}\rho_{n_k}\,\mathrm dx
  =\int_{B_R}\rho\,\mathrm dx
  \leqslant1,
  \qquad R>0.
\end{equation}
Since $\rho\geqslant0$, letting $R\to\infty$
and applying monotone convergence gives
$\rho\in L^1(\R^2;\R)$ and
$\norm{\rho}_{L^1(\R^2)}=1$, as required.

\smallskip
\noindent\emph{Step 3. Finiteness of the limiting second moment.}
We show that $\rho$ has finite second moment and,
more precisely, that this moment does not exceed
the lower limit of the moments of $\rho_{n_k}$.
Since $\abs{x}^2\notin L^2(\R^2)$, we again
truncate the test function to a bounded ball.
For each $R>0$, the function
$\abs{x}^2\mathbf1_{B_R}$ belongs to
$L^2(\R^2;\R)$. Weak convergence therefore gives
\begin{equation}
    \int_{B_R}\abs{x}^2\rho\,\mathrm dx
    =\lim_{k\to\infty}
      \int_{B_R}\abs{x}^2\rho_{n_k}\,\mathrm dx
    \leqslant\liminf_{k\to\infty}
      \int_{\R^2}\abs{x}^2\rho_{n_k}\,\mathrm dx
    \overset{\eqref{eq:compactbounds}}{\leqslant}
      \frac M2.
\end{equation}
Since $\rho\geqslant0$, monotone convergence
as $R\to\infty$ yields
\begin{equation}\label{eq:momentlsc}
  \int_{\R^2}\abs{x}^2\rho\,\mathrm dx
  \leqslant\liminf_{k\to\infty}
    \int_{\R^2}\abs{x}^2\rho_{n_k}\,\mathrm dx
  \leqslant\frac M2<\infty.
\end{equation}
Together with Steps~1 and~2, this proves
assertion~(3).
The arguments in Steps~2 and~3 apply to any
weakly convergent bounded-energy sequence;
they require no further subsequence extraction.
\end{proof}

\section{The \texorpdfstring{$\Gamma$}{Gamma}-liminf inequality:
proof of assertion (1)}
\label{sec:lower}

We prove that every family of wave functions whose
densities converge weakly to $\rho$ has limiting
energy at least $\F(\rho)$.
The main task is to obtain the interaction bound
with coefficient $\beta_*$.
We first prove this bound for LLL states by applying
the free-cell inequality on small squares and
averaging over their centers.
We then use projection to extend it to general
bounded-energy states.
Combining the interaction bound with lower
semicontinuity of the confinement term will prove
assertion~(1) in Subsection~\ref{sec:liminfproof}.

The first lemma concerns LLL states of arbitrary
mass. We need this formulation because the projection
$P_\eps v_\eps$ of a normalized wave function need
not itself have mass one.

\begin{lemma}[Averaged LLL inequality]
\label{lem:average}
For each $\eps\in(0,1)$, let
$u_\eps\in\mathcal L_\eps\cap L^4(\R^2;\C)$,
where $\mathcal L_\eps$ is the lowest Landau level
defined in \eqref{eq:lllintro}.
Thus $u_\eps\in L^2(\R^2;\C)$ and
$(2\partial_{\bar z}+z/\eps)u_\eps=0$
in $\mathcal D'(\R^2)$.

Assume that
$\sup_{\eps\in(0,1)}\norm{u_\eps}_{L^4(\R^2)}<\infty$
and that, for some $\rho\in L^2(\R^2;\R)$,
$\abs{u_\eps}^2\rightharpoonup\rho$
in $L^2(\R^2;\R)$ as $\eps\downarrow0$.
Then
\begin{equation}
  \liminf_{\eps\downarrow0}
    \int_{\R^2}\abs{u_\eps}^4\,\mathrm dx
  \geqslant
    \beta_*\int_{\R^2}\rho^2\,\mathrm dx.
\end{equation}
\end{lemma}

\begin{remark}[Why the coefficient $\beta_*$ appears]
Without the LLL constraint, weak lower semicontinuity
would give only the coefficient one:
$\liminf_{\eps\downarrow0}
\int_{\R^2}\abs{u_\eps}^4\,\mathrm dx
\geqslant\int_{\R^2}\rho^2\,\mathrm dx$.
The larger coefficient $\beta_*$ comes from applying
the holomorphic cell inequality locally.
Indeed, fix $L>0$ and consider
$S=\bm a+\sqrt\eps Q_L$, with $\bm a\in\R^2$.
Magnetic translation and rescaling transform the
restriction of $u_\eps$ to $S$ into a competitor
for the cell quotient defining $c_L$.
Consequently,
\begin{equation}
  \frac1{\abs{S}}\int_S\abs{u_\eps}^4\,\mathrm dx
  \geqslant
  c_L\left(
    \frac1{\abs{S}}\int_S\abs{u_\eps}^2\,\mathrm dx
  \right)^2.
\end{equation}
Thus $c_L$ compares the mean squared density with
the square of its mean. The inequality $c_L>1$
quantifies the local density variations required
by the holomorphic constraint.
For fixed $L$, the square $S$ shrinks as
$\eps\downarrow0$, while its side length in the
microscopic coordinate $y=(x-\bm a)/\sqrt\eps$
remains $L$.
Averaging the local inequality over $\bm a$
and passing to the weak limit gives the lemma's
lower bound with coefficient $c_L$.
Since this holds for every $L>0$, we may take
$\sup_{L>0}c_L=\beta_*$, by
Lemma~\ref{lem:cell}.
The coefficient $\beta_*$ therefore records the
optimal local quartic cost imposed by the LLL
constraint.
\end{remark}

\begin{proof}
We first average the cell inequality over all
square centers. We will then pass to the limit
$\eps\downarrow0$ at fixed cell size $L$,
and finally let $L\to\infty$ to recover $\beta_*$.

\smallskip
\smallskip
\noindent\emph{Step 1. Averaging the cell inequality.}
We apply the free-cell inequality on squares of
side $\sqrt\eps L$, then integrate over their
centers to obtain a global interaction bound.

Fix $L>0$ and $\eps\in(0,1)$, and set
$\sigma_\eps=\abs{u_\eps}^2$.
For $\sigma\in L^1_{\mathrm{loc}}(\R^2;\R)$,
define its spatial average by
\begin{equation}
  T_{\eps,L}\sigma(\bm a)
  :=\frac1{\eps L^2}
    \int_{\bm a+\sqrt\eps Q_L}
      \sigma(x)\,\mathrm dx,
  \qquad \bm a\in\R^2.
\end{equation}
The denominator $\eps L^2$ is the area of the
square over which we average.
To apply the free-cell inequality
\eqref{eq:cellquotient}, we transform the
restriction of $u_\eps$ to
$\bm a+\sqrt\eps Q_L$ into a function on $Q_L$.
For each $\bm a\in\R^2$, define
\begin{equation}
  h(y)=e^{-i\bm a^\perp\cdot y/\sqrt\eps}
       u_\eps(\bm a+\sqrt\eps\,y),
  \qquad y\in Q_L.
\end{equation}
The dilation changes the parameter in the LLL
equation from $\eps$ to $1$, while the phase
compensates for the translation by $\bm a$.
More precisely, Lemma~\ref{lem:translation} gives
\begin{equation}
  A_1h(y)
  =\sqrt\eps\,e^{-i\bm a^\perp\cdot y/\sqrt\eps}
    (A_\eps u_\eps)(\bm a+\sqrt\eps\,y)
  =0
  \qquad\text{in }\mathcal D'(Q_L).
\end{equation}
Hence \eqref{eq:local-lll-representation} gives
$h(y)=e^{-\abs{y}^2/2}f(y_1+iy_2)$
for some $f\in\mathcal O(Q_L)$.
Since the phase has modulus one, the substitution
$x=\bm a+\sqrt\eps\,y$ yields
\begin{equation}
  \int_{\bm a+\sqrt\eps Q_L}
    \abs{u_\eps(x)}^p\,\mathrm dx
  =\eps\int_{Q_L}\abs{h(y)}^p\,\mathrm dy,
  \qquad p\in\{2,4\}.
\end{equation}
In particular, $h\in L^2(Q_L;\C)\cap L^4(Q_L;\C)$.
If $h\not\equiv0$, the quotient
\eqref{eq:cellquotient} therefore gives
\begin{equation}
    \int_{\bm a+\sqrt\eps Q_L}
      \sigma_\eps^2\,\mathrm dx
    =\eps\int_{Q_L}\abs{h}^4\,\mathrm dy
    \overset{\eqref{eq:cellquotient}}{\geqslant}
      \frac{\eps c_L}{L^2}
      \left(\int_{Q_L}\abs{h}^2\,\mathrm dy\right)^2
    =\frac{c_L}{\eps L^2}
      \left(\int_{\bm a+\sqrt\eps Q_L}
        \sigma_\eps\,\mathrm dx\right)^2.
\end{equation}
If $h\equiv0$, the same inequality holds with
both sides zero.
Dividing by $\eps L^2$, we conclude that
$T_{\eps,L}(\sigma_\eps^2)(\bm a)
\geqslant c_L
\bigl(T_{\eps,L}\sigma_\eps(\bm a)\bigr)^2$
for every $\bm a\in\R^2$.
We now integrate over $\bm a$.
For fixed $x\in\R^2$, the set of centers
satisfying $x\in\bm a+\sqrt\eps Q_L$ is
$x-\sqrt\eps Q_L$, whose area is $\eps L^2$.
Tonelli's theorem thus gives
\begin{equation}
    \int_{\R^2}
      T_{\eps,L}(\sigma_\eps^2)(\bm a)\,
      \mathrm d\bm a
    =\frac1{\eps L^2}
      \int_{\R^2}\sigma_\eps^2(x)
      \left(\int_{x-\sqrt\eps Q_L}
        \mathrm d\bm a\right)\mathrm dx
    =\int_{\R^2}\sigma_\eps^2(x)\,\mathrm dx.
\end{equation}
Integrating the local inequality consequently yields
\begin{equation}\label{eq:averagedcellbound}
  \int_{\R^2}\sigma_\eps^2\,\mathrm dx
  \geqslant
  c_L\int_{\R^2}
    \bigl(T_{\eps,L}\sigma_\eps(\bm a)\bigr)^2
    \,\mathrm d\bm a
  =c_L\norm{T_{\eps,L}\sigma_\eps}_{L^2(\R^2)}^2.
\end{equation}

\smallskip
\noindent\emph{Step 2. The averages have the same weak limit.}
Fix $L>0$.
We show that
$T_{\eps,L}\sigma_\eps\rightharpoonup\rho$
in $L^2(\R^2;\R)$ as $\eps\downarrow0$.
The squares defining the averages shrink to a point.
To use this fact despite the dependence of
$\sigma_\eps$ on $\eps$, we transfer the averaging
operator to a fixed test function.
Since $Q_L=-Q_L$, the operator $T_{\eps,L}$
is convolution with the even kernel
$(\eps L^2)^{-1}\mathbf1_{\sqrt\eps Q_L}$,
whose integral is one.
It is therefore self-adjoint on $L^2(\R^2;\R)$
and, by Young's inequality, satisfies
$\norm{T_{\eps,L}\varphi}_{L^2(\R^2)}
\leqslant\norm{\varphi}_{L^2(\R^2)}$;
see \cite[Theorem~4.15]{Brezis2011}.
For each fixed $\varphi\in L^2(\R^2;\R)$,
the averaging formula gives
\begin{equation}
  \begin{aligned}
    \norm{T_{\eps,L}\varphi-\varphi}_{L^2(\R^2)}
    &\leqslant
      \frac1{\eps L^2}
      \int_{\sqrt\eps Q_L}
        \norm{\varphi(\,\cdot-\bm h)-\varphi}_{L^2(\R^2)}
        \,\mathrm d\bm h\\
    &\leqslant
      \sup_{\abs{\bm h}\leqslant\sqrt{\eps/2}L}
        \norm{\varphi(\,\cdot-\bm h)-\varphi}_{L^2(\R^2)}
      \longrightarrow0.
  \end{aligned}
\end{equation}
The convergence follows from continuity of translations
in $L^2$, since $\sqrt{\eps/2}L\to0$.

Set
$C=\sup_{\eps\in(0,1)}
\norm{\sigma_\eps}_{L^2(\R^2)}<\infty$.
By self-adjointness,
\begin{equation}
    \int_{\R^2}
      (T_{\eps,L}\sigma_\eps-\rho)\varphi\,\mathrm dx
    =
      \int_{\R^2}
        \sigma_\eps(T_{\eps,L}\varphi-\varphi)\,\mathrm dx +
      \int_{\R^2}(\sigma_\eps-\rho)\varphi\,\mathrm dx.
\end{equation}
Consequently,
\begin{equation}
    \left|\int_{\R^2}
      (T_{\eps,L}\sigma_\eps-\rho)\varphi\,\mathrm dx\right|
    \leqslant
      C\norm{T_{\eps,L}\varphi-\varphi}_{L^2(\R^2)} +
      \left|\int_{\R^2}
        (\sigma_\eps-\rho)\varphi\,\mathrm dx\right|
      \longrightarrow0.
\end{equation}
The first term tends to zero by the preceding
strong convergence of the averaged test function;
the second tends to zero because
$\sigma_\eps\rightharpoonup\rho$.
This proves
$T_{\eps,L}\sigma_\eps\rightharpoonup\rho$
in $L^2(\R^2;\R)$.

\smallskip
\noindent\emph{Step 3. Passage to the intrinsic coefficient.}
We now pass to the limit in the averaged cell
inequality. For every fixed $L>0$,
\eqref{eq:averagedcellbound} and the weak convergence
established in Step~2 give
\begin{equation}
    \liminf_{\eps\downarrow0}
      \int_{\R^2}\abs{u_\eps}^4\,\mathrm dx
    \overset{\eqref{eq:averagedcellbound}}{\geqslant}
      c_L\liminf_{\eps\downarrow0}
        \norm{T_{\eps,L}\sigma_\eps}_{L^2(\R^2)}^2
    \geqslant c_L\norm{\rho}_{L^2(\R^2)}^2,
\end{equation}
where the last inequality is weak lower
semicontinuity of the $L^2$ norm.
The left-hand side is independent of $L$,
so we may take the supremum over $L>0$:
\begin{equation}
    \liminf_{\eps\downarrow0}
      \int_{\R^2}\abs{u_\eps}^4\,\mathrm dx
    \geqslant
      \left(\sup_{L>0}c_L\right)
      \norm{\rho}_{L^2(\R^2)}^2
    \overset{\eqref{eq:intrinsiccell}}{=}
      \beta_*\int_{\R^2}\rho^2\,\mathrm dx.
\end{equation}
The same argument applies when $\eps$ ranges
over any prescribed sequence tending to zero.
\end{proof}

Lemma~\ref{lem:average} gives the interaction
lower bound
\begin{equation}
  \liminf_{\eps\downarrow0}
    \int_{\R^2}\abs{u_\eps}^4\,\mathrm dx
  \geqslant
    \beta_*\int_{\R^2}\rho^2\,\mathrm dx
\end{equation}
for families $u_\eps\in\mathcal L_\eps$
with uniformly bounded $L^4$ norms and
$\abs{u_\eps}^2\rightharpoonup\rho$
in $L^2(\R^2;\R)$.
We now extend this bound to general wave functions
with bounded energies $\E_\eps$, as defined
in \eqref{eq:modelintro}.

To do so, we use the projection estimate
\eqref{eq:r4} from Lemma~\ref{lem:projection}.
It shows that
$v_\eps-P_\eps v_\eps\to0$ in $L^4(\R^2;\C)$
for bounded-energy families.
We will deduce that projection preserves both
the weak limit of the densities and the limiting
quartic interaction.
Lemma~\ref{lem:average} will then apply to
$u_\eps=P_\eps v_\eps$.
Combining the resulting interaction bound with
the second-moment inequality \eqref{eq:momentlsc}
will give the total-energy lower bound stated
in Proposition~\ref{prop:lower}.

By the arguments in Section~\ref{sec:compactness},
any weak density limit of a normalized bounded-energy
family is nonnegative, has mass one, and has finite
second moment. We now prove the stronger interaction
bound and the resulting total-energy lower bound.

\begin{proposition}[Full GP lower bound]
\label{prop:lower}
Let $v_\eps\in\Sigma$ with
$\norm{v_\eps}_{L^2(\R^2)}=1$ for $\eps\in(0,1)$,
and let $\rho\in L^2(\R^2;\R)$.
Assume that
$\sup_{\eps\in(0,1)}\E_\eps(v_\eps)<\infty$
and
$\abs{v_\eps}^2\rightharpoonup\rho$
in $L^2(\R^2;\R)$ as $\eps\downarrow0$.
Then 
\begin{equation}\label{eq:interactionlower}
    \liminf_{\eps\downarrow0}
      \int_{\R^2}\abs{v_\eps}^4\,\mathrm dx
   \geqslant
      \beta_*\int_{\R^2}\rho^2\,\mathrm dx,\qquad
    \liminf_{\eps\downarrow0}\E_\eps(v_\eps)
    \geqslant\F(\rho).
\end{equation}
\end{proposition}

\begin{proof}
We first compare $v_\eps$ with its LLL projection.
We then apply Lemma~\ref{lem:average} to the
projected functions and combine the resulting
interaction bound with lower semicontinuity
of the confinement term.

\smallskip
\noindent\emph{Step 1. Projection preserves the limiting
density and quartic interaction.}
Set
$M=\sup_{\eps\in(0,1)}\E_\eps(v_\eps)<\infty$,
$u_\eps=P_\eps v_\eps$, and
$r_\eps=v_\eps-u_\eps$.
Since all terms in \eqref{eq:modelintro}
are nonnegative,
\begin{equation}
  K_\eps(v_\eps)\leqslant M,
  \qquad
  \norm{v_\eps}_{L^4(\R^2)}^4\leqslant M.
\end{equation}
Lemma~\ref{lem:projection} consequently gives
\begin{equation}
  \norm{r_\eps}_{L^4(\R^2)}^4
  \overset{\eqref{eq:r4}}{\leqslant}
  C\eps K_\eps(v_\eps)^2
  \leqslant C\eps M^2
  \longrightarrow0.
\end{equation}
The triangle inequality now implies
$\sup_{\eps\in(0,1)}
\norm{u_\eps}_{L^4(\R^2)}<\infty$.

To compare the densities, we use the pointwise
inequality
$\bigl|\abs{v_\eps}^2-\abs{u_\eps}^2\bigr|
\leqslant\abs{r_\eps}
(\abs{v_\eps}+\abs{u_\eps})$.
H\"older's inequality yields
\begin{equation}\label{eq:projection-density-error}
    \norm{\abs{v_\eps}^2-\abs{u_\eps}^2}_{L^2(\R^2)}
    \leqslant
      \norm{r_\eps}_{L^4(\R^2)}
      \left(
        \norm{v_\eps}_{L^4(\R^2)}
        +\norm{u_\eps}_{L^4(\R^2)}
      \right)
    \longrightarrow0.
\end{equation}
Hence $\abs{u_\eps}^2\rightharpoonup\rho$
in $L^2(\R^2;\R)$.
Moreover, factoring the difference of the squared
densities and applying Cauchy--Schwarz gives
\begin{equation}\label{eq:projection-interaction-error}
    \left|
      \int_{\R^2}\abs{v_\eps}^4\,\mathrm dx
      -\int_{\R^2}\abs{u_\eps}^4\,\mathrm dx
    \right|
    \leqslant
      \norm{\abs{v_\eps}^2-\abs{u_\eps}^2}_{L^2(\R^2)}
      \left(
        \norm{v_\eps}_{L^4(\R^2)}^2
        +\norm{u_\eps}_{L^4(\R^2)}^2
      \right)
      \longrightarrow0.
\end{equation}
Thus, projection preserves both the weak density
limit and the limiting quartic interaction.

\smallskip
\noindent\emph{Step 2. The interaction and total-energy
lower bounds.}
The functions $u_\eps\in\mathcal L_\eps$ satisfy
the hypotheses of Lemma~\ref{lem:average}.
That lemma does not require mass normalization,
so it applies even though
$\norm{u_\eps}_{L^2(\R^2)}$ need not equal one.
Together with
\eqref{eq:projection-interaction-error}, it gives
\begin{equation}\label{eq:quarticlower}
    \liminf_{\eps\downarrow0}
      \int_{\R^2}\abs{v_\eps}^4\,\mathrm dx
    =\liminf_{\eps\downarrow0}
      \int_{\R^2}\abs{u_\eps}^4\,\mathrm dx
    \overset{\text{Lemma~\ref{lem:average}}}{\geqslant}
      \beta_*\int_{\R^2}\rho^2\,\mathrm dx.
\end{equation}

We treat confinement using the original densities
$\abs{v_\eps}^2$.
The arguments in Steps~1--3 of
Section~\ref{sec:compactness} apply to any
weakly convergent bounded-energy sequence.
They show that $\rho\geqslant0$ almost everywhere,
$\int_{\R^2}\rho\,\mathrm dx=1$, and
\begin{equation}\label{eq:confinementlower}
  \int_{\R^2}\abs{x}^2\rho\,\mathrm dx
  \overset{\eqref{eq:momentlsc}}{\leqslant}
  \liminf_{\eps\downarrow0}
    \int_{\R^2}\abs{x}^2\abs{v_\eps}^2\,\mathrm dx
  \leqslant\frac M2.
\end{equation}
In particular, $\F(\rho)<\infty$.
No second-moment estimate for $u_\eps$ is needed.
Finally, we discard the nonnegative kinetic
excess in \eqref{eq:modelintro} and combine the
two lower bounds:
\begin{equation}
  \begin{aligned}
    \liminf_{\eps\downarrow0}\E_\eps(v_\eps)
    &\geqslant
      2\liminf_{\eps\downarrow0}
        \int_{\R^2}\abs{x}^2\abs{v_\eps}^2\,\mathrm dx
      +\liminf_{\eps\downarrow0}
        \int_{\R^2}\abs{v_\eps}^4\,\mathrm dx\\
    &\overset{\eqref{eq:confinementlower},
               \eqref{eq:quarticlower}}{\geqslant}
      2\int_{\R^2}\abs{x}^2\rho\,\mathrm dx
      +\beta_*\int_{\R^2}\rho^2\,\mathrm dx
    \overset{\eqref{eq:limit}}{=}\F(\rho).
  \end{aligned}
\end{equation}
Every argument remains valid when $\eps$ ranges
over a prescribed sequence tending to zero.
\end{proof}

\subsection{Proof of assertion (1) of
Theorem~\ref{thm:gamma}}
\label{sec:liminfproof}

Proposition~\ref{prop:lower} assumes a uniform
energy bound. To prove the lower bound for an
arbitrary convergent sequence, we select a
subsequence realizing its lower limiting energy.

\begin{proof}
Let $\eps_n\downarrow0$,
$v_n\in L^2(\R^2;\C)\cap L^4(\R^2;\C)$,
$n\in\mathbb N$, and suppose that
$\abs{v_n}^2\rightharpoonup\rho$
in $L^2(\R^2;\R)$.
Set
$\lambda=\liminf_{n\to\infty}\E_{\eps_n}(v_n)$.
If $\lambda=+\infty$, the required inequality
is immediate. Suppose that $\lambda<\infty$.
We can choose indices $n_k\uparrow\infty$ such
that $\E_{\eps_{n_k}}(v_{n_k})\to\lambda$.
After discarding finitely many terms, we have
\begin{equation}
  0\leqslant\E_{\eps_{n_k}}(v_{n_k})
  \leqslant\lambda+1,
  \qquad k\in\mathbb N.
\end{equation}
By \eqref{eq:modelintro}, finiteness of these
energies implies $v_{n_k}\in\Sigma$ and
$\norm{v_{n_k}}_{L^2(\R^2)}=1$.
The selected densities still converge weakly
to $\rho$.
Proposition~\ref{prop:lower}, applied along
$\eps_{n_k}\downarrow0$, therefore gives
\begin{equation}
  \liminf_{n\to\infty}\E_{\eps_n}(v_n)
  =\lim_{k\to\infty}\E_{\eps_{n_k}}(v_{n_k})
  \geqslant\F(\rho).
\end{equation}
This proves \eqref{eq:directliminf}.
In particular, if $\F(\rho)=+\infty$, the lower
limiting energy must also be $+\infty$.
\end{proof}

\section{The \texorpdfstring{$\Gamma$}{Gamma}-limsup inequality: proof of assertion (2)}\label{sec:recovery}
In this section we construct a recovery sequence for the $\Gamma-\limsup$ inequality. The construction has three parts. Lemma~\ref{lem:cutoff} localizes nearly optimal cell states; Proposition~\ref{prop:steprecover} produces LLL recovery sequences for densities constant on finitely many squares; Subsection~\ref{sec:recoveryproof} proves assertion~(2) for an arbitrary finite-energy density by approximation and an explicit diagonal selection. 

\begin{figure}[tp]
\centering
\includegraphics{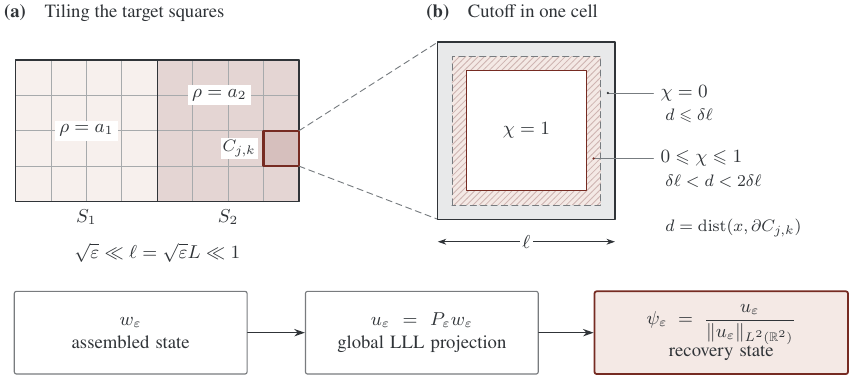}
\caption{Recovery construction. (a) Each fixed target square $S_j$ is filled by cells $C_{j,k}$ of side $\ell=\sqrt\eps L$, with $\sqrt\eps\ll\ell\ll1$. (b) One enlarged cell: $d$ denotes distance to its boundary, the cutoff is zero for $d\leqslant\delta\ell$, and it equals one for $d\geqslant2\delta\ell$. The intervening band permits the smooth transition. The local profiles are normalized and assembled before the global magnetic projection; the final normalization produces the recovery state. The diagram prescribes no vortex arrangement. Projection generally creates tails outside the original cell supports.}
\label{fig:recovery}
\end{figure}

The construction places rescaled cell profiles inside the target squares. Cutting them off destroys the LLL constraint, so we apply the global magnetic projection after assembly. The scales are chosen so that projection changes the density, interaction energy and confinement energy by quantities tending to zero. Figure~\ref{fig:recovery} summarizes the construction.

\subsection{Cutting off a nearly optimal microscopic profile}
The next lemma establishes mass, quartic-energy and magnetic-excess estimates for compactly supported cell profiles. We use them in Lemma~\ref{lem:singlesquare} to construct each local profile, and in Step~2 of Proposition~\ref{prop:steprecover} to control the projection error.
\begin{lemma}[Cutoff cell profiles]\label{lem:cutoff}
Let $L\geqslant1$ and $0<\delta<1/8$. There are smooth compactly supported functions $\tilde h_{L,\delta}$ in $Q_L$ such that
\begin{align}
 m_{L,\delta}:=L^{-2}\int_{\mathbb R^2}|\tilde h_{L,\delta}|^2
 &\in[1-C\sqrt\delta,1],\label{eq:cutmass}\\
 L^{-2}\int_{\mathbb R^2}|\tilde h_{L,\delta}|^4&\leqslant c_L+L^{-1},\label{eq:cutquartic}\\
 \norm{A_1\tilde h_{L,\delta}}_{L^2(\mathbb R^2)}^2&\leqslant C\delta^{-2}.\label{eq:cutA}
\end{align}
The constants are independent of $L$ and $\delta$. The function is extended by zero outside $Q_L$ in these formulas.
\end{lemma}
\begin{proof}
\par\medskip\noindent\textit{Step 1. A cutoff retains almost all the mass and does not increase the quartic cost.}
In this step we prove \eqref{eq:cutmass} and \eqref{eq:cutquartic}.

Since $c_L\leqslant b_\triangle<\infty$ by Lemma~\ref{lem:cell}, the definition \eqref{eq:cell} provides a holomorphic $f_L$ with $\int_{Q_L}|f_L|^2e^{-|z|^2}=L^2$ and quartic integral at most $(c_L+L^{-1})L^2\leqslant(b_\triangle+1)L^2$. Put $h_L=f_Le^{-|z|^2/2}$, so $A_1h_L=0$ in $Q_L$.

Choose a real cutoff $\chi$ with $0\leqslant\chi\leqslant1$, $\chi=1$ when $\dist(x,\partial Q_L)\geqslant2\delta L$, $\chi=0$ when the distance is at most $\delta L$, and $|\nabla\chi|\leqslant C/(\delta L)$. Such a cutoff is obtained as a product of one-dimensional smooth cutoffs. Its transition and discarded region has area at most $C\delta L^2$. With $\tilde h=\chi h_L$, Cauchy--Schwarz and this quartic bound give
\begin{equation}
 \int_{Q_L}(1-\chi^2)|h_L|^2
 \leqslant(C\delta L^2)^{1/2}\left(\int_{Q_L}|h_L|^4\right)^{1/2}
 \leqslant C\sqrt\delta L^2.
\end{equation}
This proves \eqref{eq:cutmass}; \eqref{eq:cutquartic} follows from $|\chi|\leqslant1$. 
\par\medskip\noindent\textit{Step 2. Only cutoff derivatives contribute to the magnetic excess.}
We prove the last estimate without requiring boundary regularity of $f_L$. Because $h_L$ is smooth
 in the interior and $\chi$ vanishes in a neighborhood of the boundary, $\tilde h$ extends to a smooth compactly supported function on the plane. The product rule gives $A_1\tilde h=2(\partial_{\bar z}\chi)h_L$, and hence
\begin{equation}
 \norm{A_1\tilde h}_{L^2(\mathbb R^2)}^2\leqslant\frac{C}{\delta^2L^2}\int_{Q_L}|h_L|^2\leqslant C\delta^{-2}.
\end{equation}
No derivatives of $h_L$ need to be estimated near the boundary, and no cell minimizer or boundary regularity theorem is used.
\end{proof}

\subsection{Piecewise constant target densities}
The next proposition establishes LLL recovery for unit-mass densities constant on finitely many squares. It is required in Step~2 of the proof of assertion~(2), in Subsection~\ref{sec:recoveryproof}, to construct a recovery family for each density approximation before diagonal selection. We first prove the following auxiliary lemma.
\begin{lemma}[Pre-projected profile on a single square]
  \label{lem:singlesquare}Let $S \subset \R^2$ be a fixed bounded square with
  sides parallel to the axes, and set
  \begin{equation}
    \label{eq:scales} L_{\eps} := \eps^{- 3 / 8}, \qquad
    \ell_{\eps} := \sqrt{\eps} L_{\eps} = \eps^{1 / 8},
    \qquad \delta_{\eps} := \eps^{1 / 16} .
  \end{equation}
  For all sufficiently small $\eps > 0$, there exists $w_{S, \eps} \in
  C_c^{\infty}  (\mathrm{int} (S) ; \C)$ such that, as $\eps \to 0$:
  \begin{enumerate}
    \item $|w_{S, \eps} |^2 \rightharpoonup \textbf{1}_S$ in $L^2 (\R^2 ; \R)$
    and $\left\| w_{S, \eps} \right\|_{L^2 (\R^2)}^2 \to |S|$;

    \item $\limsup_{\eps \to 0}  \int_{\R^2} |w_{S, \eps} |^4 \hspace{0.17em}
    \mathrm{d} x \leqslant \beta_{\ast} |S|$;

    \item $\left\| A_{\eps} w_{S, \eps} \right\|_{L^2 (\R^2)}^2 \leqslant C
    |S|  \eps^{- 3 / 8}$ for a constant $C > 0$ independent of $S$ and $\eps$.
  \end{enumerate}
\end{lemma}

\begin{proof}
  For every sufficiently small $\eps > 0$, let $\tilde{h}_{\eps} = \tilde{h}_{L_{\eps},
  \delta_{\eps}}$ be the function compactly supported on $Q_{L_{\eps}}$ given
  by Lemma~\ref{lem:cutoff}, and let $m_{\eps} = m_{L_{\eps}, \delta_{\eps}}$
  be its associated mass. By Lemma~\ref{lem:cutoff}, $m_{\eps} \to 1$ as
  $\eps \to 0$; hence, for small $\eps$ we have $m_{\eps} \geqslant 1 / 2$.

  Set $s := \mathrm{side} (S)$ and $n_{\eps} := \lfloor s /
  \ell_{\eps} \rfloor$, and take $\eps$ so small that $\ell_{\eps} < s$. Place
  inside $S$ an $n_{\eps} \times n_{\eps}$ grid of pairwise disjoint open
  squares $C_{\eps, k} \subset S$ ($1 \leqslant k \leqslant N_{\eps} :=
  n_{\eps}^2$) of side $\ell_{\eps}$ and centers $\mathbf{a}_k \in \R^2$:
  $C_{\eps, k} := \mathbf{a}_k + Q_{\sqrt{\eps} L_{\eps}}$. With
  $S_{\eps} := \bigcup_{k = 1}^{N_{\eps}} C_{\eps, k}$ we denote the region
  occupied by the grid. Since $n_{\eps} \ell_{\eps} \geqslant s - \ell_{\eps} >
  0$, its measure satisfies $|S_{\eps}| = N_{\eps} \ell_{\eps}^2 \geqslant (s -
  \ell_{\eps})^2$, and therefore
  \begin{equation}
    \label{eq:strip} |S \setminus S_{\eps} | = s^2 - N_{\eps}
    \ell_{\eps}^2 \leqslant s^2 - (s - \ell_{\eps})^2 = 2 s \ell_{\eps} -
    \ell_{\eps}^2 \leqslant 2 s \ell_{\eps} \xrightarrow{\eps \to 0} 0.
  \end{equation}
  Define
  \begin{equation}
    w_{S, \eps} := \sum_{k = 1}^{N_{\eps}} w_{\eps, k}, \qquad
    w_{\eps, k} (x) := \frac{1}{\sqrt{m_{\eps}}}  \hspace{0.17em}
    e^{i\mathbf{a}_k^{\perp} \cdot x / \eps}  \hspace{0.17em}
    \tilde{h}_{\eps}  \hspace{-0.17em} \left( \frac{x
    -\mathbf{a}_k}{\sqrt{\eps}} \right) .
  \end{equation}
  By Lemma~\ref{lem:cutoff}, each $w_{\eps, k}$ is smooth and supported in the
  interior of $C_{\eps, k}$, i.e., $w_{\eps, k} \in C_c^{\infty} (\mathrm{int}
  (C_{\eps, k}) ; \C)$. Since the squares $C_{\eps, k}$ have pairwise disjoint
  interiors, the $w_{\eps, k}$ have pairwise disjoint supports, and therefore
  $w_{S, \eps} \in C_c^{\infty}  (\mathrm{int} (S) ; \C)$.

  Moreover, $A_{\eps}$ is a differential operator, hence local: if $u$ vanishes
  on an open set $U$, then so do all its derivatives, and thus $A_{\eps} u = 0$
  on $U$; that is, $\mathrm{supp} (A_{\eps} u) \subset \mathrm{supp} (u)$.
  Applied to $u = w_{\eps, k}$, this gives $\mathrm{supp} (A_{\eps} w_{\eps,
  k}) \subset \mathrm{int} (C_{\eps, k})$, so the functions $A_{\eps} w_{\eps,
  k}$ are smooth with pairwise disjoint supports as well. Consequently, at
  every $x \in \R^2$ at most one term of $A_{\eps} w_{S, \eps} = \sum_k
  A_{\eps} w_{\eps, k}$ is non-zero, the mixed terms in $|A_{\eps} w_{S,
  \eps}|^2$ vanish, and
  \begin{equation}
    \label{eq:Aadd} \left\| A_{\eps} w_{S, \eps} \right\|_{L^2 (\R^2)}^2 =
    \sum_{k = 1}^{N_{\eps}} \left\| A_{\eps} w_{\eps, k} \right\|_{L^2
    (\R^2)}^2 .
  \end{equation}
  The change of variables $x =\mathbf{a}_k + \sqrt{\eps} y$,
  Lemma~\ref{lem:translation} and \eqref{eq:cutmass}--\eqref{eq:cutA} give
  \begin{equation}
    \| w_{\eps, k} \|_{L^2 (\R^2)}^2 = \frac{\eps}{m_{\eps}}
    \left\| \tilde{h}_{\eps} \right\|_{L^2 (\R^2)}^2 = \frac{\eps}{m_{\eps}}
    m_{\eps} L_{\eps}^2 = \ell_{\eps}^2 = |C_{\eps, k} |
  \end{equation}
  and, summing over the $N_{\eps}$ cells, whose supports are disjoint, and using
  \eqref{eq:Aadd} for the last line (recall that $L_{\eps} \to \infty$ as
  $\eps \to 0$), we get
  \begin{align}
    \left\| w_{S, \eps} \right\|_{L^2 (\R^2)}^2 & = N_{\eps} \ell_{\eps}^2
    = |S_{\eps} | \xrightarrow{\eps \to 0} |S|,
    \nonumber\\
    \int_{\R^2} |w_{S, \eps} |^4 \hspace{0.17em} \mathrm{d} x & =
    N_{\eps}  \frac{\eps}{m^2_{\eps}}
    \int_{Q_{L_{\eps}}} | \tilde{h}_{\eps} |^4 \hspace{0.17em}
    \mathrm{d} x \leqslant \frac{c_{L_{\eps}} + L_{\eps}^{-
    1}}{m^2_{\eps}} |S_{\eps} | \xrightarrow{\eps \to 0}
    \beta_{\ast} |S|,  \label{eq:w_quartic}\\
    \left\| A_{\eps} w_{S, \eps} \right\|_{L^2 (\R^2)}^2 & =
    \frac{N_{\eps}}{m_{\eps}} \| A_1  \tilde{h}_{\eps}
    \|_{L^2 (\R^2)}^2 \leqslant \frac{2
    C_0 N_{\eps}}{\delta^2_{\eps}} = \frac{2 C_0 |S_{\eps}
    |}{\delta^2_{\eps} \ell^2_{\eps}} \leqslant 2 C_0 |S|  \eps^{-
    3 / 8}, \nonumber
  \end{align}
  where $C_0$ is the constant in \eqref{eq:cutA}, and we used $\delta_{\eps}^2
  \ell_{\eps}^2 = \eps^{3 / 8}$ and $|S_{\eps}| \leqslant |S|$. This proves~(2),
  (3) with $C := 2 C_0$, and the second half of~(1).

  It remains to show $|w_{S, \eps} |^2 \rightharpoonup \textbf{1}_S$. By
  \eqref{eq:w_quartic}, $|w_{S, \eps} |^2$ is bounded in $L^2 (\R^2 ; \R)$, so
  by density of $C_c^1 (\R^2 ; \R)$ in $L^2 (\R^2 ; \R)$ it suffices to test
  against $\varphi \in C_c^1  (\R^2 ; \R)$. Since $|w_{S, \eps} | = |w_{\eps,
  k} |$ on $C_{\eps, k}$ and $\mathrm{supp}\, w_{\eps, k} \subset C_{\eps, k}$,
  we have
  \begin{equation}
    \int_{C_{\eps, k}} (|w_{S, \eps} |^2 - 1) \hspace{0.17em}
    \mathrm{d} x = \left( \int_{\R^2} |w_{\eps,
    k} |^2 \hspace{0.17em} \mathrm{d} x \right) - | C_{\eps, k} | = 0,
    \label{eq:tempsec6withVAl1}
  \end{equation}
  and $|x -\mathbf{a}_k | \leqslant \ell_{\eps}$ on $C_{\eps, k}$.
  Recalling that $w_{S, \eps}$ is supported in $S_{\eps}$, that is, it vanishes
  on $S \setminus S_{\eps}$, for every $\varphi \in C_c^1 (\R^2 ; \R)$ we get
  \begin{eqnarray}
    \int_{\R^2} \varphi \left( |w_{S, \eps} |^2 - \textbf{1}_S \right)
    \hspace{0.17em} \mathrm{d} x & = & \int_S \varphi \left( |w_{S, \eps} |^2
    - 1 \right) \hspace{0.17em} \mathrm{d} x \nonumber\\
    & = & \int_{S_{\eps}} \varphi \left( |w_{S, \eps} |^2 - 1 \right)
    \hspace{0.17em} \mathrm{d} x + \int_{S \setminus S_{\eps}} \varphi \left(
    |w_{S, \eps} |^2 - 1 \right)  \mathrm{d} x \nonumber\\
    & = & \sum_{k = 1}^{N_{\eps}} \int_{C_{\eps, k}} \varphi
    \left( |w_{S, \eps} |^2 - 1 \right) \mathrm{d} x - \int_{S \setminus S_{\eps}}
    \varphi \hspace{0.17em} \mathrm{d} x
  \end{eqnarray}
  and therefore, by \eqref{eq:tempsec6withVAl1},
  \begin{eqnarray}
    \int_{\R^2} \varphi \left( |w_{S, \eps} |^2 - \textbf{1}_S \right)
    \hspace{0.17em} \mathrm{d} x & = & \sum_{k = 1}^{N_{\eps}}
    \int_{C_{\eps, k}} (\varphi (x) - \varphi (\mathbf{a}_k))  \left(
    |w_{S, \eps} (x) |^2 - 1 \right) \mathrm{d} x \nonumber\\
    &  & \qquad \qquad \qquad \qquad \qquad \qquad \qquad - \int_{S \setminus
    S_{\eps}} \varphi \hspace{0.17em} \mathrm{d} x.
  \end{eqnarray}
  Since $\int_{S_{\eps}} (|w_{S, \eps} |^2 + 1) \hspace{0.17em} \mathrm{d} x =
  2 |S_{\eps}| \leqslant 2 |S|$, using \eqref{eq:strip} we conclude
  \begin{equation} \left| \int_{\R^2} \varphi \left( |w_{S, \eps} |^2 - \textbf{1}_S \right)
     \hspace{0.17em} \mathrm{d} x \right| \leqslant \ell_{\eps}  \|
     \nabla \varphi \|_{L^{\infty} (\R^2)} \cdot 2| S | + 2 s
     \ell_{\eps} \| \varphi \|_{L^{\infty} (\R^2)}
     \xrightarrow{\eps \to 0} 0. \end{equation}
  This completes~(1).
\end{proof}
We are now ready to construct a recovery sequence for elementary functions.
\begin{proposition}[Recovery on finitely many squares]
  \label{prop:steprecover}Suppose $\rho = \sum_{j = 1}^J a_j
  \textbf{1}_{S_j}$, where $a_j > 0$, the $S_j \subset \R^2$ are bounded
  axis-parallel squares with disjoint interiors, and $\sum_{j = 1}^J a_j |S_j
  | = 1$. Then there exist normalized $\psi_{\eps} \in \LLL \cap \Sigma$ such
  that
  \begin{equation}
    \psi_{\eps} \xrightarrow{\mathrm{dens}} \rho, \qquad \limsup_{\eps \to 0}
    \E_{\eps} (\psi_{\eps}) \leqslant \F (\rho) .
  \end{equation}
\end{proposition}

\begin{proof}
  We glue the profiles of Lemma~\ref{lem:singlesquare}, project onto $\LLL$
  and normalize. In Step 1 we collect the limits of the
  glued function $w_{\eps}$; in Step 2 we show that the
  projection changes none of them; in Step 3 we conclude.

  \smallskip\noindent\emph{Step 1.} Let $w_{S_j, \eps}$ be given by
  Lemma~\ref{lem:singlesquare}, fix $R > 0$ with $\bigcup_{j = 1}^J S_j
  \subset B_R$, and set $w_{\eps} := \sum_{j = 1}^J \sqrt{a_j}
  \hspace{0.17em} w_{S_j, \eps} \in C_c^{\infty}  (B_R ; \C)$. Since the $S_j$
  have disjoint interiors, the summands have pairwise disjoint supports, and
  so do their images under the differential operator $A_{\eps}$ (which is
  local, see the proof of Lemma~\ref{lem:singlesquare}). Hence $|w_{\eps} |^2 =
  \sum_j a_j |w_{S_j, \eps} |^2$, $\int |w_{\eps} |^4 = \sum_j a_j^2  \int
  |w_{S_j, \eps} |^4$ and $\norm{A_{\eps} w_{\eps}}^2 = \sum_j a_j
  \norm{A_{\eps} w_{S_j, \eps}}^2$. Since $\rho^2 = \sum_j a_j^2
  \textbf{1}_{S_j}$ a.e., $\sum_j a_j |S_j | = 1$, and the $\limsup$ of a
  finite sum is at most the sum of the $\limsup$'s, Lemma~\ref{lem:singlesquare}
  gives
  \begin{equation}
    \label{eq:w_props} |w_{\eps} |^2 \rightharpoonup \rho \quad \text{in } L^2
    (\R^2 ; \R), \qquad \limsup_{\eps \to 0}  \int_{\R^2} |w_{\eps} |^4
    \hspace{0.17em} \mathrm{d} x \leqslant \beta_{\ast}  \int_{\R^2} \rho^2
    \hspace{0.17em} \mathrm{d} x, \qquad K_{\eps} (w_{\eps}) \leqslant C
    \eps^{- 3 / 8},
  \end{equation}
  where the last bound follows from $\sum_j a_j C |S_j| \eps^{-3/8} = C
  \eps^{-3/8}$, with $C$ the constant of Lemma~\ref{lem:singlesquare} statement (3).
  Likewise, $\norm{w_{\eps}}_{L^2}^2 = \sum_j a_j \norm{w_{S_j,
  \eps}}_{L^2}^2 \to \sum_j a_j |S_j | = 1$. Testing $|w_{\eps} |^2
  \rightharpoonup \rho$ against $|x|^2  \textbf{1}_{B_R} \in L^2 (\R^2 ; \R)$,
  and recalling $\mathrm{supp} \hspace{0.17em} w_{\eps} \subset B_R$, we also
  get
  \begin{equation}
    \label{eq:w_moments} \norm{w_{\eps}}_{L^2 (\R^2)}^2 \longrightarrow 1,
    \qquad \norm{|x| w_{\eps}}_{L^2 (\R^2)}^2 \longrightarrow \int_{\R^2}
    |x|^2 \rho \hspace{0.17em} \mathrm{d} x.
  \end{equation}
  \smallskip\noindent\emph{Step 2.} Put $u_{\eps} := P_{\eps} w_{\eps}$ and
  $r_{\eps} := w_{\eps} - u_{\eps}$. By \eqref{eq:w_moments},
  \eqref{eq:w_props}, Lemmas~\ref{lem:projection} and~\ref{lem:kernel}
  \begin{equation}
    \label{eq:r_vanish} \norm{r_{\eps}}_{L^2 (\R^2)}^2 \leqslant C \eps^{5 /
    8}, \qquad \norm{r_{\eps}}_{L^4 (\R^2)}^4 \leqslant C \eps^{1 / 4}, \qquad
    \norm{|x| r_{\eps}}_{L^2 (\R^2)} \leqslant C_R  \left( \eps^{1 / 2} +
    \eps^{5 / 16} \right),
  \end{equation}
  and all three right-hand sides tend to $0$. By the triangle inequality,
  $\norm{u_{\eps}}_{L^2}$, $\norm{|x| u_{\eps}}_{L^2}$ and
  $\norm{u_{\eps}}_{L^4}$ differ from the corresponding norms of $w_{\eps}$ by
  at most $\norm{r_{\eps}}_{L^2}$, $\norm{|x| r_{\eps}}_{L^2}$ and
  $\norm{r_{\eps}}_{L^4}$; in particular $|x| u_{\eps} \in L^2 (\R^2 ; \C)$,
  and $\norm{u_{\eps}}_{L^4}$ is bounded because $\norm{w_{\eps}}_{L^4}$ is bounded by
  \eqref{eq:w_props}. Moreover, since $\left| |u_{\eps} | - |w_{\eps} |
  \right| \leqslant |r_{\eps} |$, H{\"o}lder's inequality gives
  $\norm{|u_{\eps} |^2 - |w_{\eps} |^2}_{L^2} \leqslant \norm{r_{\eps}}_{L^4}
  \left( \norm{u_{\eps}}_{L^4} + \norm{w_{\eps}}_{L^4} \right) \to 0$.
  Together with \eqref{eq:w_props} and \eqref{eq:w_moments}, this yields
  \begin{equation}
    \label{eq:u_limits} |u_{\eps} |^2 \rightharpoonup \rho \text{ in } L^2
    (\R^2 ; \R), \quad \norm{u_{\eps}}_{L^2 (\R^2)} \to 1, \quad \norm{|x|
    u_{\eps}}_{L^2 (\R^2)}^2 \to \int_{\R^2} |x|^2 \rho \hspace{0.17em}
    \mathrm{d} x, 
  \end{equation}
  and 
  \begin{align}
      \limsup_{\eps \to 0}  \int_{\R^2} |u_{\eps} |^4
    \hspace{0.17em} \mathrm{d} x \leqslant \beta_{\ast}  \int_{\R^2} \rho^2
    \hspace{0.17em} \mathrm{d} x.
  \end{align}
  \smallskip\noindent\emph{Step 3.} By \eqref{eq:u_limits}, $\norm{u_{\eps}}_{L^2}> 1/2$ for
  small $\eps$, so we set $\psi_{\eps} :=
  \norm{u_{\eps}}_{L^2}^{- 1} u_{\eps}$.
  Since $A_{\eps} u_{\eps} = 0$, we have $u_{\eps} \in \mathcal H_{\eps}$,
  and $|x| u_{\eps} \in L^2$ by Step~2 and hence $u_{\eps} \in \Sigma$ by
  \eqref{eq:energy-domain}. Thus $\psi_{\eps} \in \LLL \cap \Sigma$,
  $\norm{\psi_{\eps}}_{L^2} = 1$ and $K_{\eps} (\psi_{\eps}) = 0$. Dividing
  \eqref{eq:u_limits} by powers of $\norm{u_{\eps}}_{L^2} \to 1$,
  \begin{equation}
    | \psi_{\eps} |^2 \rightharpoonup \rho \text{ in } L^2 (\R^2 ; \R), \quad
    \int_{\R^2} |x|^2 | \psi_{\eps} |^2 \hspace{0.17em} \mathrm{d} x \to
    \int_{\R^2} |x|^2 \rho \hspace{0.17em} \mathrm{d} x, \quad \limsup_{\eps
    \to 0}  \int_{\R^2} | \psi_{\eps} |^4 \hspace{0.17em} \mathrm{d} x
    \leqslant \beta_{\ast}  \int_{\R^2} \rho^2 \hspace{0.17em} \mathrm{d} x.
  \end{equation}
  The first relation is the statement $\psi_{\eps} \xrightarrow{\mathrm{dens}}
  \rho$. Since $\norm{\psi_{\eps}}_{L^2} = 1$ and $K_{\eps}(\psi_{\eps}) = 0$,
  \eqref{eq:modelintro} gives $\E_{\eps}(\psi_{\eps}) = 2 \int |x|^2
  |\psi_{\eps}|^2 + \int |\psi_{\eps}|^4$, and therefore
  \begin{equation} \limsup_{\eps \to 0} \E_{\eps}(\psi_{\eps}) \leqslant
     2 \int_{\R^2} |x|^2 \rho \, \mathrm{d} x
     + \beta_{\ast} \int_{\R^2} \rho^2 \, \mathrm{d} x = \F(\rho). \end{equation}
\end{proof}

\begin{remark}[The scale conditions]
The particular powers in \eqref{eq:scales} are convenient rather than optimal. The argument needs $L\to\infty$, $\ell=\sqrt\eps L\to0$, $\delta\to0$, and $\eps/(\delta^4\ell^4)\to0$ to make the fourth power of the $L^4$ projection error vanish. The last condition is stronger than the $L^2$ requirement. The displayed choices satisfy all of them. Only the qualitative cell limit \eqref{eq:intrinsiccell}, proved in Subsection~\ref{sec:cells}, is used here. A non-sharp rate for the cell constants will be obtained separately in Section~\ref{sec:certificate}. 
\end{remark}

\subsection{Proof of assertion (2) of Theorem~\ref{thm:gamma}}\label{sec:recoveryproof}
The proof extends Proposition~\ref{prop:steprecover} to every density in the finite domain of $\F$. 

\begin{proof}
Let $\rho\in L^2(\mathbb R^2;\mathbb R)$ satisfy $\F(\rho)<\infty$. Thus $\rho\geqslant0$, $\int_{\mathbb R^2}\rho\,\mathrm dx=1$ and $\int_{\mathbb R^2}|x|^2\rho\,\mathrm dx<\infty$.
Step~1 approximates $\rho$ by step densities with convergent effective
energies. Step~2 selects a single recovery family by a diagonal choice, and Step~3 applies the lower bound.

\par\medskip\noindent\textit{Step 1. Approximation with convergence of both effective energy terms.}
For $R>0$ set $m_R:=\int_{Q_R}\rho\,\mathrm dx$. Monotone convergence gives $m_R\uparrow1$, so $m_R>0$ for all sufficiently large $R$. Define $\rho_R=m_R^{-1}\rho\mathbf1_{Q_R}$ for these $R$. Then $\rho_R\geqslant0$, $\int_{\mathbb R^2}\rho_R\,\mathrm dx=1$, and
\begin{align}
 \norm{\rho_R-\rho}_{L^2(\mathbb R^2)}
 &\leqslant|m_R^{-1}-1|\norm{\rho\mathbf1_{Q_R}}_{L^2(\mathbb R^2)}
 +\norm{\rho\mathbf1_{\mathbb R^2\setminus Q_R}}_{L^2(\mathbb R^2)}\longrightarrow0,\\
 \int_{\mathbb R^2}(1+|x|^2)|\rho_R-\rho|\,\mathrm dx
 &\leqslant|m_R^{-1}-1|\int_{Q_R}(1+|x|^2)\rho\,\mathrm dx
 +\int_{\mathbb R^2\setminus Q_R}(1+|x|^2)\rho\,\mathrm dx\longrightarrow0.
\end{align}
For fixed $R$ and $k\in\mathbb N$, partition $Q_R$ into $k^2$ squares of side $h=R/k$. For $w \in L^2(Q_R)$ let $T_{R,k}w$ equal the average of $w$ on each square and vanish outside $Q_R$. Jensen's inequality gives
\begin{equation}
 \norm{T_{R,k}w}_{L^2(Q_R)}\leqslant\norm{w}_{L^2(Q_R)}.
\end{equation}
For $\varphi\in C(\overline{Q_R})$, let $\omega_\varphi$ be its modulus of continuity. Then
\begin{equation}
 \norm{T_{R,k}\varphi-\varphi}_{L^2(Q_R)}
 \leqslant |Q_R|^{1/2}\omega_\varphi(\sqrt2R/k)\longrightarrow0.
\end{equation}
Density of $C(\overline{Q_R})$ in $L^2(Q_R)$ and the contraction estimate imply, for $w\in L^2(Q_R)$,
\begin{equation}
 \limsup_{k\to\infty}\norm{T_{R,k}w-w}_{L^2(Q_R)}
 \leqslant2\norm{w-\varphi}_{L^2(Q_R)}\qquad(\forall \varphi\in C(\overline{Q_R})),
\end{equation}
and hence $T_{R,k}w\to w$ in $L^2(Q_R)$. Since both functions are supported in $Q_R$, Cauchy--Schwarz also gives
\begin{equation}
 \int_{\mathbb R^2}(1+|x|^2)|T_{R,k}w-w|\,\mathrm dx
 \leqslant(1+R^2/2)|Q_R|^{1/2}
 \norm{T_{R,k}w-w}_{L^2(Q_R)}\longrightarrow0.
\end{equation}
The preceding estimates allow us to choose $R_n\uparrow\infty$, $k_n\in\mathbb N$, with
\begin{equation}\label{eq:rhoapprox}
 \rho^{(n)}:=T_{R_n,k_n}\rho_{R_n},\qquad
 \norm{\rho^{(n)}-\rho}_{L^2(\mathbb R^2)}
 +\int_{\mathbb R^2}(1+|x|^2)|\rho^{(n)}-\rho|\,\mathrm dx\leqslant\frac1n.
\end{equation}
For example, select $R_n$ so that the sum of the two truncation errors is at most $1/(2n)$, and then $k_n$ so that the sum of the two averaging errors is at most $1/(2n)$. Positivity and preservation of the cell integrals yield
\begin{equation}
 \rho^{(n)}=\sum_{j=1}^{J_n}a_{n,j}\mathbf1_{S_{n,j}},
 \qquad a_{n,j}>0,\qquad
 \sum_{j=1}^{J_n}a_{n,j}|S_{n,j}|=1,
\end{equation}
after zero-height squares are omitted. The remaining squares have disjoint interiors. Thus every $\rho^{(n)}$ satisfies the hypotheses of Proposition~\ref{prop:steprecover}. Finally,
\begin{align}
 |\F(\rho^{(n)})-\F(\rho)|
 &\leqslant2\int_{\mathbb R^2}|x|^2|\rho^{(n)}-\rho|\,\mathrm dx\\
 &\quad+\beta_*\norm{\rho^{(n)}-\rho}_{L^2(\mathbb R^2)}
 \big(\norm{\rho^{(n)}}_{L^2(\mathbb R^2)}+\norm{\rho}_{L^2(\mathbb R^2)}\big)
 \overset{\eqref{eq:rhoapprox}}{\longrightarrow}0.
\end{align}

\par\medskip\noindent\textit{Step 2. Diagonal selection gives a single recovery family.}
For each $n$, Proposition~\ref{prop:steprecover} gives normalized finite-energy
$\psi_\eps^{(n)}\in\mathcal L_\eps$ recovering $\rho^{(n)}$.
We select one family by controlling the energy and successively more test functions. The construction below spells out the diagonal choice.

Let $(e_j)_{j\geqslant1}$ be an orthonormal basis of
$L^2(\mathbb R^2;\mathbb R)$. For each fixed $n$, the energy upper bound
and weak convergence in Proposition~\ref{prop:steprecover} hold below a
positive threshold. Reducing these thresholds recursively yields
$\eta_{n+1}<\eta_n/2$ and
\begin{equation}\label{eq:recovery-thresholds}
 0<\eps<\eta_n\quad\Longrightarrow\quad
 \begin{cases}
 \E_\eps(\psi_\eps^{(n)})\leqslant\F(\rho^{(n)})+1/n,\\[2mm]
 \displaystyle\max_{1\leqslant j\leqslant n}
 \left|\int_{\mathbb R^2}
 (|\psi_\eps^{(n)}|^2-\rho^{(n)})e_j\,\mathrm dx\right|<1/n.
 \end{cases}
\end{equation}
For $\eta_{n+1}\leqslant\eps<\eta_n$, define
$\psi_\eps:=\psi_\eps^{(n)}$.

For $\eps\in[\eta_1,1)$, we may choose any normalized Gaussian in $\Sigma\cap\mathcal L_\eps$; this completes the definition on $(0,1)$ without affecting the limit. Since $n\to\infty$ as $\eps\downarrow0$,
\begin{equation}\label{eq:diagonalenergy}
 \limsup_{\eps\downarrow0}\E_\eps(\psi_\eps)
 \overset{\eqref{eq:recovery-thresholds}}{\leqslant}
 \lim_{n\to\infty}\bigl(\F(\rho^{(n)})+1/n\bigr)=\F(\rho).
\end{equation}
For each fixed $j$, whenever $n\geqslant j$ and
$\eta_{n+1}\leqslant\eps<\eta_n$, we also have
\begin{equation}
 \left|\int_{\mathbb R^2}(|\psi_\eps|^2-\rho)e_j\,\mathrm dx\right|
 \overset{\eqref{eq:recovery-thresholds}}{\leqslant}
 \frac1n+\norm{\rho^{(n)}-\rho}_{L^2(\mathbb R^2)}
 \overset{\eqref{eq:rhoapprox}}{\leqslant}\frac2n.
\end{equation}
Thus the density difference converges to zero against every finite linear
combination of the $e_j$. Moreover,
$\||\psi_\eps|^2\|_{L^2(\mathbb R^2)}^2\leqslant\E_\eps(\psi_\eps)$
by \eqref{eq:modelintro} and \eqref{eq:Keps}, so
\eqref{eq:diagonalenergy} gives a uniform $L^2$ bound. Consequently, for
$\varphi\in L^2(\mathbb R^2;\mathbb R)$ and
$\varphi_0\in\operatorname{span}\{e_j:j\geqslant1\}$,
\begin{equation}
 \limsup_{\eps\downarrow0}
 \left|\int_{\mathbb R^2}(|\psi_\eps|^2-\rho)\varphi\,\mathrm dx\right|
 \leqslant C\norm{\varphi-\varphi_0}_{L^2(\mathbb R^2)}.
\end{equation}
Taking the infimum over this dense span proves
$|\psi_\eps|^2\rightharpoonup\rho$ in
$L^2(\mathbb R^2;\mathbb R)$.

\par\medskip\noindent\textit{Step 3. The lower bound gives equality of the limiting energies.}
Assertion~(1), proved in Subsection~\ref{sec:liminfproof}, now implies
\begin{equation}
 \F(\rho)\overset{\eqref{eq:directliminf}}{\leqslant}
 \liminf_{\eps\downarrow0}\E_\eps(\psi_\eps)
 \leqslant\limsup_{\eps\downarrow0}\E_\eps(\psi_\eps)
 \overset{\eqref{eq:diagonalenergy}}{\leqslant}\F(\rho).
\end{equation}
This proves \eqref{eq:directrecovery}.
\end{proof}

\section{Minimization and original variables}\label{sec:assembly}

We now prove Theorem~\ref{thm:min},
Proposition~\ref{prop:original}, and
Corollary~\ref{thm:intro}, in that order.
The compactness, lower bound, and LLL recovery in
Theorem~\ref{thm:gamma} give the limits of the minimum
energies and minimizing densities. The separate
interaction bound \eqref{eq:interactionlower} will
also give convergence of each energy contribution.
We then compare $\E_\eps$ with $\widehat\E_\eps$
and return to the original variables.
Until the proof of Corollary~\ref{thm:intro}, we use
only the intrinsic coefficient $\beta_*$ from
Lemma~\ref{lem:cell}; its identification with
$b_\triangle$ is proved in
Section~\ref{sec:certificate}.

\subsection{Minima, almost minimizers and persistent oscillations: Proof of Theorem~\ref{thm:min}}
\label{sec:min}

\begin{proof}[Proof of Theorem~\ref{thm:min}]
We first prove existence for each fixed $\eps>0$
and compute the minimum of $\F$.
We then use Theorem~\ref{thm:gamma} to pass to the
limit in the two minimization problems.
Finally, we prove the separate energy limits and
\eqref{eq:defect}.

\smallskip
\noindent\emph{Step 1. Existence in both admissible classes.}
Fix $\eps>0$ and recall that
$\Sigma=\{v\in H^1(\R^2;\C):
\abs{x}v\in L^2(\R^2;\C)\}$.
The function
$g_\eps(z)=(\pi\eps)^{-1/2}e^{-\abs{z}^2/(2\eps)}$
belongs to $\Sigma\cap\mathcal L_\eps$, has
unitary $L^2$ norm, and has finite energy.
Thus both infima in \eqref{eq:minvalue} are finite.

Let $(v_j)_{j\in\mathbb N}$ be a minimizing
sequence in either admissible class. By
\eqref{eq:modelintro} and \eqref{eq:Keps},
\begin{equation}
  \sup_{j\in\mathbb N}
  \left(
    \norm{A_\eps v_j}_{L^2(\R^2)}^2
    +\norm{\abs{x}v_j}_{L^2(\R^2)}^2
    +\norm{v_j}_{L^4(\R^2)}^4
  \right)<\infty.
\end{equation}
The factorization \eqref{eq:factor} and
$\nabla v_j=D_\eps v_j+i\eps^{-1}x^\perp v_j$
then give a uniform $H^1(\R^2;\C)$ bound,
with a constant that may depend on $\eps$.
After passing to a subsequence, we therefore have
$v_j\rightharpoonup v$ in $H^1(\R^2;\C)$ and
in $L^4(\R^2;\C)$, and
$\abs{x}v_j\rightharpoonup\abs{x}v$
in $L^2(\R^2;\C)$.
In particular, $v\in\Sigma$.

To retain the mass constraint, we need strong
$L^2$ convergence on the whole plane.
Rellich compactness gives this convergence on
bounded balls, while the moment bound gives
\begin{equation}
  \int_{\R^2\setminus B_R}\abs{v_j}^2\,\mathrm dx
  \leqslant R^{-2}
    \norm{\abs{x}v_j}_{L^2(\R^2)}^2
  \leqslant C_\eps R^{-2},
  \qquad R>0,\quad j\in\mathbb N.
\end{equation}
The same bound holds for $v$ by weak lower
semicontinuity. Combining these tail estimates
with local convergence gives
$v_j\to v$ in $L^2(\R^2;\C)$; hence
$\norm{v}_{L^2(\R^2)}=1$.
Moreover, $A_\eps v_j\rightharpoonup A_\eps v$
in $L^2(\R^2;\C)$.
Using $K_\eps=\norm{A_\eps\cdot}_{L^2(\R^2)}^2$
in \eqref{eq:modelintro}, weak lower
semicontinuity gives
$\E_\eps(v)\leqslant\liminf_{j\to\infty}\E_\eps(v_j)$.
If the sequence belongs to $\mathcal L_\eps$,
so does $v$, since this subspace is closed in
$L^2(\R^2;\C)$ by \eqref{eq:lllintro}.
Thus both minima are attained.

\smallskip
\noindent\emph{Step 2. The minimum and unique minimizer of $\F$.}
Let $s=\sqrt{2\beta_*/\pi}$ and let $\rho_*$
be the profile in \eqref{eq:TF}.
It is nonnegative, bounded, and compactly supported.
Radial integration gives
\begin{equation}
  \int_{\R^2}\rho_*\,\mathrm dx
  =\frac{2\pi}{\beta_*}
    \int_0^{\sqrt{s}}(s-r^2)r\,\mathrm dr
  =\frac{\pi s^2}{2\beta_*}=1.
\end{equation}
To prove its optimality, let
$\rho\in L^2(\R^2;\R)$ satisfy $\F(\rho)<\infty$.
Expanding the squared $L^2$ distance and using
$\int_{\R^2}(\rho-\rho_*)\,\mathrm dx=0$, we obtain
\begin{equation}\label{eq:TF-minimum-identity}
    \F(\rho)-\F(\rho_*)
    =\beta_*\norm{\rho-\rho_*}_{L^2(\R^2)}^2
   +2\int_{\{x\in\R^2:\abs{x}^2>s\}}
      (\abs{x}^2-s)\rho(x)\,\mathrm dx
    \geqslant0.
\end{equation}
Here we used
$\abs{x}^2+\beta_*\rho_*(x)=s$
when $\abs{x}^2\leqslant s$, whereas
$\rho_*(x)=0$ when $\abs{x}^2>s$.
Since $\beta_*>0$, equality in
\eqref{eq:TF-minimum-identity} forces
$\rho=\rho_*$ almost everywhere on $\R^2$.
This proves uniqueness. A further radial integration gives
\begin{equation}\label{eq:TF-energy-parts}
  2\int_{\R^2}\abs{x}^2\rho_*\,\mathrm dx
  =\beta_*\int_{\R^2}\rho_*^2\,\mathrm dx
  =\frac{\pi s^3}{3\beta_*}=\frac{2s}{3}.
\end{equation}
Consequently,
$\min\F=\F(\rho_*)=4s/3
=\frac43\sqrt{2\beta_*/\pi}=e_*$.

\smallskip
\noindent\emph{Step 3. Convergence of both minimum energies.}
Apply the LLL recovery assertion of
Theorem~\ref{thm:gamma} to $\rho_*$.
It gives mass-one functions
$u_\eps\in\Sigma\cap\mathcal L_\eps$ such that
$\E_\eps(u_\eps)\to e_*$.
The inclusion of the admissible classes implies
\begin{equation}\label{eq:wave-minima-upper}
  \min_{\substack{v\in\Sigma\\
                  \norm{v}_{L^2(\R^2)}=1}}\E_\eps(v)
  \leqslant
  \min_{\substack{v\in\Sigma\cap\mathcal L_\eps\\
                  \norm{v}_{L^2(\R^2)}=1}}\E_\eps(v)
  \leqslant\E_\eps(u_\eps)
  \longrightarrow e_*.
\end{equation}

For the opposite inequality, choose
$\eps_n\downarrow0$ along which the unrestricted
minimum approaches its lower limit, and let
$\psi_n\in\Sigma$ be a normalized minimizer
of $\E_{\eps_n}$.
By \eqref{eq:wave-minima-upper}, these energies
are bounded. Compactness in
Theorem~\ref{thm:gamma}(3) gives, after passage
to a subsequence,
$\abs{\psi_n}^2\rightharpoonup\rho$
in $L^2(\R^2;\R)$.
The lower bound in Theorem~\ref{thm:gamma}(1)
then yields
\begin{equation}
    \liminf_{\eps\downarrow0}
      \min_{\substack{v\in\Sigma
                      \norm{v}_{L^2(\R^2)}=1}}
        \E_\eps(v)
    =\lim_{n\to\infty}\E_{\eps_n}(\psi_n)
    \geqslant\F(\rho)\geqslant e_*.
\end{equation}
Together with \eqref{eq:wave-minima-upper},
this proves both limits in \eqref{eq:minvalue}.
This is the minimization argument associated with
$\Gamma$-convergence; see
\cite[Chapter~7]{DalMaso1993} and
\cite[Chapter~1]{Braides2002}.

\smallskip
\noindent\emph{Step 4. Convergence of minimizing densities.}
We prove a statement that also covers the almost
minimizers discussed after Theorem~\ref{thm:min}.
Suppose that $\psi_\eps\in\Sigma$ has mass one
and satisfies $\E_\eps(\psi_\eps)\to e_*$.
For every sequence $\eps_n\downarrow0$,
Theorem~\ref{thm:gamma}(3) provides a subsequence
whose densities converge weakly to some $\rho$
in $L^2(\R^2;\R)$.
The lower bound gives
$e_*\geqslant\F(\rho)\geqslant e_*$,
so uniqueness from Step~2 forces $\rho=\rho_*$.
Thus every parameter sequence has a subsequence
with this same weak density limit, which proves
$\abs{\psi_\eps}^2\rightharpoonup\rho_*$
for the whole family as $\eps\downarrow0$.

By \eqref{eq:minvalue}, the hypothesis
$\E_\eps(\psi_\eps)\to e_*$ holds for
minimizers of either problem. It also holds
whenever their energy excess above the corresponding
minimum tends to zero. This proves the stated
density convergence and the accompanying assertion
for almost minimizers in Section~\ref{sec:results}.

\smallskip
\noindent\emph{Step 5. Separate energy limits and the $L^2$ defect.}
Let $\psi_\eps$ be a minimizer in either admissible
class and set $\rho_\eps=\abs{\psi_\eps}^2$.
We have proved that
$\rho_\eps\rightharpoonup\rho_*$ and
$\E_\eps(\psi_\eps)\to e_*=\F(\rho_*)$.
The moment inequality \eqref{eq:momentlsc} and
the interaction bound \eqref{eq:interactionlower}
of Proposition~\ref{prop:lower} give
\begin{equation}
    \liminf_{\eps\downarrow0}
      2\int_{\R^2}\abs{x}^2\rho_\eps\,\mathrm dx
    \geqslant2\int_{\R^2}\abs{x}^2\rho_*\,\mathrm dx,\qquad
    \liminf_{\eps\downarrow0}
      \int_{\R^2}\rho_\eps^2\,\mathrm dx
    \geqslant\beta_*\int_{\R^2}\rho_*^2\,\mathrm dx.
\end{equation}
Since $K_\eps(\psi_\eps)\geqslant0$, subtracting
these lower bounds from the limiting total energy
in \eqref{eq:modelintro} yields
\begin{align}
    0\leqslant\limsup_{\eps\downarrow0}
      K_\eps(\psi_\eps)
    &\leqslant e_*
      -2\int_{\R^2}\abs{x}^2\rho_*\,\mathrm dx
      -\beta_*\int_{\R^2}\rho_*^2\,\mathrm dx=0,\\
    \limsup_{\eps\downarrow0}
      2\int_{\R^2}\abs{x}^2\rho_\eps\,\mathrm dx
    &\leqslant e_*-
      \beta_*\int_{\R^2}\rho_*^2\,\mathrm dx
      =2\int_{\R^2}\abs{x}^2\rho_*\,\mathrm dx,\\
    \limsup_{\eps\downarrow0}
      \int_{\R^2}\rho_\eps^2\,\mathrm dx
    &\leqslant e_*-
      2\int_{\R^2}\abs{x}^2\rho_*\,\mathrm dx
      =\beta_*\int_{\R^2}\rho_*^2\,\mathrm dx.
\end{align}
These inequalities prove
\eqref{eq:partK}--\eqref{eq:partI}.
Finally, weak convergence tested against
$\rho_*\in L^2(\R^2;\R)$ gives
\begin{equation}
    \norm{\rho_\eps-\rho_*}_{L^2(\R^2)}^2
    =\int_{\R^2}\rho_\eps^2\,\mathrm dx
      -2\int_{\R^2}\rho_\eps\rho_*\,\mathrm dx
      +\int_{\R^2}\rho_*^2\,\mathrm dx
    \longrightarrow
      (\beta_*-1)\norm{\rho_*}_{L^2(\R^2)}^2>0.
\end{equation}
Strict positivity follows from $\beta_*>1$ in
Lemma~\ref{lem:cell} and the fact that $\rho_*$
has mass one. This proves \eqref{eq:defect}
and completes the proof of Theorem~\ref{thm:min}.
\end{proof}

\begin{remark}[Projection onto an almost-minimizing LLL state]
Projection onto the LLL changes both the mass and
the confinement energy. We show that, after
normalization, it nevertheless gives an
almost-minimizing family for the restricted problem.

Let $\psi_\eps\in\Sigma$, $\eps>0$, be normalized
minimizers in either admissible class of
Theorem~\ref{thm:min}, and set
$p_\eps=P_\eps\psi_\eps$.
By \eqref{eq:partK}, \eqref{eq:r2}, and \eqref{eq:r4},
\begin{equation}\label{eq:minimizer-projection-closeness}
  \norm{\psi_\eps-p_\eps}_{L^2(\R^2)}
    =o(\eps^{1/2}),
  \qquad
  \norm{\psi_\eps-p_\eps}_{L^4(\R^2)}
    =o(\eps^{1/4}).
\end{equation}
Orthogonality also gives
$\norm{p_\eps}_{L^2(\R^2)}^2
=1-\norm{\psi_\eps-p_\eps}_{L^2(\R^2)}^2
=1-o(\eps)$.
Thus, for sufficiently small $\eps$, we may define
$q_\eps=p_\eps/\norm{p_\eps}_{L^2(\R^2)}$.
This normalized projection belongs to
$\Sigma\cap\mathcal L_\eps$ and is almost minimizing
for the LLL-restricted problem, in the sense that
\begin{equation}\label{eq:projected-almost-minimizer}
  \E_\eps(q_\eps)
  -\min_{\substack{v\in\Sigma\cap\mathcal L_\eps\\
                   \norm{v}_{L^2(\R^2)}=1}}
      \E_\eps(v)
  \longrightarrow0
  \qquad\text{as }\eps\downarrow0.
\end{equation}
To justify these assertions, we first control
the second moment after projection.
For almost every
$x\in\R^2$, the kernel representation gives
\begin{equation}
  xp_\eps(x)-[P_\eps(y\psi_\eps(y))](x)
  =\int_{\R^2}(x-y)P_\eps(x,y)\psi_\eps(y)\,\mathrm dy,
\end{equation}
where $P_\eps$ acts componentwise on
vector-valued functions.
Its kernel satisfies
$\abs{P_\eps(x,y)}
=(\pi\eps)^{-1}e^{-\abs{x-y}^2/(2\eps)}$
for $x,y\in\R^2$.
Since the function
$x\mapsto(\pi\eps)^{-1}\abs{x}e^{-\abs{x}^2/(2\eps)}$
has $L^1(\R^2)$ norm $C\sqrt\eps$,
Young's inequality and $L^2$ contractivity imply
\begin{equation}\label{eq:projected-moment-bound}
  \norm{\abs{x}p_\eps}_{L^2(\R^2)}
  \leqslant
    \norm{\abs{x}\psi_\eps}_{L^2(\R^2)}
    +C\sqrt\eps.
\end{equation}
Hence $p_\eps$ has finite second moment.
Since $p_\eps\in\mathcal L_\eps\subset\mathcal H_\eps$,
\eqref{eq:energy-domain} yields
$p_\eps,q_\eps\in\Sigma\cap\mathcal L_\eps$.

The energies $\E_\eps(\psi_\eps)$ are bounded by
Theorem~\ref{thm:min}, so their confinement and
quartic terms are bounded as well.
Consequently, \eqref{eq:projected-moment-bound},
\eqref{eq:minimizer-projection-closeness}, and
$\norm{p_\eps}_{L^2(\R^2)}=1-o(\eps)$ give
\begin{equation}
  \begin{aligned}
    \E_\eps(q_\eps)
    &=2\int_{\R^2}\abs{x}^2\abs{q_\eps}^2\,\mathrm dx
      +\int_{\R^2}\abs{q_\eps}^4\,\mathrm dx\\
    &\leqslant
      2\int_{\R^2}\abs{x}^2\abs{\psi_\eps}^2\,\mathrm dx
      +\int_{\R^2}\abs{\psi_\eps}^4\,\mathrm dx+o(1)\\
    &=\E_\eps(\psi_\eps)-K_\eps(\psi_\eps)+o(1)
    \leqslant\E_\eps(\psi_\eps)+o(1).
  \end{aligned}
\end{equation}
Because $\psi_\eps$ minimizes in one of the two
admissible classes, its energy is no larger than
the LLL-restricted minimum. Since $q_\eps$ is
admissible for that restricted problem,
\eqref{eq:projected-almost-minimizer} follows.

Finally, normalization changes $p_\eps$ by
$o(\eps)$ in both $L^2(\R^2;\C)$ and
$L^4(\R^2;\C)$. Thus the rates in
\eqref{eq:minimizer-projection-closeness} also
hold with $p_\eps$ replaced by $q_\eps$.
Each minimizer is therefore close to a normalized
almost-minimizing LLL state. Closeness in norm
to an exact LLL minimizer would require an
additional stability estimate uniform in $\eps$.
\end{remark}

\subsection{Return to the original rotating energy: Proof of Proposition~\ref{prop:original} and Corollary~\ref{thm:intro}}
\label{sec:original}

\begin{proof}[Proof of Proposition~\ref{prop:original}]
We first prove \eqref{eq:comparison} under a bound
on either rescaled energy. We then transfer the
lower bound and recovery families to
$\widehat\E_\eps$, and compare the infima and
almost minimizers of the two energies.

\smallskip
\noindent\emph{Step 1. The exact difference of the energies.}
Let $\eps\in(0,1)$, $\Omega=1-\eps^2$, and
$\psi\in\Sigma$ with
$\norm{\psi}_{L^2(\R^2)}=1$.
Since $\eps/\Omega>0$, the operators with this
parameter are defined as in
Subsection~\ref{sec:intro-scaling}, and
$D_{\eps/\Omega}=D_\eps+i\eps x^\perp$.
Expanding the square in \eqref{eq:exactscaled}
and subtracting \eqref{eq:modelintro}, we find
\begin{equation}\label{eq:rescaled-energy-difference}
  \begin{aligned}
    \widehat\E_\eps(\psi)-\E_\eps(\psi)
    &=2\Rea\int_{\R^2}
       D_\eps\psi\cdot
       \overline{i\eps x^\perp\psi}\,\mathrm dx
       +2\eps.
  \end{aligned}
\end{equation}
Indeed, the squared derivative adds
$\eps^2\int_{\R^2}\abs{x}^2\abs{\psi}^2\,\mathrm dx$,
whereas the confinement coefficient changes by
$1+\Omega-2=-\eps^2$; these terms cancel.
The constant term changes by
$2(1-\Omega)/\eps=2\eps$.
Cauchy--Schwarz therefore gives
\begin{equation}\label{eq:rescaled-difference-bound}
  \abs{\widehat\E_\eps(\psi)-\E_\eps(\psi)}
  \leqslant
    2\eps\norm{D_\eps\psi}_{L^2(\R^2)}
      \norm{\abs{x}\psi}_{L^2(\R^2)}+2\eps.
\end{equation}

\smallskip
\noindent\emph{Step 2. Uniform comparison on either bounded sublevel.}
Fix $M>0$.
If $\E_\eps(\psi)\leqslant M$,
\eqref{eq:modelintro} and \eqref{eq:factor} imply
\begin{equation}
  \norm{\abs{x}\psi}_{L^2(\R^2)}^2\leqslant\frac M2,
  \qquad
  \norm{D_\eps\psi}_{L^2(\R^2)}^2
    =K_\eps(\psi)+\frac2\eps
    \leqslant M+\frac2\eps.
\end{equation}
Substitution into \eqref{eq:rescaled-difference-bound}
gives \eqref{eq:comparison}.

If instead $\widehat\E_\eps(\psi)\leqslant M$,
apply \eqref{eq:factor} at parameter $\eps/\Omega$
to rewrite \eqref{eq:exactscaled} as
\begin{equation}\label{eq:exact-rescaled-positive}
  \widehat\E_\eps(\psi)
  =\norm{A_{\eps/\Omega}\psi}_{L^2(\R^2)}^2
   +(1+\Omega)\int_{\R^2}
      \abs{x}^2\abs{\psi}^2\,\mathrm dx
   +\int_{\R^2}\abs{\psi}^4\,\mathrm dx.
\end{equation}
All three terms are nonnegative. Hence
\begin{equation}
  \norm{\abs{x}\psi}_{L^2(\R^2)}^2\leqslant M,
  \qquad
  \norm{D_{\eps/\Omega}\psi}_{L^2(\R^2)}^2
    \leqslant M+\frac{2\Omega}{\eps}
    \leqslant M+\frac2\eps.
\end{equation}
Using $D_\eps=D_{\eps/\Omega}-i\eps x^\perp$,
we conclude that
$\norm{D_\eps\psi}_{L^2(\R^2)}
\leqslant(M+2/\eps)^{1/2}+\eps\sqrt M$.
Again \eqref{eq:rescaled-difference-bound}
gives \eqref{eq:comparison}, with a constant
$C_M$ depending only on $M$.
This proves the estimate whenever the minimum
of the two energies is at most $M$, as asserted
in Proposition~\ref{prop:original}; we may choose,
for example, $\eps_M=1/2$.

\smallskip
\noindent\emph{Step 3. The $\Gamma$-limit of the exact rescaled energy.}
Let $\eps_n\downarrow0$ and
$v_n\in L^2(\R^2;\C)\cap L^4(\R^2;\C)$
satisfy $\abs{v_n}^2\rightharpoonup\rho$
in $L^2(\R^2;\R)$.
If $\liminf_{n\to\infty}\widehat\E_{\eps_n}(v_n)
=+\infty$, the lower bound is immediate.
Otherwise, choose a subsequence $n_k$ realizing
this finite lower limit. Its exact energies are
bounded, so its terms belong to $\Sigma$ and
have mass one. By \eqref{eq:comparison},
\begin{equation}
    \liminf_{n\to\infty}\widehat\E_{\eps_n}(v_n)
    =\lim_{k\to\infty}\E_{\eps_{n_k}}(v_{n_k})
    \overset{\eqref{eq:directliminf}}{\geqslant}
      \F(\rho).
\end{equation}
This also shows that a limit density with
$\F(\rho)=+\infty$ cannot have a finite energy
lower limit.

For recovery, let $\F(\rho)<\infty$ and take
$\psi_\eps\in\Sigma\cap\mathcal L_\eps$
from \eqref{eq:directrecovery}.
Its model energies converge to $\F(\rho)$,
so \eqref{eq:comparison} gives
$\widehat\E_\eps(\psi_\eps)\to\F(\rho)$,
with the same weak density limit $\rho$.
We have proved both conditions in
Definition~\ref{def:density-convergence}.
Thus $\widehat\E_\eps$ has $\Gamma$-limit $\F$.

\smallskip
\noindent\emph{Step 4. Infima and almost minimizers.}
We also justify the consequences used after
Proposition~\ref{prop:original} in Section~\ref{sec:results}.
Let $\psi_\eps\in\Sigma$ be a normalized
minimizer of $\E_\eps$.
Its energy tends to $e_*$ by \eqref{eq:minvalue},
so \eqref{eq:comparison} gives
\begin{equation}
  \inf_{\substack{\psi\in\Sigma\\
                  \norm{\psi}_{L^2(\R^2)}=1}}
      \widehat\E_\eps(\psi)
  \leqslant\widehat\E_\eps(\psi_\eps)
  =\min_{\substack{\psi\in\Sigma\\
                   \norm{\psi}_{L^2(\R^2)}=1}}
      \E_\eps(\psi)+o(1).
\end{equation}
In particular, the exact infima are bounded
for all sufficiently small $\eps$.
Choose mass-one functions $\varphi_\eps\in\Sigma$
whose exact energies are within $\eps$ of these
infima. These energies are bounded, and
\eqref{eq:comparison} now gives the reverse bound:
\begin{equation}
    \min_{\substack{\psi\in\Sigma\\
                    \norm{\psi}_{L^2(\R^2)}=1}}
        \E_\eps(\psi)
    \leqslant\E_\eps(\varphi_\eps)
    \leqslant
      \inf_{\substack{\psi\in\Sigma
                      \norm{\psi}_{L^2(\R^2)}=1}}
          \widehat\E_\eps(\psi)+\eps+o(1).
\end{equation}
Using \eqref{eq:minvalue} in these two inequalities,
we obtain
\begin{equation}\label{eq:exact-minima-limit}
  \lim_{\eps\downarrow0}
    \inf_{\substack{\psi\in\Sigma\\
                    \norm{\psi}_{L^2(\R^2)}=1}}
      \widehat\E_\eps(\psi)=e_*.
\end{equation}
If a mass-one family $\psi_\eps\in\Sigma$ is
almost minimizing for $\widehat\E_\eps$, then
\eqref{eq:exact-minima-limit} implies
$\widehat\E_\eps(\psi_\eps)\to e_*$.
The comparison \eqref{eq:comparison} gives
$\E_\eps(\psi_\eps)\to e_*$ as well.
By \eqref{eq:minvalue}, it is therefore almost
minimizing for $\E_\eps$.
Step~4 of the proof of Theorem~\ref{thm:min}
then gives
$\abs{\psi_\eps}^2\rightharpoonup\rho_*$
in $L^2(\R^2;\R)$.
\end{proof}

\begin{proof}[Proof of Corollary~\ref{thm:intro}]
We first use the exact relation between the energies
to prove assertion~(i). We then rescale the family
in assertion~(ii) and apply the almost-minimizer
conclusion proved above.
Here we use $\beta_*=b_\triangle$ from
Theorem~\ref{thm:abrikosov}, whose proof in
Section~\ref{sec:certificate} is independent of
this corollary.

\smallskip
\noindent\emph{Step 1. The minimum-energy expansion.}
Let $\Omega\in(0,1)$ and $\eps=\sqrt{1-\Omega}$.
The dilation
$\Psi(y)=\sqrt\eps\,\psi(\sqrt\eps\,y)$,
$y\in\R^2$, is a bijection between the
mass-one functions in $\Sigma$ in the original
and rescaled coordinates.
Thus \eqref{eq:exactscaled} and
\eqref{eq:exact-minima-limit} imply
\begin{equation}
  \begin{aligned}
    \inf_{\substack{\Psi\in\Sigma\\
                    \norm{\Psi}_{L^2(\R^2)}=1}}
        E_\Omega^{\mathrm{GP}}(\Psi)
    &=2\Omega+\eps
      \inf_{\substack{\psi\in\Sigma\\
                      \norm{\psi}_{L^2(\R^2)}=1}}
          \widehat\E_\eps(\psi)\\
    &=2\Omega+\eps e_*+o(\eps).
  \end{aligned}
\end{equation}
This proves \eqref{eq:intro-energy}, with
$\eps=\sqrt{1-\Omega}$ and
$e_*=\frac43\sqrt{2b_\triangle/\pi}$ by
\eqref{eq:minvalue} and Theorem~\ref{thm:abrikosov}.

\smallskip
\noindent\emph{Step 2. Densities and mass on rescaled regions.}
Let $(\Psi_\Omega)_{\Omega\in(0,1)}$ satisfy
the hypotheses of assertion~(ii), including
\eqref{eq:intro-almost-minimizers}.
For $\eps\in(0,1)$ and $\Omega=1-\eps^2$, set
$\psi_\eps(x)=\eps^{-1/2}
\Psi_\Omega(x/\sqrt\eps)$ for almost every
$x\in\R^2$.
Then $\psi_\eps\in\Sigma$ and
$\norm{\psi_\eps}_{L^2(\R^2)}=1$.
The energy excess satisfies the exact identity
\begin{equation}
    \widehat\E_\eps(\psi_\eps)
      -\inf_{\substack{\psi\in\Sigma\\
                       \norm{\psi}_{L^2(\R^2)}=1}}
          \widehat\E_\eps(\psi)
    =\frac1\eps\left(
      E_\Omega^{\mathrm{GP}}(\Psi_\Omega)
      -\inf_{\substack{\Psi\in\Sigma
                       \norm{\Psi}_{L^2(\R^2)}=1}}
          E_\Omega^{\mathrm{GP}}(\Psi)\right)
      \longrightarrow0,
\end{equation}
where convergence follows from
\eqref{eq:intro-almost-minimizers}.
Step~4 of the proof of Proposition~\ref{prop:original}
therefore gives
$\rho_\eps:=\abs{\psi_\eps}^2
\rightharpoonup\rho_*$ in $L^2(\R^2;\R)$.
The formula for $\rho_\eps$ is precisely the one
in the corollary, and \eqref{eq:TF}, with
$\beta_*=b_\triangle$, gives
\eqref{eq:intro-TF-profile}.

Finally, let $B\subset\R^2$ be bounded and
Lebesgue-measurable.
Since $\mathbf1_B\in L^2(\R^2;\R)$,
weak convergence and the change of variables
$x=\sqrt\eps\,y$ yield
\begin{equation}
  \int_{\eps^{-1/2}B}
    \abs{\Psi_\Omega(y)}^2\,\mathrm dy
  =\int_B\rho_\eps(x)\,\mathrm dx
  \longrightarrow\int_B\rho_*(x)\,\mathrm dx.
\end{equation}
This is \eqref{eq:intro-mass-regions} and
completes the proof of assertion~(ii).
\end{proof}

\section{Identification of the triangular coefficient}\label{sec:certificate}

In Sections~\ref{sec:compactness}--\ref{sec:assembly}, we proved the density compactness, the lower-bound, and the LLL recovery statements of
Theorem~\ref{thm:gamma}. The limiting functional \eqref{eq:limit} is
\begin{equation}\label{cert:eq:recalled-limit}
  \F(\rho)
  =2\int_{\R^2}|x|^2\rho\,\mathrm dx
   +\beta_*\int_{\R^2}\rho^2\,\mathrm dx,
\end{equation}
defined on nonnegative densities of mass one in $L^2(\R^2)$ with finite second moment. To complete the proof of the theorem, it remains to show that $\beta_*=b_\triangle$. To that end, let us briefly recall how the coefficients $\beta_*$ and $b_\triangle$ are defined. For $L>0$, we let $Q_L=(-L/2,L/2)^2$ and identify $x=(x_1,x_2)$ with $z=x_1+ix_2$.
The local optimization problem associated with the free-cell constant \eqref{eq:cell} is given by
\begin{equation}\label{eq:recallcL}
c_L=\inf\left\{
\frac{1}{L^2}\int_{Q_L}|u|^4\,\mathrm dx:
u(z)=e^{-|z|^2/2}f(z),\quad
f\in\mathcal O(Q_L),\quad
\frac{1}{L^2}\int_{Q_L}|u|^2\,\mathrm dx=1
\right\},
\end{equation}
where $\mathcal O(Q_L)$ denotes the space of holomorphic functions on the open square. We impose no boundary conditions, require no holomorphic extension beyond $Q_L$, and impose no restrictions on the zeros of $f$.
For every fixed $L>0$, a standard application of the direct method
in the calculus of variations shows that the infimum in
\eqref{eq:recallcL} is attained and can therefore be replaced by a
minimum. 
The minimizing wave function is never unique: replacing a minimizer
$u_L$ by $e^{i\theta}u_L$, for any $\theta\in\mathbb R$, leaves both
its normalization and its energy unchanged. Uniqueness of the minimizer up to multiplication by a constant phase
lies beyond the scope of this paper.

As established in Lemma~\ref{lem:cell}, $\beta_*$ admits the intrinsic
characterization
\begin{equation}
\beta_*=\lim_{L\to\infty}c_L=\sup_{L>0}c_L.
\end{equation}
The comparison value is the triangular Abrikosov ratio
\begin{equation}
b_\triangle
=\sum_{k,l\in\mathbb Z}
\exp\left(-\frac{\pi}{t}|l-k\tau|^2\right)
\simeq1.1596,\qquad 
\tau=\frac12+i\frac{\sqrt3}{2},
\qquad t=\frac{\sqrt3}{2}.
\end{equation}
As recalled in \eqref{eq:btri} and realized in
Lemma~\ref{lem:periodic}, $b_\triangle$ is the normalized quartic
interaction ratio of a lowest-Landau-level state whose density has
triangular lattice periodicity. The only assertion of
Theorem~\ref{thm:gamma} that remains to be proved is the identity
$\beta_*=b_\triangle$ in \eqref{eq:celllimit}. This section proves
that identity together with a quantitative rate of convergence of
$c_L$ to $b_\triangle$. These arguments are independent of the
compactness, lower-bound, and recovery proofs of the preceding sections.

\begin{theorem}[Coefficient identification and finite-cell bound]
  \label{thm:abrikosov}
  The free-cell constants $c_L$ of \eqref{eq:cell} and the intrinsic
  coefficient $\beta_{\ast}$ satisfy
  \begin{equation}
    \label{eq:identified}
    \beta_{\ast} = \lim_{L\to\infty} c_L
      = \sup_{L>0} c_L = b_{\triangle}.
  \end{equation}
  More precisely, there are constants $C, L_0 > 0$ such that
  \begin{equation}
    \label{eq:cellrate}
    b_{\triangle} - CL^{-1/5} \leqslant c_L \leqslant
    b_{\triangle}  \qquad (L \geqslant L_0) .
  \end{equation}
\end{theorem}

Theorem~\ref{thm:abrikosov} shows that, as $L\to\infty$, this unrestricted
local minimization gives exactly $b_\triangle$, the normalized quartic
interaction ratio of the triangular periodic lowest-Landau-level state
constructed in Lemma~\ref{lem:periodic}. Thus the triangular value remains
optimal when the periodicity assumption is removed in the large-square
limit. This identifies the optimal energy coefficient; it does not
classify the microscopic configurations attaining or approaching that
value. For the underlying vortex-lattice
problem and the triangular lowest-Landau-level comparison, see
\cite{Abr57,ABD05,AB06,ABN06}.

Lemma~\ref{lem:cell} gives $c_L\leqslant b_\triangle$
for every $L>0$. We therefore need only prove
$c_L\geqslant b_\triangle-CL^{-1/5}$ for large $L$.
Fourier analysis and a weighted Poisson identity
first give a lower bound on magnetic tori in terms
of an auxiliary Gaussian minorant
(Proposition~\ref{cert:prop:periodic}).
Cutoff and projection transfer this bound to free
cells, with an error of order $L^{-1/5}$
(Proposition~\ref{cert:prop:free}).
The minorant's lattice contact conditions then
identify the lower-bound constant with
$b_\triangle$
(Proposition~\ref{cert:prop:contact}).
The required minorant is supplied by the sharp
Gaussian-minorant theorem proved by OpenAI in
\cite[Theorem~1.2]{Atomic}.

Subsection~\ref{cert:sec:framework} specifies the Fourier convention,
magnetic normalization, and periodic function spaces, and states the
external minorant theorem. Subsection~\ref{cert:sec:criterion} collects
the periodic lower bound, its transfer to free cells, and the
identification of the constant into a three-part criterion, then
deduces Theorem~\ref{thm:abrikosov}. These three assertions are proved,
respectively, in Subsections~\ref{cert:sec:lll},
\ref{cert:sec:transfer}, and~\ref{cert:sec:contact}.

\subsection{Normalization, magnetic tori, and the external
input}\label{cert:sec:framework}

To apply the sharp Gaussian-minorant theorem of
{\cite[Theorem~1.2]{Atomic}} directly, we first rescale the spatial variable
by $\sqrt{\pi}$ so that the Gaussian weight $e^{- |z|^2}$ of the free-cell
problem defining $c_L$ in \eqref{eq:cell} becomes $G (x) = e^{- \pi |x|^2}$ and fundamental
lattice cells have unit area. In these units, the effective magnetic field is
$2 \pi$ (one flux quantum per unit area), and $G$ is self-dual under the
unitary Fourier transform \eqref{cert:eq:fourier-convention} used in
{\cite{Atomic}}. Below we fix the Fourier and lattice conventions, rewrite the
free-cell constants $c_L$ in this normalization, introduce the auxiliary
magnetic-periodic spaces on the torus, and state the external minorant
theorem.

\subsubsection{Fourier and lattice conventions}\label{cert:sec:conventions}

We begin by fixing our Fourier convention and recording a few
lattice identities. Throughout this section, Fourier series and Parseval's
identity use the standard cell-normalized coefficients; see
\cite[Chapters~2--3]{SteinShakarchi2003}.
For $g$ in the Schwartz space $\Sch (\mathbb{R}^2)$, we define its Fourier
transform by
\begin{equation}
  \label{cert:eq:fourier-convention} \hat{g} (\xi) = \int_{\mathbb{R}^2} g (x)
  e^{- 2 \pi ix \cdot \xi} \hspace{0.17em} \mathrm{d} x, \qquad \xi \in
  \mathbb{R}^2 .
\end{equation}
With this normalization, the Gaussian
\begin{equation}
  \label{cert:eq:G} G (x) = e^{- \pi |x|^2}
\end{equation}
satisfies $\hat{G} = G$, and for any orthogonal map $R$ on $\mathbb{R}^2$ a
change of variables gives
\begin{equation}
  \label{cert:eq:rotation-covariance} \widehat{g \circ R} = \hat{g} \circ R.
\end{equation}
In particular, the Fourier transform of a radial function is radial.

A lattice $\Lambda = M\mathbb{Z}^2$, with $M \in \mathrm{GL}_2 (\mathbb{R})$,
has fundamental cell $M [0, 1)^2$, covolume $\mathrm{covol} (\Lambda) = | \det
M|$, and reciprocal (or dual) lattice\footnote{Here $M [0, 1)^2$ is shorthand
for the image of the half-open unit square $[0, 1)^2 = [0, 1) \times [0, 1)$
under the linear transformation defined by $M$: $M [0, 1)^2 = \{Mx : x \in [0,
1)^2 \}$. If $w_1, w_2 \in \mathbb{R}^2$ are the column vectors of the matrix
$M$ (which form a $\mathbb{Z}$-basis of the lattice $\Lambda =
M\mathbb{Z}^2$), then writing $x = (t_1, t_2) \in [0, 1)^2$ yields $M [0, 1)^2
= \{t_1 w_1 + t_2 w_2 : 0 \leqslant t_1, t_2 < 1\}$. Geometrically, this is
the half-open fundamental parallelogram spanned by $w_1$ and $w_2$, whose area
is given by $| \det M| = \mathrm{covol} (\Lambda)$. Because the upper and
right boundaries ($t_1 = 1$ or $t_2 = 1$) are excluded, the translates $p + M
[0, 1)^2$ over all lattice points $p \in \Lambda$ form a tiling that
partitions $\mathbb{R}^2$ without overlap. Strictly speaking, a dual basis
consists of linear forms acting on $\mathbb{R}^2$; however, via the Riesz
representation theorem under the standard Euclidean inner product, we identify
these linear forms with their corresponding vector representations
$w_1^{\ast}, w_2^{\ast} \in \mathbb{R}^2$. The dual lattice $\Lambda^{\ast}$
can then be characterized as the lattice generated by this dual basis,
uniquely defined by the biorthogonality conditions $w_i^{\ast} \cdot w_j =
\delta_{ij}$ (i.e., $w_1^{\ast} \cdot w_1 = w_2^{\ast} \cdot w_2 = 1$ and
$w_1^{\ast} \cdot w_2 = w_2^{\ast} \cdot w_1 = 0$), with its covolume
satisfying $\mathrm{covol} (\Lambda^{\ast}) = \mathrm{covol} (\Lambda)^{-
1}$.}
\begin{equation}
  \label{cert:eq:dual-definition} \Lambda^{\ast} = \{\xi \in \mathbb{R}^2 :
  \xi \cdot p \in \mathbb{Z} \text{ for all } p \in \Lambda\} = M^{- T}
  \mathbb{Z}^2 .
\end{equation}
We will repeatedly use the Poisson summation formula
\begin{equation}
  \label{cert:eq:poisson-general} \sum_{p \in \Lambda} g (p) =
  \frac{1}{\mathrm{covol} (\Lambda)}  \sum_{q \in \Lambda^{\ast}} \hat{g} (q),
  \qquad g \in \Sch (\mathbb{R}^2),
\end{equation}
where both sums converge absolutely because $g$ and $\hat{g}$ are Schwartz
functions {\cite[Chapter~5, Section~3]{SteinShakarchi2003}.

Let $J : \mathbb{R}^2 \to \mathbb{R}^2$ denote the linear
transformation corresponding to a counterclockwise rotation by $\pi / 2$,
defined explicitly for all $(x_1, x_2) \in \mathbb{R}^2$ by $J (x_1, x_2) = (-
x_2, x_1)$. Then, for all vectors $\xi, \eta \in \mathbb{R}^2$, the
determinant of the matrix formed by $\xi$ and $\eta$ satisfies
\begin{equation}
  \label{cert:eq:detJ} \det (\xi, \eta) = (- J \eta) \cdot \xi .
\end{equation}
When $\mathrm{covol} (\Lambda) = 1$, any positively oriented basis $w_1, w_2$
of $\Lambda$ has $\det (w_1, w_2) = 1$, so \eqref{cert:eq:detJ} implies that
its dual basis is $- Jw_2, Jw_1$. For instance, from \eqref{cert:eq:detJ} one
immediately gets $(- Jw_2) \cdot w_1 = \det (w_1, w_2) = 1$ and $(- Jw_2)
\cdot w_2 = \det (w_2, w_2) = 0$; similarly one gets $(Jw_1) \cdot w_1 = 0$
and $(Jw_1) \cdot w_2 = 1$. Thus, the vectors $- Jw_2$ and $Jw_1$ satisfy the
biorthogonality conditions $w_i^{\ast} \cdot w_j = \delta_{ij}$, confirming
them as the explicit vector representation of the dual basis.

Summarizing, when $\mathrm{covol} (\Lambda) = 1$, we have
\begin{equation}
  \label{cert:eq:dual-rotation} \Lambda^{\ast} = J \Lambda .
\end{equation}
In particular, the unit-covolume triangular lattice from
\eqref{eq:intro-btri},
\begin{equation}
  \label{cert:eq:triangular-lattice} \Lambda_{\triangle} = t^{- 1 / 2}  \{(m +
  n / 2, tn) : m, n \in \mathbb{Z}\}, \qquad t = \sqrt{3} / 2,
\end{equation}
satisfies $\Lambda_{\triangle}^{\ast} = J \Lambda_{\triangle}$. 
Substituting $(m,n)=(l,-k)$ into \eqref{eq:btri} identifies
$b_\triangle$ with the sum of $G$ over $\Lambda_\triangle$.
Since $\mathrm{covol}(\Lambda_\triangle)=1$ and $\hat G=G$,
Poisson summation gives
\begin{equation}
  \label{cert:eq:btri-lattice} b_{\triangle} \overset{\eqref{eq:btri},
  \eqref{cert:eq:triangular-lattice}}{=} \sum_{p \in \Lambda_{\triangle}} G
  (p) \overset{\eqref{cert:eq:dual-rotation}}{=} \sum_{q \in
  \Lambda_{\triangle}^{\ast}} G (q),
\end{equation}
We will use this identity in Proposition~\ref{cert:prop:contact}.

\subsubsection{Free cells and the change of gauge}
\label{cert:sec:normalization}
To prove the lower bound in Theorem~\ref{thm:abrikosov},
\begin{equation}
  c_L\geqslant b_\triangle-CL^{-1/5}
  \qquad (L\geqslant L_0),
\end{equation}
with $C,L_0>0$ as in \eqref{eq:cellrate}, we first estimate the
quartic quotient on magnetic tori and then transfer the estimate
to free cells.
We recall that the free-cell problem is posed on a square, with
$u=e^{-\abs{z}^2/2}f$ and $f$ holomorphic, and no boundary
conditions are imposed. On a magnetic torus, the values of $u$
on opposite sides of a rectangle are related by prescribed
phase factors, so that $\abs{u}^2$ is periodic. These boundary
conditions allow us to compute the quartic quotient by Fourier
series.

We begin by rescaling the free-cell problem to match the
normalization of the Gaussian minorant used below. We then pass
to Landau gauge, which allows ordinary periodicity in one
coordinate and makes the Fourier calculation simpler.
Both transformations leave the normalized quartic quotient
unchanged, as we now verify.

\smallskip

\paragraph{\emph{Rescaled free cells.}}
On $Q_L=(-L/2,L/2)^2$, define
\begin{equation}\label{eq:cell_unit_density}
\tilde c_L
=\inf_{{\begin{array}{c}
f\in\mathcal O(Q_L)\\
0<\int_{Q_L}\abs{u}^2<\infty
\end{array}}}
\frac{L^2\int_{Q_L}\abs{u}^4}
     {\left(\int_{Q_L}\abs{u}^2\right)^2},
\qquad
u(z)=e^{-\pi|z|^2/2}f(z).
\end{equation}
The quotient is understood as $+\infty$ if the quartic integral
diverges. No boundary condition is imposed, and multiplication
of $u$ by a nonzero constant leaves the quotient unchanged.

The dilation $\zeta=\sqrt\pi\,z$ maps $Q_L$ onto
$Q_{\sqrt\pi L}$ and gives a bijection between the admissible
classes in \eqref{eq:cell_unit_density} and
\eqref{eq:cellquotient}:
\begin{equation}\label{cert:eq:cell-dilation}
h(\zeta)=u(\zeta/\sqrt\pi)
=e^{-|\zeta|^2/2}f(\zeta/\sqrt\pi),
\qquad
\int_{Q_{\sqrt\pi L}}|h|^p
=\pi\int_{Q_L}|u|^p
\quad(p=2,4).
\end{equation}
Consequently,
\begin{equation}\label{cert:eq:quotient-dilation}
\frac{(\sqrt\pi L)^2\int_{Q_{\sqrt\pi L}}|h|^4}
     {\left(\int_{Q_{\sqrt\pi L}}|h|^2\right)^2}
\overset{\eqref{cert:eq:cell-dilation}}{=}
\frac{L^2\int_{Q_L}|u|^4}
     {\left(\int_{Q_L}|u|^2\right)^2}.
\end{equation}
Taking infima and applying Lemma~\ref{lem:cell}, we obtain
\begin{equation}\label{eq:certscale}
\tilde c_L
\overset{\eqref{cert:eq:quotient-dilation}}{=}
c_{\sqrt\pi L}
\overset{\text{Lemma~\ref{lem:cell}}}{\leqslant}
b_\triangle
\qquad(L>0).
\end{equation}
Thus, any explicit lower bound for $\tilde c_L$ yields the
corresponding bound for $c_L$ by substituting $L/\sqrt\pi$
for $L$, since $c_L=\tilde c_{L/\sqrt\pi}$ by
\eqref{eq:certscale}.

\begin{remark}[Meaning of the quotient]
\label{cert:rem:quotient}
For an admissible density $\rho=|u|^2\in L^2(Q_L)$, let
$\bar\rho=L^{-2}\int_{Q_L}\rho$ denote its spatial mean.
Expanding the square and using
$\int_{Q_L}(\rho-\bar\rho)=0$ gives
\begin{equation}\label{eq:quotienttwotofourtemp}
\frac{L^2\int_{Q_L}\rho^2}
     {\left(\int_{Q_L}\rho\right)^2}
=1+\frac1{L^2}\int_{Q_L}
\left(\frac{\rho}{\bar\rho}-1\right)^2.
\end{equation}
Thus minimizing the quotient amounts to minimizing the
$L^2$ distance of $\rho/\bar\rho$ from the constant function one.

By \eqref{eq:certscale} and Lemma~\ref{lem:cell},
we have $\tilde c_L=c_{\sqrt\pi L}>1$.
Combining this strict inequality with
\eqref{eq:quotienttwotofourtemp}, we obtain
\begin{equation}\label{cert:eq:density-deviation}
\frac1{L^2}\int_{Q_L}
\left(\frac{\rho}{\bar\rho}-1\right)^2
\overset{\eqref{eq:cell_unit_density},
         \eqref{eq:quotienttwotofourtemp}}{\geqslant}
\tilde c_L-1>0
\end{equation}
for every admissible density on the fixed square $Q_L$.
\end{remark}

\smallskip

\smallskip
\smallskip
\paragraph{\emph{Gauge convention.}}\label{cert:sec:gauge}
We have so far used symmetric gauge because it describes LLL
states as a Gaussian times a holomorphic function. This is the
description used in the free-cell problem \eqref{eq:cell} and
in the preceding holomorphic estimates and compactness arguments.
For the torus calculation, we pass to Landau gauge, which allows
us to use Fourier series in one coordinate.

We first recall the relation with the magnetic kinetic energy.
In the rescaled variables, the symmetric-gauge potential is
$\mathcal A(x)=\pi(-x_2,x_1)$. For
$u\in C_c^\infty(\R^2;\C)$, expansion of the square and integration
by parts give
\begin{equation}\label{cert:eq:symmetric-factorization}
    \int_{\R^2}\abs{(\nabla-i\mathcal A)u}^2
    ={}2\pi\int_{\R^2}\abs{u}^2
    +\int_{\R^2}
    \abs{\bigl[(\partial_{x_1}-i\mathcal A_1)
    +i(\partial_{x_2}-i\mathcal A_2)\bigr]u}^2.
\end{equation}
By a standard density argument, this identity extends to every
$u\in L^2(\R^2;\C)$ such that
$(\nabla-i\mathcal A)u\in L^2(\R^2;\C^2)$, with derivatives
understood in the distributional sense. The minimum magnetic kinetic
energy per unit mass is therefore $2\pi$, attained, for example,
by $u=e^{-\pi\abs{z}^2/2}$. The lowest Landau level is therefore achieved
by all the states for which the last integral vanishes.
Given that with $2\partial_{\bar z}=\partial_{x_1}+i\partial_{x_2}$, one has $A:=(\partial_{x_1}-i\mathcal A_1)+i(\partial_{x_2}-i\mathcal A_2)=(2\partial_{\bar z}+\pi z)$, the (global) LLL states are defined as the solutions of the distributional equation
\begin{equation}\label{cert:eq:symmetric-lll}
  A u=(2\partial_{\bar z}+\pi z)u=0
\end{equation}
We now consider this equation on a bounded open rectangle
$Q\subset\R^2$, in particular on $Q_L$. Given
$u\in L^2(Q;\C)$, set $f=e^{\pi\abs{z}^2/2}u$.
The Gaussian and its reciprocal are bounded on $Q$, so
$f\in L^2(Q;\C)$. No differentiability or holomorphicity
is assumed at this stage. The distributional product rule gives
\begin{equation}\label{cert:eq:local-holomorphic-identity}
  (2\partial_{\bar z}+\pi z)u
  =2e^{-\pi\abs{z}^2/2}\partial_{\bar z}f
  \qquad\text{in }\mathcal D'(Q).
\end{equation}
Consequently, $(2\partial_{\bar z}+\pi z)u=0$ if and only if
$\partial_{\bar z}f=0$ in $\mathcal D'(Q)$. By the distributional
Cauchy--Riemann characterization
\eqref{eq:local-lll-representation}, this is equivalent to
$f$ having a representative in $\mathcal O(Q)$. Thus the
$L^2(Q;\C)$ solutions are precisely
$u=e^{-\pi\abs{z}^2/2}f$ with
$f\in\mathcal O(Q)\cap L^2(Q;\C)$; in particular,
$u\in C^\infty(Q;\C)$.

For $Q=Q_L$, the nonzero solutions form exactly the admissible
class in \eqref{eq:cell_unit_density}. The quartic quotient is
finite precisely when $u\in L^4(Q_L;\C)$, equivalently
$f\in L^4(Q_L;\C)$. This characterization uses only the
equation inside $Q_L$, without boundary conditions.
The whole-plane factorization
\eqref{cert:eq:symmetric-factorization} explains why
$(2\partial_{\bar z}+\pi z)u=0$ is called the LLL equation.
On $Q_L$, we use this differential equation simply as an
equivalent formulation of the holomorphic constraint in
\eqref{eq:cell_unit_density}. The change of gauge below
uses this local equation and the product rule.

For the torus calculation, we use the Landau potential
\begin{equation}\label{cert:eq:landau-potential}
  \mathcal A^{\mathbb T}(x)=(0,2\pi x_1),
  \qquad
  \mathcal A^{\mathbb T}-\mathcal A=\nabla\chi,
  \qquad
  \chi(x)=\pi x_1x_2.
\end{equation}
On any bounded rectangle $Q$, multiplication by $e^{i\chi}$
defines a unitary map on $L^2(Q;\C)$. Since $\chi$ is real, setting $v=e^{i\chi}u$,
we have $\abs{v}^2=\abs{u}^2$ almost everywhere, and the
distributional product rule gives
\begin{equation}\label{cert:eq:gaugecov}
  (\nabla-i\mathcal A^{\mathbb T})v
  =e^{i\chi}(\nabla-i\mathcal A)u
  \qquad\text{in }\mathcal D'(Q;\C^2).
\end{equation}
Since $e^{\pm i\chi}\in W^{1,\infty}(Q)$,
$u\in H^1(Q;\C)$ if and only if $v\in H^1(Q;\C)$.
For such states, \eqref{cert:eq:gaugecov} gives
\begin{equation}\label{cert:eq:gauge-kinetic}
  \int_Q\abs{(\nabla-i\mathcal A^{\mathbb T})v}^2
  \overset{\eqref{cert:eq:gaugecov}}{=}
  \int_Q\abs{e^{i\chi}(\nabla-i\mathcal A)u}^2
  =\int_Q\abs{(\nabla-i\mathcal A)u}^2.
\end{equation}
Thus the density and magnetic kinetic energy remain unchanged under the
simultaneous replacement $(\mathcal{A}, u) \mapsto (\mathcal{A}^{\mathbb{T}},v)$.

In Landau gauge, the same combination of magnetic derivatives that led to \eqref{cert:eq:symmetric-lll} now gives
\begin{equation}\label{cert:eq:covariant}
  A^{\mathbb T}
  :=(\partial_{x_1}-i\mathcal A^{\mathbb T}_1)
   +i(\partial_{x_2}-i\mathcal A^{\mathbb T}_2)
  =2\partial_{\bar z}+2\pi x_1.
\end{equation}
For $u\in L^2(Q;\C)$, gauge covariance gives
\begin{equation}
  A^{\mathbb T}(e^{i\chi}u)=e^{i\chi}Au.
\end{equation}
Thus $u$ is an LLL state in symmetric gauge if and only if
$v=e^{i\chi}u$ is an LLL state in Landau gauge. In particular,
the symmetric-gauge state $u(z)=e^{-\pi\abs{z}^2/2}f(z)$
corresponds to
\begin{equation}
  v(z)=e^{i\pi x_1x_2}e^{-\pi\abs{z}^2/2}f(z),
  \qquad f\in\mathcal O(Q)\cap L^2(Q;\C),
\end{equation}
in Landau gauge.
Multiplication by $e^{i\chi}$ therefore identifies the local
LLL states in the two gauges and preserves their quartic
quotient \eqref{eq:cell_unit_density}.

The advantage of Landau gauge is that $\mathcal A^{\mathbb T}$
is independent of $x_2$, allowing ordinary periodicity in
that variable. Let $L_1,L_2>0$, with $L_1L_2\in\mathbb N$,
and set $Q=(0,L_1)\times(0,L_2)$. For the torus problem,
we consider $v\in L^2_{\mathrm{loc}}(\R^2;\C)$ satisfying
$A^{\mathbb T}v=0$ in $\mathcal D'(\R^2)$ and
\begin{equation}\label{cert:eq:landau-bc-preview}
  v(x_1+L_1,x_2)=e^{2\pi iL_1x_2}v(x_1,x_2),
  \qquad
  v(x_1,x_2+L_2)=v(x_1,x_2)
\end{equation}
almost everywhere. The two orders of translation differ
by the factor $e^{2\pi iL_1L_2}=1$, so these conditions
are compatible.

The local characterization above gives
$v\in C^\infty(\R^2;\C)$. By periodicity in $x_2$, we can write
\begin{equation}\label{cert:eq:landau-fourier-preview}
  \begin{aligned}
    v(x_1,x_2)
    &=\sum_{m\in\Z}d_m(x_1)e^{2\pi imx_2/L_2},\\
    d_m(x_1)
    &=\frac{1}{L_2}\int_0^{L_2}
      v(x_1,x_2)e^{-2\pi imx_2/L_2}\,\mathrm dx_2.
  \end{aligned}
\end{equation}
Here $d_m\in C^\infty(\R;\C)$, and the series converges in
$C^\infty_{\mathrm{loc}}(\R^2)$. Substituting into
$A^{\mathbb T}v=0$ and comparing Fourier coefficients gives
\begin{equation}\label{cert:eq:landau-mode-equation}
  d_m'(x_1)+2\pi\left(x_1-\frac{m}{L_2}\right)d_m(x_1)=0,
  \qquad x_1\in\R,\quad m\in\Z.
\end{equation}
We will use these equations together with the first boundary
condition in \eqref{cert:eq:landau-bc-preview} to construct
the torus LLL basis and derive the quartic Fourier formula
in Lemma~\ref{cert:lem:basis}.

\subsubsection{Magnetic tori}\label{cert:sec:torus}

We now introduce the periodic problem used to prove the lower
bound in Proposition~\ref{cert:prop:periodic}. The magnetic
boundary conditions allow us to describe the LLL by finitely
many coefficients and compute its quartic energy by Fourier
series. We will then transfer the estimate to free cells by
cutoff and projection in Proposition~\ref{cert:prop:free}.
We recall the magnetic-periodic setting at field $2\pi$,
fixing the notation used in the Fourier calculation and
the projection estimates below. The facts recalled here
are standard; see \cite[Sections~4--5]{Onofri2001} and
\cite[Sections~2 and~6--7]{SakamotoTanimura2003}.

\smallskip
\noindent\emph{Boundary conditions.}
Let $L_1,L_2>0$ and $Q=(0,L_1)\times(0,L_2)$.
Throughout this subsection, we work in Landau gauge and
denote the wave function by $u$. We call
$u\in L^2_{\mathrm{loc}}(\R^2;\C)$ magnetic-periodic if
\begin{equation}\label{cert:eq:bc}
  u(x_1+L_1,x_2)=e^{2\pi iL_1x_2}u(x_1,x_2),
  \qquad
  u(x_1,x_2+L_2)=u(x_1,x_2)
\end{equation}
almost everywhere. These are translation identities on
$\R^2$; no boundary trace is required for an $L^2$ state.
Applying the translations in the two possible orders gives
\begin{equation}\label{cert:eq:translation-compatibility}
  e^{2\pi iL_1(x_2+L_2)}u(x_1,x_2)
  =e^{2\pi iL_1x_2}u(x_1,x_2).
\end{equation}
For a nonzero state, we must therefore have
$e^{2\pi iL_1L_2}=1$, or equivalently
\begin{equation}\label{cert:eq:fluxquant}
  N:=L_1L_2\in\mathbb N.
\end{equation}
We assume this condition from now on. Iterating
\eqref{cert:eq:bc}, we obtain
\begin{equation}\label{cert:eq:bciter}
  u(x_1+kL_1,x_2+\ell L_2)
  =e^{2\pi ikL_1x_2}u(x_1,x_2),
  \qquad k,\ell\in\Z.
\end{equation}
Conversely, \eqref{cert:eq:bciter} extends any function in
$L^2(Q;\C)$ uniquely, up to equality almost everywhere,
to a magnetic-periodic function on $\R^2$.
We equip these functions with the inner product and norm
\begin{equation}\label{cert:eq:torus-inner-product}
  \langle u,v\rangle=\int_Q\bar u v\,\mathrm dx,
  \qquad
  \norm{u}^2=\int_Q\abs{u}^2\,\mathrm dx.
\end{equation}
Restriction to $Q$ then identifies this Hilbert space
isometrically with $L^2(Q;\C)$. The phase factors cancel
in $\abs{u}^2$ and $\bar u v$, so these functions are
periodic and their integrals are independent of the
chosen fundamental rectangle. We call $Q$, with opposite
sides identified through \eqref{cert:eq:bc}, a
\emph{magnetic torus}.

\smallskip
\smallskip
\noindent\emph{The lowest Landau level.}
We now impose the LLL equation on the magnetic-periodic
functions defined above. The resulting space is
finite-dimensional, which will allow us to express the
quartic energy in terms of finitely many coefficients.

\begin{definition}[Torus LLL]\label{cert:def:lll}
The space $\mathcal L^{\mathbb T}_{L_1,L_2}$ consists of
the functions $u\in L^2_{\mathrm{loc}}(\R^2;\C)$ satisfying
the magnetic-periodic conditions \eqref{cert:eq:bc} and
\begin{equation}
  A^{\mathbb T}u
  =(2\partial_{\bar z}+2\pi x_1)u=0
  \qquad\text{in }\mathcal D'(\R^2).
\end{equation}
\end{definition}
We identify these functions with their restrictions to $Q$
and use the $L^2(Q;\C)$ inner product.
Under this identification, $\mathcal L^{\mathbb T}_{L_1,L_2}$
is a closed subspace of $L^2(Q;\C)$: convergence in
$L^2(Q;\C)$ implies convergence of the magnetic-periodic
extensions in $L^2_{\mathrm{loc}}(\R^2;\C)$, and the
equation passes to the limit in distributions.
Moreover, as already observed in previous circumstances, every element of $\mathcal L^{\mathbb T}_{L_1,L_2}$ has a
smooth representative. Indeed, for
$f=e^{\pi\abs{z}^2/2-i\pi x_1x_2}u
\in L^2_{\mathrm{loc}}(\R^2;\C)$, the product rule gives
\begin{equation}\label{cert:eq:torus-holomorphic}
  2\partial_{\bar z}f
  =e^{\pi\abs{z}^2/2-i\pi x_1x_2}A^{\mathbb T}u
  =0
  \qquad\text{in }\mathcal D'(\R^2).
\end{equation}
Thus $f$ has an entire holomorphic representative, and
$u=e^{i\pi x_1x_2}e^{-\pi\abs{z}^2/2}f$ is smooth.

The dimension of this space is a standard consequence
of the theta-function construction; see
\cite[Section~5]{Onofri2001} and
\cite[Section~VII, equation~(7.11)]{SakamotoTanimura2003}.
More concretely, periodicity in $x_2$ allows us to expand
$u$ in Fourier modes. The LLL equation determines each
mode up to a constant, while the magnetic boundary
condition in $x_1$ identifies the constants whose indices
differ by $N$. We are therefore left with $N$ independent
coefficients, and
\begin{equation}\label{cert:eq:torus-dimension}
  \dim\mathcal L^{\mathbb T}_{L_1,L_2}=N=L_1L_2.
\end{equation}
In Lemma~\ref{cert:lem:basis}, we recall the resulting
orthonormal basis \eqref{cert:eq:basis}. Every
$u\in\mathcal L^{\mathbb T}_{L_1,L_2}$ then has a unique
representation
\begin{equation}
  u=\sum_{r=0}^{N-1}c_r\varphi_r,
  \qquad (c_0,\ldots,c_{N-1})\in\C^N.
\end{equation}
We use this representation to derive the quartic
Fourier formula \eqref{cert:eq:quartic}.

For the projection estimates, we also need the magnetic
kinetic-energy identity. With
$D^{\mathbb T}=\nabla-i\mathcal A^{\mathbb T}
=(\partial_{x_1},\partial_{x_2}-2\pi ix_1)$,
every smooth magnetic-periodic function $h$ satisfies
\begin{equation}\label{cert:eq:torus-factorization}
  \norm{D^{\mathbb T}h}_{L^2(Q)}^2
  =\norm{A^{\mathbb T}h}_{L^2(Q)}^2
   +2\pi\norm{h}_{L^2(Q)}^2.
\end{equation}
This follows by integration by parts: the magnetic
boundary conditions cancel the contributions from
opposite sides. We use this identity in
Lemma~\ref{cert:lem:projection} to estimate the projection
error uniformly in the size of a square torus.

\begin{remark}[Zeros and boundary conditions]
\label{cert:rem:zeros}
Let $u\in\mathcal L^{\mathbb T}_{L_1,L_2}\setminus\{0\}$.
Near each zero $z_j$, the holomorphic representation gives
$u(z)=(z-z_j)^{m_j}a_j(z)$, where $m_j\in\mathbb N$
and $a_j$ is smooth and nonvanishing. Thus, for small $r>0$,
\begin{equation}
  \deg\left(\frac{u}{\abs{u}},\partial B_r(z_j)\right)=m_j:
\end{equation}
the phase winds $m_j$ times around the zero.
If $z_1,\ldots,z_k$ represent the distinct zeros modulo
$L_1\Z\times L_2\Z$, then
\begin{equation}
  \sum_{j=1}^k m_j
  =\frac{1}{2\pi}\int_Q
    \operatorname{curl}\mathcal A^{\mathbb T}\,\mathrm dx
  =L_1L_2=N;
\end{equation}
see \cite[Theorem~7.6]{DFSNotes}.
Indeed, on a fundamental rectangle whose boundary avoids
the zeros, the horizontal phase changes cancel, while
the vertical boundary factor $e^{2\pi iL_1x_2}$ produces
a total phase change of $2\pi N$. The boundary conditions
therefore fix the total multiplicity, but not the locations
of the zeros.

No such restriction holds on a free cell:
$f\equiv1$ gives the zero-free state $u=e^{-\pi\abs{z}^2/2}$.
To apply the periodic lower bound, we therefore use cutoff
and projection in Proposition~\ref{cert:prop:free} to
construct a magnetic-periodic LLL state, controlling the
change in mass and quartic energy without assumptions
on the original zeros.
\end{remark}

\subsubsection{The external Gaussian minorant}
\label{cert:sec:minorant}

The identification $\beta_*=b_\triangle$ proceeds in two steps.
First, we prove a conditional statement: if there exists a real
radial Schwartz function satisfying the sign and contact conditions
below, then $\beta_*=b_\triangle$
(Theorem~\ref{cert:thm:criterion}).
Second, we invoke the external theorem
\cite[Theorem~1.2]{Atomic} to establish the existence of such
a function. We now state the required conditions and the
existence result.

\smallskip
\noindent\emph{Conditions on the minorant.}
For a real radial function $g\in\Sch(\R^2)$, define
\begin{equation}\label{cert:eq:Bg}
\kappa_g:=1+\hat g(0)-g(0).
\end{equation}
With $G(x)=e^{-\pi|x|^2}$ as in \eqref{cert:eq:G}, we require
the \emph{sign conditions}
\begin{equation}\label{eq:signconds}
g(x)\leqslant G(x),
\qquad
\hat g(\bm\xi)\geqslant0
\qquad(x,\bm\xi\in\R^2),
\end{equation}
and the \emph{contact conditions}
\begin{equation}\label{eq:latticecontactconds}
\begin{aligned}
g(\bm p)&=G(\bm p)
&&(\bm p\in\Lambda_\triangle\setminus\{0\}),\\
\hat g(\bm q)&=0
&&(\bm q\in\Lambda_\triangle^*\setminus\{0\}).
\end{aligned}
\end{equation}
These two groups of conditions have distinct roles.
The sign conditions allow us to show that $\kappa_g$ is a lower bound for the quartic quotient on magnetic tori, which we then transfer
to free cells. The contact conditions give
$\kappa_g=b_\triangle$ by Poisson summation.
Together with the upper bound
$\beta_*\leqslant b_\triangle$ from Lemma~\ref{lem:cell},
they therefore imply
\begin{equation}\label{cert:eq:minorant-implication}
\text{\eqref{eq:signconds} and \eqref{eq:latticecontactconds}}
\quad\Longrightarrow\quad
\beta_*\geqslant\kappa_g=b_\triangle
\quad\Longrightarrow\quad
\beta_*=b_\triangle.
\end{equation}
Theorem~\ref{cert:thm:criterion} states this implication
precisely and also gives a quantitative bound for the
finite-cell constants.

\smallskip
\noindent\emph{Existence of the required minorant.}
The following external result supplies a function satisfying
both groups of conditions.

\begin{theorem}[Sharp Gaussian minorant,
{\cite[Theorem~1.2]{Atomic}}]
\label{thm:gaussiancertificate}
There exists a real radial function $g\in\Sch(\R^2)$
satisfying \eqref{eq:signconds} and
\eqref{eq:latticecontactconds}.
\end{theorem}

This is the case $\alpha=1$ of
\cite[Theorem~1.2]{Atomic}, with its $f_\alpha$ denoted
here by $g$. Its Fourier convention and covolume-one
triangular lattice agree with
\eqref{cert:eq:fourier-convention} and
\eqref{cert:eq:triangular-lattice}, so no further change
of scale is needed.

The sign conditions are required on the whole plane,
because the weighted Poisson identity in
Lemma~\ref{cert:lem:poisson} evaluates $g$ and $\hat g$
on rectangular lattices of arbitrary side lengths.
The proof of the conditional statement uses only the
stated properties of $g$; the external theorem supplies
its existence. The variational limit proved in
Sections~\ref{sec:proofs}--\ref{sec:assembly}, with the
intrinsic coefficient $\beta_*$, is independent of this
existence result.

}

\subsection{The lower-bound criterion and deduction of
Theorem~\ref{thm:abrikosov}}
\label{cert:sec:criterion}\label{sec:coefficientproof}

We now state the conditional result that connects the
Gaussian minorant to the coefficient $\beta_*$.
The argument has three steps: the sign conditions
\eqref{eq:signconds} give a lower bound on magnetic tori;
cutoff and projection transfer this bound to free cells;
and the contact conditions \eqref{eq:latticecontactconds}
identify the resulting constant with $b_\triangle$.
The external theorem is used only to supply a function
satisfying these conditions.

\begin{theorem}[Fourier lower bound from a Gaussian minorant]
\label{cert:thm:criterion}
Let $g\in\Sch(\R^2)$ be real and radial, satisfy
\eqref{eq:signconds}, and let
$\kappa_g=1+\hat g(0)-g(0)$. Then:
\begin{enumerate}
\item[(i)] \emph{Periodic lower bound.}
We have $\kappa_g\geqslant0$. For every $L_1,L_2>0$
with $N=L_1L_2\in\mathbb N$ and every nonzero
$u\in\mathcal L^{\mathbb T}_{L_1,L_2}$,
\begin{equation}\label{cert:eq:periodic-bound}
\frac{N\int_Q\abs{u}^4}
     {\left(\int_Q\abs{u}^2\right)^2}
\geqslant\kappa_g,
\qquad Q=(0,L_1)\times(0,L_2).
\end{equation}

\item[(ii)] \emph{Free-cell lower bound.}
The constants defined in \eqref{eq:cell_unit_density}
satisfy
\begin{equation}\label{cert:eq:free-liminf}
\liminf_{L\to\infty}\tilde c_L\geqslant\kappa_g.
\end{equation}

\item[(iii)] \emph{Identification of the coefficient.}
If $g$ also satisfies \eqref{eq:latticecontactconds},
then $\kappa_g=b_\triangle$, and there are constants
$C,L_0>0$ such that
\begin{equation}\label{cert:eq:rate}
b_\triangle-CL^{-1/5}
\leqslant\tilde c_L\leqslant b_\triangle
\qquad(L\geqslant L_0).
\end{equation}
In particular, $\tilde c_L\to b_\triangle$ as $L\to\infty$,
and $\beta_*=b_\triangle$.
\end{enumerate}
\end{theorem}

\smallskip
\noindent\emph{Structure of the proof.}
The next three subsections prove assertions
(i)--(iii), respectively.
For (i), Lemma~\ref{cert:lem:basis} expresses the quartic
energy as a Gaussian sum with nonnegative weights.
The weighted Poisson identity of
Lemma~\ref{cert:lem:poisson} then allows both sign conditions
\eqref{eq:signconds} to be applied to this sum, giving
Proposition~\ref{cert:prop:periodic}.

For (ii), Proposition~\ref{cert:prop:free} cuts off a
free-cell state and projects it onto the lowest Landau
level of a slightly larger square torus.
Lemma~\ref{cert:lem:projection} controls the projection
error uniformly in the torus size. Applying (i) to the
projected state gives the free-cell bound, with a
quantitative error.

\begin{figure}[tp]
\centering
\includegraphics{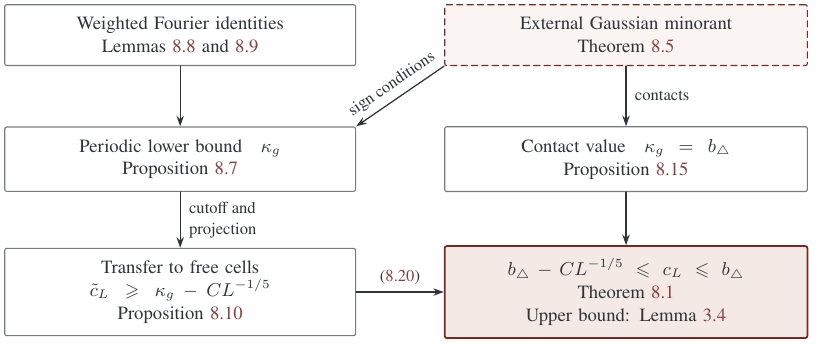}
\caption{Dependencies in the coefficient identification.
The sign conditions give the periodic and free-cell lower
bounds; the contact conditions identify their coefficient
with $b_\triangle$. The external existence theorem,
shown in the dashed box, supplies a function satisfying
both groups of conditions. The variational limit of
Sections~\ref{sec:proofs}--\ref{sec:assembly} is established
independently, with the intrinsic coefficient $\beta_*$.}
\label{fig:cert-roadmap}
\end{figure}

Finally, Proposition~\ref{cert:prop:contact} uses the
contact conditions and Poisson summation to obtain
$\kappa_g=b_\triangle$. Combining this identity with
the quantitative free-cell bound proves (iii).
Figure~\ref{fig:cert-roadmap} summarizes these dependencies.

\smallskip
\noindent\emph{Deduction of Theorem~\ref{thm:abrikosov}.}
Before proving the criterion, we show how it gives the
coefficient identity and finite-cell estimate announced
in Theorem~\ref{thm:abrikosov}.

\begin{proof}[Proof of Theorem~\ref{thm:abrikosov}]
Theorem~\ref{thm:gaussiancertificate} supplies a real
radial Schwartz function satisfying both
\eqref{eq:signconds} and \eqref{eq:latticecontactconds}.
Hence Theorem~\ref{cert:thm:criterion}(iii) applies.
Returning to the original cell scale through
\eqref{eq:certscale}, we obtain, for all sufficiently
large $L$,
\begin{equation}\label{cert:eq:original-cell-rate}
b_\triangle-C\pi^{1/10}L^{-1/5}
\overset{\eqref{cert:eq:rate}}{\leqslant}
\tilde c_{L/\sqrt\pi}
\overset{\eqref{eq:certscale}}{=}
c_L
\overset{\eqref{eq:certscale}}{\leqslant}
b_\triangle.
\end{equation}
Absorbing $\pi^{1/10}$ into the constant gives
\eqref{eq:cellrate}. Letting $L\to\infty$ and recalling
the intrinsic characterization of $\beta_*$ yields
\begin{equation}\label{cert:eq:coefficient-conclusion}
\beta_*
\overset{\eqref{eq:intrinsiccell}}{=}
\sup_{L>0}c_L
=\lim_{L\to\infty}c_L
\overset{\eqref{cert:eq:original-cell-rate}}{=}
b_\triangle,
\end{equation}
which proves \eqref{eq:identified}.
\end{proof}

\smallskip
\noindent\emph{Consequences for the variational limit.}
The compactness, lower-bound, and recovery assertions
of Theorem~\ref{thm:gamma} were proved in
Section~\ref{sec:compactness} and
Subsections~\ref{sec:liminfproof}
and~\ref{sec:recoveryproof}, with the intrinsic coefficient
$\beta_*$. Once the criterion is proved below,
Theorem~\ref{thm:abrikosov} completes the identification
of their limiting functional.
In particular, \eqref{eq:limit} and \eqref{eq:TF} become
\begin{equation}\label{cert:eq:explicit-limit}
\begin{aligned}
\F(\rho)
&=2\int_{\R^2}|x|^2\rho\,\mathrm dx
+b_\triangle\int_{\R^2}\rho^2\,\mathrm dx,\\
\rho_*(x)
&=\frac1{b_\triangle}
\left(\sqrt{\frac{2b_\triangle}{\pi}}-|x|^2\right)_+,
\end{aligned}
\end{equation}
where $\F$ is finite on the admissible class specified
in \eqref{eq:limit}.
The conclusions of Theorem~\ref{thm:min} and
Proposition~\ref{prop:original} therefore hold with
$\beta_*=b_\triangle$.

\subsection{Periodic bound: proof of assertion (i) in Theorem~\ref{cert:thm:criterion}}
\label{cert:sec:lll}

We prove the lower bound on magnetic tori that will be
transferred to free cells in Proposition~\ref{cert:prop:free}.
The argument expresses the quartic energy as a Gaussian
sum and uses the sign conditions \eqref{eq:signconds}
to bound this sum from below.

\begin{proposition}[Periodic lower bound]
\label{cert:prop:periodic}
Let $g\in\Sch(\R^2)$ be real and radial and satisfy
the sign conditions \eqref{eq:signconds}, namely (recall that $G(x)=e^{-\pi|x|^2}$)
\begin{equation}
g(x)\leqslant G(x),
\qquad
\hat g(\bm\xi)\geqslant0
\qquad(x,\bm\xi\in\R^2).
\end{equation}
Set $\kappa_g:=1+\hat g(0)-g(0)$, as in
\eqref{cert:eq:Bg}. Then $\kappa_g\geqslant0$.
Moreover, for every $L_1,L_2>0$ with
$N=L_1L_2\in\mathbb N$, every nonzero
$u\in\mathcal L^{\mathbb T}_{L_1,L_2}$
satisfies the periodic lower bound
\eqref{cert:eq:periodic-bound}:
\begin{equation}
\frac{N\int_Q\abs{u}^4\,\mathrm dx}
     {\left(\int_Q\abs{u}^2\,\mathrm dx\right)^2}
\geqslant\kappa_g,
\qquad Q=(0,L_1)\times(0,L_2).
\end{equation}
\end{proposition}

The proof uses two lemmas. We first state them and use
them to prove the proposition; their proofs are deferred
to the next two subsubsections.

Both lemmas involve the following weights. Given
$c_0,\ldots,c_{N-1}\in\C$ with
$\sum_{j=0}^{N-1}|c_j|^2=1$, extend the coefficients 
$N$-periodically, i.e. $c_{N+j}=c_j$ for all $j \in \Z$, and set
\begin{equation}\label{cert:eq:A}
W_{a,b}:=\left|\sum_{j=0}^{N-1}
c_{j+b}\overline{c_j}e^{-2\pi iaj/N}\right|^2,
\qquad a,b\in\Z.
\end{equation}
These weights are $N$-periodic in both indices, and
Cauchy--Schwarz gives
\begin{equation}\label{cert:eq:weight-bounds}
0\leqslant W_{a,b}\leqslant W_{0,0}=1.
\end{equation}
For $Q=(0,L_1)\times(0,L_2)$, write
\begin{equation}\label{cert:eq:freq}
\bm k_{a,b}:=\left(\frac a{L_1},\frac b{L_2}\right),
\qquad a,b\in\Z.
\end{equation}
Although $W_{a,b}$ depends on the sequence $(c_r)_{r\in \Z}$, we suppress this dependence in the notation as it poses no risk of confusion.

The following lemma constructs an orthonormal basis of the
torus LLL and expresses the quartic energy in its coefficients.
The construction also proves that the dimension is $N$.
\begin{lemma}[Basis and quartic identity]
\label{cert:lem:basis}
Let $L_1,L_2>0$, $N=L_1L_2\in\mathbb N$,
$Q=(0,L_1)\times(0,L_2)$, and $\mu=2^{1/4}L_2^{-1/2}$.
The functions
\begin{equation}\label{cert:eq:basis}
  \varphi_r(x_1,x_2)
  =\mu\sum_{\ell\in\Z}
  e^{2\pi i(r+\ell N)x_2/L_2}
  e^{-\pi(x_1-(r+\ell N)/L_2)^2},
  \qquad 0\leqslant r<N,
\end{equation}
form an orthonormal basis of
$\mathcal L^{\mathbb T}_{L_1,L_2}$.
In particular, $\dim\mathcal L^{\mathbb T}_{L_1,L_2}=N$.
Every $u\in\mathcal L^{\mathbb T}_{L_1,L_2}$ with
$\norm{u}_{L^2(Q)}=1$ has a unique expansion
\begin{equation}
  u=\sum_{r=0}^{N-1}c_r\varphi_r,
  \qquad
  \sum_{r=0}^{N-1}\abs{c_r}^2=1.
\end{equation}
Extending the coefficients by $c_{r+N}=c_r$, we have
\begin{equation}\label{cert:eq:quartic}
    N\int_Q\abs{u}^4\,\mathrm dx
    = \sum_{a,b\in\Z}G(\bm k_{a,b})
    \left| \sum_{j=0}^{N-1}
    c_{j+b}\overline{c_j}e^{-2\pi iaj/N} \right|^2
    =\sum_{a,b\in\Z}G(\bm k_{a,b})W_{a,b},
\end{equation}
where $G(x)=e^{-\pi\abs{x}^2}$,
$\bm k_{a,b}=(a/L_1,b/L_2)$, and $W_{a,b}$ are the
weights \eqref{cert:eq:A} associated with the $N$ periodic coefficients $(c_r)_{r\in \Z}$.
\end{lemma}

The second lemma replaces the weighted sum of $g$
by the corresponding sum of $\hat g$. The condition
$\hat g\geqslant0$ then allows us to retain only
the contribution at the origin.

\begin{lemma}[Weighted Poisson identity]
\label{cert:lem:poisson}
Let $L_1L_2=N\in\mathbb N$, with $W_{a,b}$ and
$\bm k_{a,b}$ defined by \eqref{cert:eq:A} and
\eqref{cert:eq:freq}. For every radial $g\in\Sch(\R^2)$,
\begin{equation}\label{cert:eq:weighted}
\sum_{a,b\in\Z}g(\bm k_{a,b})W_{a,b}
=\sum_{a,b\in\Z}\hat g(\bm k_{a,b})W_{a,b}.
\end{equation}
Both sums converge absolutely.
\end{lemma}

We now combine Lemmas~\ref{cert:lem:basis}
and~\ref{cert:lem:poisson} to prove
Proposition~\ref{cert:prop:periodic}, which gives
assertion~(i) of Theorem~\ref{cert:thm:criterion}. The proof of 
Lemmas~\ref{cert:lem:basis} and~\ref{cert:lem:poisson} will be provided afterward.

\begin{proof}[Proof of Proposition~\ref{cert:prop:periodic}]
Because of the sign conditions \eqref{eq:signconds}, we have $g(0)\leqslant G(0)=1$
and $\hat g(0)\geqslant0$.
Hence, $\kappa_g=1-g(0)+\hat g(0)\geqslant0$.

By homogeneity, we may assume $\norm{u}_{L^2(Q)}=1$.
Combining the quartic formula \eqref{cert:eq:quartic}
with the weighted Poisson identity \eqref{cert:eq:weighted},
and separating the term at the origin, where
$G(0)=W_{0,0}=1$, we obtain
\begin{equation}\label{cert:eq:periodic-defect}
  N\int_Q\abs{u}^4\,\mathrm dx-\kappa_g
  =
  \sum_{(a,b)\neq(0,0)}
  \bigl[(G-g)(\bm k_{a,b})+\hat g(\bm k_{a,b})\bigr]
  W_{a,b}
  \geqslant0.
\end{equation}
Indeed, every summand is nonnegative by
\eqref{eq:signconds} and \eqref{cert:eq:weight-bounds}.
All sums converge absolutely by Schwartz decay and
$0\leqslant W_{a,b}\leqslant1$.
Restoring the normalization proves
\eqref{cert:eq:periodic-bound}.
\end{proof}

\subsubsection{Proof of Lemma~\ref{cert:lem:basis}:
the quartic Fourier formula}

We use the standard theta-function construction of the
torus LLL; see \cite[Section~5]{Onofri2001} and
\cite[Section~VII, equation~(7.11)]{SakamotoTanimura2003}.
These constructions give the basis \eqref{cert:eq:basis}
after the corresponding change of gauge and normalization.
We therefore restrict ourselves to checking the normalization
and deriving the quartic identity in our notation.

\begin{proof}
We first check that \eqref{cert:eq:basis} is an orthonormal system.
Integrating in $x_2$ and then combining the translated
intervals in $x_1$, we obtain
\begin{equation}\label{cert:eq:basis-normalization}
  \ip{\varphi_r}{\varphi_s}_{L^2(Q)}
  =\delta_{rs}\mu^2L_2
   \int_{\R}e^{-2\pi t^2}\,\mathrm dt
  =\delta_{rs},
\end{equation}
where we used $\mu^2L_2=\sqrt2$.
Next, let $u=\sum_{r=0}^{N-1}c_r\varphi_r$ satisfy
$\norm{u}_{L^2(Q)}=1$, and extend the coefficients
$N$-periodically. We can then write
\begin{equation}\label{cert:eq:state-series}
  u(x_1,x_2)
  =\mu\sum_{m\in\Z}c_m
   e^{2\pi imx_2/L_2}e^{-\pi(x_1-m/L_2)^2}.
\end{equation}
We will compute the Fourier coefficients of the periodic
density $\abs{u}^2$ and then apply Parseval's identity.
We deniote by $I_{a,b}$ the unnormalized Fourier coefficient
of the periodic density $|u|^2$ at frequency $(a/L_1,b/L_2)$:
\begin{equation}\label{cert:eq:density-coefficients}
  I_{a,b}
  =\int_Q\abs{u(x)}^2
   e^{-2\pi i\bm k_{a,b}\cdot x}\,\mathrm dx,
  \qquad
  \bm k_{a,b}=\left(\frac a{L_1},\frac b{L_2}\right).
\end{equation}
Gaussian decay and boundedness of the periodic sequence
$(c_m)_{m\in \Z}$ justify the termwise integrations below.

Fix $(a,b)\in\Z^2$. We substitute
\eqref{cert:eq:state-series} into
\eqref{cert:eq:density-coefficients} and integrate in $x_2$.
By orthogonality of the exponentials, only the pairs
$m+b,m$, with $m\in\Z$, remain.
We write each $m\in\Z$ uniquely as
$m=j+\ell N$, where $j\in\{0,\ldots,N-1\}$ and $\ell\in\Z$,
and set $t=x_1-m/L_2=x_1-j/L_2-\ell L_1$.
As $x_1$ ranges over $(0,L_1)$, the variable $t$ ranges over
\begin{equation}
  \left(-\frac j{L_2}-\ell L_1,\,
        -\frac j{L_2}-(\ell-1)L_1\right).
\end{equation}
For each fixed $j\in\{0,\ldots,N-1\}$, these intervals,
indexed by $\ell\in\Z$, are consecutive and cover $\R$
up to their endpoints.
Since $(c_m)_{m\in\Z}$ is $N$-periodic and $N=L_1L_2$,
we have
\begin{equation}
  c_{m+b}\overline{c_m}
  =c_{j+b}\overline{c_j},
  \qquad
  e^{-2\pi iax_1/L_1}
  =e^{-2\pi iat/L_1}e^{-2\pi iaj/N},
\end{equation}
where the factor $e^{-2\pi ia\ell}$ equals $1$ because
$a,\ell\in\Z$. The transformed integrand is therefore
independent of $\ell$, and we can combine the corresponding
integrals into one integral over $\R$.
Using $\mu^2L_2=\sqrt2$, we obtain
\begin{equation}\label{cert:eq:unfolded-density}
  I_{a,b}
  =\sum_{j=0}^{N-1}
   c_{j+b}\overline{c_j}e^{-2\pi iaj/N}
   \sqrt2\int_{\R}
   e^{-\pi[t^2+(t-b/L_2)^2]}
   e^{-2\pi iat/L_1}\,\mathrm dt.
\end{equation}
To evaluate the integral in \eqref{cert:eq:unfolded-density},
we rewrite the quadratic expression in $t$ as
\begin{equation}
  t^2+\left(t-\frac b{L_2}\right)^2
  =2\left(t-\frac b{2L_2}\right)^2
   +\frac{b^2}{2L_2^2}.
\end{equation}
We can thus take $e^{-\pi b^2/(2L_2^2)}$ outside the
integral. Setting $s=t-b/(2L_2)$, we obtain
\begin{equation}
    \sqrt2\int_{\R}
    e^{-\pi[t^2+(t-b/L_2)^2]}
    e^{-2\pi iat/L_1}\,\mathrm dt
   =
    e^{-\pi b^2/(2L_2^2)}e^{-\pi iab/N}
    \sqrt2\int_{\R}
    e^{-2\pi s^2}e^{-2\pi ias/L_1}\,\mathrm ds.
\end{equation}
Here the phase $e^{-\pi iab/N}$ comes from substituting
$t=s+b/(2L_2)$ into the oscillatory exponential and
using $L_1L_2=N$. 
We now apply the Gaussian Fourier formula
\begin{equation}\label{cert:eq:gaussian-fourier}
  \sqrt2\int_{\R}e^{-2\pi s^2}e^{-2\pi i\xi s}\,\mathrm ds
  =e^{-\pi\xi^2/2},
  \qquad \xi\in\R,
\end{equation}
with $\xi=a/L_1$. Since
$\abs{\bm k_{a,b}}^2=a^2/L_1^2+b^2/L_2^2$,
substitution into \eqref{cert:eq:unfolded-density} gives
\begin{equation}\label{cert:eq:Ifourier}
  I_{a,b}
  =e^{-\pi\abs{\bm k_{a,b}}^2/2}e^{-\pi iab/N}
   \sum_{j=0}^{N-1}
   c_{j+b}\overline{c_j}e^{-2\pi iaj/N},
  \qquad (a,b)\in\Z^2.
\end{equation}
Finally, the Fourier coefficients of the periodic density
$\abs{u}^2$, with the usual normalization by
$\abs{Q}=N$, are $I_{a,b}/N$. We therefore apply
Parseval's identity and use \eqref{cert:eq:Ifourier}
and the definition \eqref{cert:eq:A} of $W_{a,b}$:
\begin{equation}\label{cert:eq:quartic-parseval}
  N\int_Q\abs{u}^4\,\mathrm dx
  =\sum_{a,b\in\Z}\abs{I_{a,b}}^2
  =\sum_{a,b\in\Z}G(\bm k_{a,b})W_{a,b}.
\end{equation}
This proves \eqref{cert:eq:quartic}.
\end{proof}

\subsubsection{Proof of Lemma~\ref{cert:lem:poisson}:
the weighted Poisson identity}

We now prove \eqref{cert:eq:weighted}, which allows us
to use $\hat g\geqslant0$ in the periodic lower bound.
The finite Fourier identity \eqref{cert:eq:selfdual}
is a special case of Sussman's identity on finite groups;
see \cite[Proposition~4.4]{FeichtingerKozekLuef2009}.
We combine it with ordinary Poisson summation and the
radiality of $g$ to get
\eqref{cert:eq:weighted}, i.e., that
for every radial $g\in\Sch(\R^2)$,
\begin{equation}
\sum_{a,b\in\Z}g(\bm k_{a,b})W_{a,b}
=\sum_{a,b\in\Z}\hat g(\bm k_{a,b})W_{a,b}.
\end{equation}. We give the proof in our notation
for the reader's convenience.

\begin{proof}
We extend $(c_j)_{j=0}^{N-1}$ $N$-periodically to $\Z$.
All infinite sums below converge absolutely, since
$g,\hat g\in\Sch(\R^2)$ and $0\leqslant W_{a,b}\leqslant1$.

\smallskip
\noindent\emph{Step 1. A finite Fourier identity for the weights.}
We fix $(a,b)\in\Z^2$ and use the orthogonality relation
\begin{equation}\label{cert:eq:characters}
  \frac1N\sum_{r=0}^{N-1}e^{2\pi irk/N}
  =
  \begin{cases}
    1,&k\in N\Z,\\
    0,&k\notin N\Z,
  \end{cases}
  \qquad k\in\Z.
\end{equation}
By the definition \eqref{cert:eq:A}, we have
\begin{equation}
  W_{r,s}
  =\sum_{j,k=0}^{N-1}
   c_{j+s}\overline{c_j}\,
   \overline{c_{k+s}}c_k
   e^{2\pi ir(k-j)/N},
  \qquad r,s\in\{0,\ldots,N-1\}.
\end{equation}
After multiplication by $e^{2\pi i(br-as)/N}$,
the sum in $r$ vanishes unless $k-j+b\in N\Z$.
We can therefore replace $k$ by $j-b$, reading all
coefficient indices modulo $N$, and obtain
\begin{equation}
  \frac1N\sum_{r,s=0}^{N-1}
  W_{r,s}e^{2\pi i(br-as)/N}
  =
  \sum_{s,j=0}^{N-1}
  c_{j+s}\overline{c_j}\,
  \overline{c_{j-b+s}}c_{j-b}e^{-2\pi ias/N}.
\end{equation}
We now set $k=j+s$ modulo $N$. Since $k$ ranges over
$\{0,\ldots,N-1\}$ as $s$ does, we can factor the sum:
\begin{equation}
  \begin{aligned}
    \sum_{s,j=0}^{N-1}
    c_{j+s}\overline{c_j}\,
    \overline{c_{j-b+s}}c_{j-b}e^{-2\pi ias/N}
    &=
    \left(\sum_{k=0}^{N-1}
    c_k\overline{c_{k-b}}e^{-2\pi iak/N}\right)
    \overline{\left(\sum_{j=0}^{N-1}
    c_j\overline{c_{j-b}}e^{-2\pi iaj/N}\right)}\\
    &=
    \left|\sum_{j=0}^{N-1}
    c_j\overline{c_{j-b}}e^{-2\pi iaj/N}\right|^2
    =W_{a,b}.
  \end{aligned}
\end{equation}
For the last equality, we shift $j$ by $b$ and use
\eqref{cert:eq:A}; the resulting factor
$e^{-2\pi iab/N}$ has modulus one. We have proved
\begin{equation}\label{cert:eq:selfdual}
  W_{a,b}
  =\frac1N\sum_{r,s=0}^{N-1}
   W_{r,s}e^{2\pi i(br-as)/N},
  \qquad (a,b)\in\Z^2.
\end{equation}

\smallskip
\noindent\emph{Step 2. Poisson summation.}
Fix $r,s\in\{0,\ldots,N-1\}$.
We apply the Poisson sumamtion formula \eqref{cert:eq:poisson-general} to the lattice
$L_1^{-1}\Z\times L_2^{-1}\Z$, whose covolume is $1/N$.
This gives
\begin{equation}\label{cert:eq:rectangular-poisson}
  \frac1N\sum_{a,b\in\Z}
  g(\bm k_{a,b})e^{2\pi i(br-as)/N}
  =
  \sum_{\ell,m\in\Z}
  \hat g\left(L_1\ell+\frac s{L_2},
              L_2m-\frac r{L_1}\right).
\end{equation}
The prefactor $1/N$ cancels the reciprocal covolume
arising in Poisson summation.
Since $g$ is radial, so is $\hat g$.
We can therefore rotate each argument by $\pi/2$:
\begin{equation}
  \begin{aligned}
    \hat g\left(L_1\ell+\frac s{L_2},
                L_2m-\frac r{L_1}\right)
    &=
    \hat g\left(\frac r{L_1}-L_2m,\,
                \frac s{L_2}+L_1\ell\right)\\
    &=\hat g\left(\bm k_{r-Nm,s+N\ell}\right),
    \qquad (\ell,m)\in\Z^2.
  \end{aligned}
\end{equation}
Reindexing the sum by $(\ell,m)\mapsto(-m,\ell)$,
we conclude that
\begin{equation}
  \frac1N\sum_{a,b\in\Z}
  g(\bm k_{a,b})e^{2\pi i(br-as)/N}
  =
  \sum_{\ell,m\in\Z}
  \hat g(\bm k_{r+N\ell,s+Nm}).
\end{equation}
We now insert \eqref{cert:eq:selfdual} into the weighted
sum of $g$ and apply the preceding identity:
\begin{equation}\label{cert:eq:weighted-calculation}
  \begin{aligned}
    \sum_{a,b\in\Z}g(\bm k_{a,b})W_{a,b}
    &=
    \sum_{r,s=0}^{N-1}W_{r,s}
    \left[
    \frac1N\sum_{a,b\in\Z}
    g(\bm k_{a,b})e^{2\pi i(br-as)/N}
    \right]\\
    &=
    \sum_{r,s=0}^{N-1}\sum_{\ell,m\in\Z}
    W_{r,s}\hat g(\bm k_{r+N\ell,s+Nm})\\
    &=\sum_{a,b\in\Z}\hat g(\bm k_{a,b})W_{a,b}.
  \end{aligned}
\end{equation}
Indeed, every $(a,b)\in\Z^2$ has a unique representation
$a=r+N\ell$, $b=s+Nm$, with
$r,s\in\{0,\ldots,N-1\}$ and $\ell,m\in\Z$,
and periodicity gives $W_{a,b}=W_{r,s}$.
This proves \eqref{cert:eq:weighted}.
\end{proof}

\subsection{Transfer to free cells: proof of assertion (ii)}
\label{cert:sec:transfer}

To prove assertion~(ii) of Theorem~\ref{cert:thm:criterion},
we transfer the periodic lower bound to free cells.
We first cut off a free-cell state near the boundary
so that it admits a magnetic-periodic extension.
We then project this extension onto the torus LLL
to restore the first-order equation.
Before stating the quantitative estimates, we explain
the construction and the errors we need to control.

Let $g\in\Sch(\R^2)$ be real and radial, satisfying
the sign conditions \eqref{eq:signconds}, and set
$\kappa_g:=1+\hat g(0)-g(0)$.
For $L_1,L_2\in(0,\infty)$ with
$N=L_1L_2\in\mathbb N$, the periodic lower bound
\eqref{cert:eq:periodic-bound} gives
\begin{equation}
  \frac{N\int_Q\abs{v}^4\,\mathrm dx}
       {\left(\int_Q\abs{v}^2\,\mathrm dx\right)^2}
  \geqslant\kappa_g,
  \qquad
  v\in\mathcal L^{\mathbb T}_{L_1,L_2}\setminus\{0\},
  \qquad Q=(0,L_1)\times(0,L_2).
\end{equation}
We want to obtain the corresponding lower bound for
the free-cell constant $\tilde c_L$.
For $L>0$, its competitors are
\begin{equation}
  u(z)=e^{-\pi\abs{z}^2/2}f(z),
  \qquad
  f\in\mathcal O(Q_L)\cap L^2(Q_L;\C)\setminus\{0\},
  \qquad Q_L=(-L/2,L/2)^2,
\end{equation}
where we identify $x=(x_1,x_2)\in\R^2$ with
$z=x_1+ix_2\in\C$.
There are two obstructions: these states satisfy no boundary conditions, and
the area of $Q_L$, i.e., $L^2$, need not be an integer, as required for
magnetic periodicity. We therefore cannot
apply the periodic estimate directly to $u$.

We first choose a containing square of integer area\footnote{Any $S\geqslant L$ with $S^2\in\mathbb N$
would suffice. Our choice gives
$0\leqslant S^2-L^2<1$ and hence
$1\leqslant S^2/L^2<1+L^{-2}$.
The alternative $S=\lceil L\rceil$ gives
$S^2/L^2=1+O(L^{-1})$ as $L\to\infty$,
which also suffices for the final $L^{-1/5}$ bound.}:
\begin{equation}
  Q_S=(-S/2,S/2)^2,
  \qquad S=\sqrt{\lceil L^2\rceil},
  \qquad S^2\in\mathbb N.
\end{equation}
This is the smallest such enlargement of $Q_L$.

The free-cell states are written in symmetric gauge,
whereas the torus estimate and projection are formulated
in Landau gauge. We therefore multiply $u$ by
$e^{i\pi x_1x_2}$, preserving its modulus and quartic
quotient. The transformed state satisfies
\begin{equation}
  A^{\mathbb T}\bigl(e^{i\pi x_1x_2}u\bigr)=0
  \qquad\text{in }\mathcal D'(Q_L).
\end{equation}
We then multiply it by a cutoff
$\chi\in C_c^\infty(Q_L;\R)$, extend the product by zero
to $Q_S$, and extend it magnetically to $\R^2$.
This gives a smooth magnetic-periodic function $v$,
which generally no longer satisfies $A^{\mathbb T}v=0$.
We correct this error by projecting $v$ onto
$\mathcal L^{\mathbb T}_{S,S}$.

The estimates below show that, for a suitable cutoff
and large $L$, the projected state is nonzero and its
quartic quotient is close to that of $u$.
We first state the quantitative lower bound and the
projection lemma. We then prove the proposition using
the lemma, whose proof follows afterwards.
Only the sign conditions \eqref{eq:signconds} are used.

\begin{proposition}[Quantitative free-cell lower bound]
\label{cert:prop:free}
Let $g\in\Sch(\R^2)$ be real and radial, with
$g\leqslant G$ and $\hat g\geqslant0$ on $\R^2$.
Set $\kappa_g=1+\hat g(0)-g(0)$ as in
\eqref{cert:eq:Bg}.
There exist constants $C_g>0$ and $L_g\geqslant1$
such that
\begin{equation}\label{cert:eq:generic-rate}
  \tilde c_L\geqslant\kappa_g-C_gL^{-1/5}
  \qquad\text{for every }L\in[L_g,\infty).
\end{equation}
Consequently,
$\liminf_{L\to\infty}\tilde c_L\geqslant\kappa_g$.
\end{proposition}

To control the projection error, we estimate the
remainder $v-P^{\mathbb T}v$ through
$\norm{A^{\mathbb T}v}_{L^2(Q_S^+)}$.
This quantity vanishes precisely when $v$ already
belongs to the torus LLL.

\begin{lemma}[Projection estimates]
\label{cert:lem:projection}
Let $N\in\mathbb N$, $S=\sqrt N\geqslant1$, and
$Q_S^+=(0,S)^2$.
We identify magnetic-periodic functions with their
restrictions to $Q_S^+$ and denote by $P^{\mathbb T}$
the $L^2(Q_S^+;\C)$-orthogonal projection onto
$\mathcal L^{\mathbb T}_{S,S}$.

For every $v\in C^\infty(\R^2;\C)$ satisfying
the magnetic-periodic conditions \eqref{cert:eq:bc}
with $L_1=L_2=S$, the remainder
$r=(I-P^{\mathbb T})v$ satisfies
\begin{equation}\label{cert:eq:gap-est}
    \norm{r}_{L^2(Q_S^+)}
    \leqslant(4\pi)^{-1/2}
      \norm{A^{\mathbb T}v}_{L^2(Q_S^+)},\qquad
    \norm{D^{\mathbb T}r}_{L^2(Q_S^+)}^2
    \leqslant\frac32
      \norm{A^{\mathbb T}v}_{L^2(Q_S^+)}^2,
\end{equation}
and
\begin{equation}\label{cert:eq:L4proj}
  \norm{r}_{L^4(Q_S^+)}
  \leqslant C\norm{A^{\mathbb T}v}_{L^2(Q_S^+)},
\end{equation}
where $C>0$ is independent of $S$ and $v$.
Here
$D^{\mathbb T}=(\partial_{x_1},
\partial_{x_2}-2\pi ix_1)$ and
$A^{\mathbb T}=D^{\mathbb T}_1+iD^{\mathbb T}_2$.
\end{lemma}

\begin{remark}[Uniformity in the torus size]
\label{cert:rem:spectral-gap}
We will apply the lemma with $S\to\infty$, so its
constants must be independent of $S$.
The estimates \eqref{cert:eq:gap-est} follow from
the Landau-level gap $4\pi$, which does not depend
on the torus size.
To obtain \eqref{cert:eq:L4proj}, we apply a periodic
interpolation inequality to $\abs{r}$ and use
the diamagnetic inequality
\begin{equation}
  \abs{\nabla\abs{r}(x)}
  \leqslant\abs{D^{\mathbb T}r(x)}
  \qquad\text{for almost every }x\in\R^2.
\end{equation}
The resulting interpolation constant is uniform
for $S\geqslant1$.
\end{remark}

\begin{figure}[tp]
\centering
\includegraphics[]{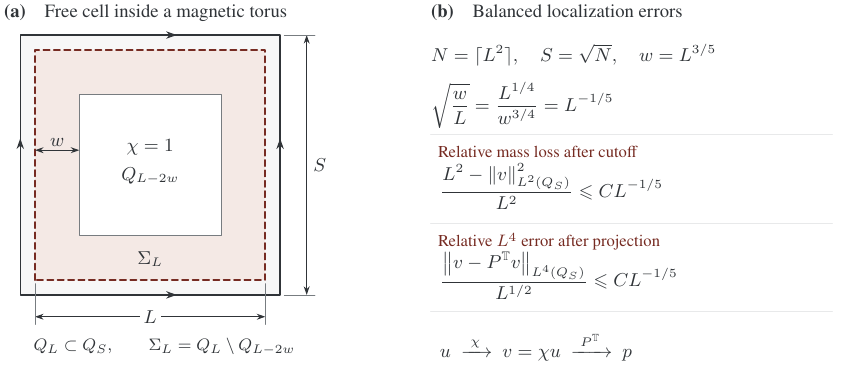}
\caption{Transfer from a free cell to a magnetic torus.
(a) The cutoff equals one on $Q_{L-2w}$ and vanishes
near $\partial Q_L$, so the truncated state extends
smoothly to the magnetic torus.
The geometric widths are schematic.
(b) The bounds for the relative cutoff and projection
errors are proportional to $(w/L)^{1/2}$ and
$L^{1/4}w^{-3/4}$, respectively.
We plot these expressions without their multiplicative
constants, using horizontal variable $w/L^{3/5}$ and
multiplying the vertical values by $L^{1/5}$.
Both expressions equal $L^{-1/5}$ when $w=L^{3/5}$.}
\label{fig:boundary}
\end{figure}

\subsubsection{Proof of Proposition~\ref{cert:prop:free}}

We construct a torus LLL state whose mass and quartic
integral differ little from those of the free-cell
competitor. We then apply the periodic lower bound
to this state. Related cutoff and projection arguments
appear in \cite[proof of Proposition~2.1 and
Section~4.2]{AS07}.

\begin{proof}
We take $L$ sufficiently large.
The constant $C>0$ may change from line to line,
but is independent of $L$ and of the competitor $u$.

\smallskip
\noindent\emph{Step 1. Cutoff and extension.}
By homogeneity and the upper bound
$\tilde c_L\leqslant b_\triangle$ in
\eqref{eq:certscale}, we can restrict the infimum
defining $\tilde c_L$ to admissible states satisfying
\begin{equation}\label{cert:eq:nearmin}
  \int_{Q_L}\abs{u}^2\,\mathrm dx=L^2,
  \qquad
  \int_{Q_L}\abs{u}^4\,\mathrm dx
  \leqslant(b_\triangle+1)L^2.
\end{equation}
Indeed, states with normalized quartic integral larger
than $b_\triangle+1$ cannot lower this infimum.
We pass to Landau gauge by replacing $u$ with
$e^{i\pi x_1x_2}u$ and, with a harmless abuse of notation,
continue to denote the transformed function by $u$. Its modulus is unchanged,
and $A^{\mathbb T}u=0$ on $Q_L$.

We choose a square containing $Q_L$ whose area is an integer:
\begin{equation}\label{cert:eq:transfer-scales}
  S=\sqrt{\lceil L^2\rceil},
  \qquad L^2\leqslant S^2<L^2+1,
  \qquad Q_S=(-S/2,S/2)^2.
\end{equation}
We also set
\begin{equation}\label{cert:eq:transfer-error-scales}
  w=L^{3/5},
  \qquad
  \sqrt{\frac wL}=L^{-1/5}.
\end{equation}
For large $L$, we have $1\leqslant w<L/2$.
The reason for this choice of $w$ is explained in
Remark~\ref{cert:ex:optimal-w} below (see also 
Remark~\ref{cert:rem:flux-quantization} for a comment
on the choice of $S$).

Take $\chi\in C_c^\infty(Q_L;\R)$ such that
\begin{equation}\label{cert:eq:cutoff}
  \begin{aligned}
    &0\leqslant\chi(x)\leqslant1,\qquad
      \abs{\nabla\chi(x)}\leqslant Cw^{-1}
      &&(x\in Q_L),\\
    &\chi(x)=1
      &&(x\in Q_{L-2w}).
  \end{aligned}
\end{equation}
We define $v=\chi u$ on $Q_L$, extend it by zero
to $Q_S$, and then extend it magnetically to $\R^2$
by \eqref{cert:eq:bciter}, with $L_1=L_2=S$.
The extension is compatible because $S^2\in\mathbb N$.
It is smooth because $u$ is smooth inside $Q_L$
and $\chi u$ vanishes near $\partial Q_S$.

We use $Q_S$ for the torus integrals below.
This gives the same inner product as integration over
$Q_S^+=(0,S)^2$, since the boundary phases cancel
in the product of one state with the complex conjugate
of another. The moduli of covariant derivatives are
also periodic. Thus Lemma~\ref{cert:lem:projection}
applies with all its norms taken over $Q_S$.

Let $\Sigma_L:=Q_L\setminus Q_{L-2w}$.
The cutoff changes $u$ only in this strip, whose area
is $4Lw-4w^2$, which is bounded above by $4Lw$. Cauchy--Schwarz
therefore gives
\begin{equation}\label{cert:eq:strip-mass}
  \int_{\Sigma_L}\abs{u}^2\,\mathrm dx
  \leqslant\abs{\Sigma_L}^{1/2}
    \left(\int_{Q_L}\abs{u}^4\,\mathrm dx\right)^{1/2}
  \overset{\eqref{cert:eq:nearmin}}{\leqslant}
    CL^{3/2}w^{1/2}
    =CL^{9/5}.
\end{equation}
Hence, since $\chi=1$ on $Q_{L-2w}$ and $v=\chi u$ vanishes on $Q_S\setminus Q_L$, we obtain
\begin{equation}\label{cert:eq:massloss}
  0\leqslant L^2-\norm{v}_{L^2(Q_S)}^2
  =\int_{Q_L}(1-\chi^2)\abs{u}^2\,\mathrm dx
  \leqslant\int_{\Sigma_L}\abs{u}^2\,\mathrm dx
  \leqslant CL^{9/5}.
\end{equation}
By the product rule we see that the cutoff
changes the LLL equation as follows:
\begin{equation}
  A^{\mathbb T}v(x)
  =\bigl(\partial_{x_1}\chi(x)
         +i\partial_{x_2}\chi(x)\bigr)u(x),
  \qquad x\in Q_L,
\end{equation}
and since $\nabla\chi=0$ on $Q_{L-2w}$, we have
\begin{equation}\label{cert:eq:defect}
  \norm{A^{\mathbb T}v}_{L^2(Q_S)}^2
  =\int_{\Sigma_L}\abs{\nabla\chi}^2
    \abs{u}^2\,\mathrm dx
  \overset{\eqref{cert:eq:cutoff},
             \eqref{cert:eq:strip-mass}}{\leqslant}
    Cw^{-2}L^{9/5}
    =CL^{3/5}.
\end{equation}

\smallskip
\noindent\emph{Step 2. Projection and comparison.}
We first check that $v$ satisfies the hypotheses of
Lemma~\ref{cert:lem:projection}.
Since $u\in C^\infty(Q_L;\C)$ and
$\chi\in C_c^\infty(Q_L;\R)$, we have
$\chi u\in C_c^\infty(Q_L;\C)$.
Its zero extension to $Q_S$ vanishes near
$\partial Q_S$, so its magnetic-periodic extension
$v$ belongs to $C^\infty(\R^2;\C)$.
We set $p=P^{\mathbb T}v$ and $r=v-p$.
Lemma~\ref{cert:lem:projection}, applied on $Q_S$,
and \eqref{cert:eq:defect} give
\begin{equation}\label{eq:sec8restimates24}
  \norm{r}_{L^2(Q_S)}^2\leqslant CL^{3/5},
  \qquad
  \norm{r}_{L^4(Q_S)}\leqslant CL^{3/10}.
\end{equation}
Since $p$ and $r$ are orthogonal in $L^2(Q_S;\C)$,
we obtain, using \eqref{cert:eq:massloss},
\begin{equation}\label{cert:eq:pmass}
  \norm{p}_{L^2(Q_S)}^2
  =\norm{v}_{L^2(Q_S)}^2-\norm{r}_{L^2(Q_S)}^2
  \geqslant L^2-CL^{9/5}-CL^{3/5}
  \geqslant L^2(1-CL^{-1/5}).
\end{equation}
In particular, $p\neq0$ for large $L$. For the quartic integral, we use the triangle inequality, the second estimate in \eqref{eq:sec8restimates24}, and the fact that $\norm{v}_{L^4(Q_S)}\leqslant\norm{u}_{L^4(Q_L)}$ to get:
\begin{equation}\label{cert:eq:pquartic}
  \norm{p}_{L^4(Q_S)}
  \leqslant\norm{u}_{L^4(Q_L)}+CL^{3/10}.
\end{equation}
By \eqref{cert:eq:nearmin}, namely $\norm{u}^4_{L^4(Q_L)}\leqslant (b_\triangle+1) L^{2}$, we have the bound $\norm{u}_{L^4(Q_L)}\leqslant CL^{1/2}$. Expanding the fourth power of \eqref{cert:eq:pquartic} yields:
\begin{equation}\label{cert:eq:transfer-normalized}
\begin{aligned}
  \norm{p}_{L^4(Q_S)}^4
  &\leqslant \left(\norm{u}_{L^4(Q_L)} + CL^{3/10}\right)^4 
  = \int_{Q_L}\abs{u}^4\,\mathrm dx + \mathcal{O}\!\left(\norm{u}_{L^4(Q_L)}^3 L^{3/10}\right) \\
  &\leqslant \int_{Q_L}\abs{u}^4\,\mathrm dx + CL^{9/5},
\end{aligned}
\end{equation}
where we used that the largest cross term is bounded by $C(L^{1/2})^3L^{3/10}=CL^{9/5}$ and all remaining lower-order terms grow strictly slower.
We now apply the periodic lower bound to $p\in\mathcal L^{\mathbb T}_{S,S}\setminus\{0\}$. Using \eqref{cert:eq:transfer-scales}, \eqref{cert:eq:pmass}, and \eqref{cert:eq:transfer-normalized}, we obtain
\begin{equation}\label{cert:eq:projected-lower}
  \begin{split}
    \kappa_g
    &\overset{\eqref{cert:eq:periodic-bound}}{\leqslant}
      \frac{S^2\norm{p}_{L^4(Q_S)}^4}
           {\norm{p}_{L^2(Q_S)}^4}\\
    &\leqslant
      \frac{(1+L^{-2})
        \left(L^{-2}\int_{Q_L}\abs{u}^4\,\mathrm dx
              +CL^{-1/5}\right)}
           {(1-CL^{-1/5})^2}\\
    &\leqslant
      \frac1{L^2}\int_{Q_L}\abs{u}^4\,\mathrm dx
      +CL^{-1/5}.
  \end{split}
\end{equation}
For the last inequality, we observe that the prefactor satisfies $(1+L^{-2})(1-CL^{-1/5})^{-2} = 1 + \mathcal{O}(L^{-1/5})$ for large $L$. Since the integral term $L^{-2}\int_{Q_L}\abs{u}^4 \leqslant b_\triangle+1$ is uniformly bounded, expanding the product yields the main integral term plus an error of order $\mathcal{O}(L^{-1/5})$, which is simply absorbed into the generic constant $C$.

Taking the infimum over all competitors satisfying \eqref{cert:eq:nearmin} yields $\tilde c_L\geqslant\kappa_g-CL^{-1/5}$. This establishes \eqref{cert:eq:generic-rate}, and letting $L\to\infty$ proves the last assertion.
\end{proof}

\begin{remark}[Choice of the cutoff width]
\label{cert:ex:optimal-w}
We explain the choice $w=L^{3/5}$ by keeping
$w\in[1,L/2)$ arbitrary in the cutoff construction
\eqref{cert:eq:cutoff}, with $u$ satisfying
\eqref{cert:eq:nearmin}. As in the proof, we extend
$v=\chi u$ by zero to $Q_S$ and then magnetically
to $\R^2$, and set $r=v-P^{\mathbb T}v$.
We compare two errors.

The first is the \emph{relative mass loss}
$(L^2-\norm{v}_{L^2(Q_S)}^2)/L^2$.
The strip estimate \eqref{cert:eq:strip-mass} and
the identity \eqref{cert:eq:massloss}, before
substituting $w=L^{3/5}$, give
\begin{equation}\label{cert:eq:width-mass-loss}
  0\leqslant
  \frac{L^2-\norm{v}_{L^2(Q_S)}^2}{L^2}
  \leqslant C\left(\frac wL\right)^{1/2}.
\end{equation}
The second is the \emph{relative projection error}
$\norm{r}_{L^4(Q_S)}/\norm{u}_{L^4(Q_L)}$,
which measures the projection remainder relative
to the original state in $L^4$.
The cutoff derivative bound \eqref{cert:eq:cutoff}
and the strip estimate \eqref{cert:eq:strip-mass},
used in \eqref{cert:eq:defect}, give
$\norm{A^{\mathbb T}v}_{L^2(Q_S)}^2
\leqslant CL^{3/2}w^{-3/2}$.
Combining this with \eqref{cert:eq:L4proj}
and the lower bound
$\norm{u}_{L^4(Q_L)}\geqslant L^{1/2}$,
which follows from Cauchy--Schwarz and the
normalization in \eqref{cert:eq:nearmin}, we obtain
\begin{equation}\label{cert:eq:width-projection-error}
  \frac{\norm{r}_{L^4(Q_S)}}{\norm{u}_{L^4(Q_L)}}
  \leqslant CL^{1/4}w^{-3/4}.
\end{equation}
Increasing $w$ worsens the mass-loss bound
\eqref{cert:eq:width-mass-loss} but improves the
projection bound \eqref{cert:eq:width-projection-error}.
For $w=L^\alpha$, their right-hand sides are
$CL^{(\alpha-1)/2}$ and $CL^{(1-3\alpha)/4}$,
respectively. Both tend to zero when
$\alpha\in(1/3,1)$.
We obtain the best common decay exponent by
setting $(\alpha-1)/2=(1-3\alpha)/4$,
which gives $\alpha=3/5$.
Thus $w=L^{3/5}$ bounds both relative errors
by $CL^{-1/5}$. This choice optimizes the two
bounds above; it does not establish that the
convergence rate is sharp.
\end{remark}

\begin{remark}[The area adjustment]
\label{cert:rem:flux-quantization}
Magnetic periodicity requires $S^2\in\mathbb N$,
while the inclusion $Q_L\subseteq Q_S$ requires
$S\geqslant L$. Our choice
$S=\sqrt{\lceil L^2\rceil}$ gives the smallest
such square and satisfies
$0\leqslant S^2-L^2<1$, hence
$0\leqslant S^2/L^2-1<L^{-2}$.

We could instead take $S=\lceil L\rceil$.
Then $S\geqslant L$ and $S^2\in\mathbb N$,
but the relative area increase satisfies only
$0\leqslant S^2/L^2-1<2L^{-1}+L^{-2}$.
The proof would still work: in
\eqref{cert:eq:projected-lower}, the area factor
would contribute an error of order $L^{-1}$
rather than $L^{-2}$, because the normalized
quartic integral is uniformly bounded by
\eqref{cert:eq:nearmin}.
Both errors are smaller than $L^{-1/5}$, so
either choice gives the same final rate.
\end{remark}

\subsubsection{Proof of Lemma~\ref{cert:lem:projection}}

We estimate the error caused by projection onto the torus LLL.
We first use the spectral gap to control its $L^2$ norm
and magnetic derivatives. We then obtain the $L^4$ estimate
by interpolation, keeping the constants independent of $S$.

\begin{proof}
Throughout the proof, all norms are taken over
$Q_S^+=(0,S)^2$. 
Let $r=v-P^{\mathbb T}v$.
Since $P^{\mathbb T}v$ is a finite linear combination
of the smooth basis functions \eqref{cert:eq:basis},
the function $r$ is smooth. By the definition of
the orthogonal projection, we also have
\begin{equation}
  r\perp\mathcal L^{\mathbb T}_{S,S}
  \quad\text{in }L^2(Q_S^+;\C),
  \qquad
  A^{\mathbb T}r=A^{\mathbb T}v.
\end{equation}

\smallskip
\noindent\emph{Step 1. The $L^2$ estimates.}
We recall the factorization
\eqref{cert:eq:torus-factorization}: every smooth
magnetic-periodic function $h$ satisfies
\begin{equation}\label{cert:eq:Dsplit}
  \norm{D^{\mathbb T}h}_{L^2}^2
  =\norm{A^{\mathbb T}h}_{L^2}^2
   +2\pi\norm{h}_{L^2}^2.
\end{equation}
Thus every smooth function $h$ satisfying the magnetic-periodic boundary conditions obeys
\begin{equation}
  \int_{Q_S^+}\abs{D^{\mathbb T}h}^2\,\mathrm dx
  \geqslant
  2\pi\int_{Q_S^+}\abs{h}^2\,\mathrm dx,
\end{equation}
with equality if and only if $A^{\mathbb T}h=0$,
that is, $h\in\mathcal L^{\mathbb T}_{S,S}$.

We now use the spectral decomposition of the magnetic
Laplacian associated with the quadratic form
$\norm{D^{\mathbb T}h}_{L^2}^2$ on magnetic-periodic
$H^1$ functions. Its eigenfunctions form a complete
orthonormal basis, and its eigenvalues are
$2\pi+4\pi n$, with $n\in\{0,1,\ldots\}$;
see \cite[Section~VII, equation~(7.12)]{SakamotoTanimura2003},
where the operator includes an additional factor $1/2$.
By \eqref{cert:eq:Dsplit}, the eigenspace at $2\pi$
is exactly $\mathcal L^{\mathbb T}_{S,S}$.
For each $n\in\{0,1,\ldots\}$, let $r_n$ be the component
of $r$ in the eigenspace at $2\pi+4\pi n$.
Since $r=v-P^{\mathbb T}v$ is orthogonal to
$\mathcal L^{\mathbb T}_{S,S}$, we have $r_0=0$.
Consequently,
\begin{equation}
  r=\sum_{n=1}^{\infty}r_n
  \quad\text{in }L^2(Q_S^+;\C),
  \qquad
  \norm{r}_{L^2}^2
  =\sum_{n=1}^{\infty}\norm{r_n}_{L^2}^2.
\end{equation}
The same spectral decomposition gives
\begin{equation}
  \norm{D^{\mathbb T}r}_{L^2}^2
  =\sum_{n=1}^{\infty}(2\pi+4\pi n)
    \norm{r_n}_{L^2}^2
  \geqslant6\pi\norm{r}_{L^2}^2.
\end{equation}
Here we used $n\geqslant1$: after removing the lowest
Landau level, the smallest remaining eigenvalue is $6\pi$.

We subtract $2\pi\norm{r}_{L^2}^2$ and apply
\eqref{cert:eq:Dsplit}. Observe that
$A^{\mathbb T}P^{\mathbb T}v=0$, so
$A^{\mathbb T}r=A^{\mathbb T}v$. We therefore obtain
\begin{equation}\label{cert:eq:projection-spectral-gap}
    4\pi\norm{r}_{L^2}^2
    \leqslant
      \norm{D^{\mathbb T}r}_{L^2}^2
      -2\pi\norm{r}_{L^2}^2
    \overset{\eqref{cert:eq:Dsplit}}{=}
      \norm{A^{\mathbb T}r}_{L^2}^2
      =\norm{A^{\mathbb T}v}_{L^2}^2.
\end{equation}
Taking square roots proves the first estimate in
\eqref{cert:eq:gap-est}. For the second, we return to
\eqref{cert:eq:Dsplit}:
\begin{equation}
    \norm{D^{\mathbb T}r}_{L^2}^2
    =\norm{A^{\mathbb T}v}_{L^2}^2
      +2\pi\norm{r}_{L^2}^2
    \overset{\eqref{cert:eq:projection-spectral-gap}}
      {\leqslant}
      \norm{A^{\mathbb T}v}_{L^2}^2
      +\frac12\norm{A^{\mathbb T}v}_{L^2}^2
      =\frac32\norm{A^{\mathbb T}v}_{L^2}^2.
\end{equation}
Both constants are independent of $S$, because the
Landau-level spacing is $4\pi$ for every admissible
torus size.

\smallskip
\noindent\emph{Step 2. The $L^4$ estimate.}
For every $h\in H^1_{\mathrm{loc}}(\R^2;\C)$
that is $S$-periodic in each coordinate, we have
\begin{equation}\label{cert:eq:ladyzhenskaya}
  \norm{h}_{L^4}^2
  \leqslant C\norm{h}_{L^2}
  \left(\norm{\nabla h}_{L^2}
        +S^{-1}\norm{h}_{L^2}\right),
\end{equation}
with $C$ independent of $S$.
For $S=1$, we obtain this estimate by applying the planar
Gagliardo--Nirenberg inequality to a fixed smooth cutoff
of the periodic extension of $h$;
see \cite[Section~8.5]{LiebLoss2001}.
For general $S>0$, we apply the unit-square estimate
to $h(S\,\cdot)$ and change variables.

The magnetic boundary conditions imply that $\abs{r}$
is $S$-periodic in each coordinate:
\begin{equation}
  \abs{r(x_1+S,x_2)}
  =\abs{r(x_1,x_2+S)}
  =\abs{r(x_1,x_2)},
  \qquad (x_1,x_2)\in\R^2.
\end{equation}
Moreover, the diamagnetic inequality gives
\begin{equation}
  \abs{r}\in H^1_{\mathrm{loc}}(\R^2;\R),
  \qquad
  \abs{\nabla\abs{r}(x)}
  \leqslant\abs{D^{\mathbb T}r(x)}
  \quad\text{for almost every }x\in\R^2.
\end{equation}
We can therefore apply \eqref{cert:eq:ladyzhenskaya}
to $\abs{r}$ on $Q_S^+=(0,S)^2$.
Using $S\geqslant1$ and \eqref{cert:eq:gap-est}, we obtain
\begin{equation}
    \norm{r}_{L^4(Q_S^+)}^2
    \leqslant C\norm{r}_{L^2(Q_S^+)}
      \left(\norm{D^{\mathbb T}r}_{L^2(Q_S^+)}
            +S^{-1}\norm{r}_{L^2(Q_S^+)}\right)
    \overset{\eqref{cert:eq:gap-est}}{\leqslant}
      C\norm{A^{\mathbb T}v}_{L^2(Q_S^+)}^2.
\end{equation}
Taking square roots proves \eqref{cert:eq:L4proj},
with $C$ independent of $S$ and $v$.
\end{proof}

\subsection{Contact value and rate: proof of assertion (iii)}
\label{cert:sec:contact}\label{cert:sec:rate}

We complete the criterion by showing that the contact conditions identify
$\kappa_g$ with $b_\triangle$. The free-cell estimate of
Proposition~\ref{cert:prop:free} then gives the convergence rate.
The identification of $\kappa_g$ uses only Poisson summation and the
contact conditions.

\begin{proposition}[Contact value and rate]
\label{cert:prop:contact}\label{cert:prop:rate}
Let $g\in\Sch(\R^2)$ satisfy the contact conditions
\eqref{eq:latticecontactconds}, namely
\begin{equation}
  \begin{aligned}
    g(\bm p)&=G(\bm p)
      &&(\bm p\in\Lambda_\triangle\setminus\{0\}),\\
    \hat g(\bm q)&=0
      &&(\bm q\in\Lambda_\triangle^*\setminus\{0\}).
  \end{aligned}
\end{equation}
Then $\kappa_g=1+\hat g(0)-g(0)=b_\triangle$.

If, in addition, $g$ is real and radial and satisfies the sign conditions
\eqref{eq:signconds}, that is, $g\leqslant G$ and $\hat g\geqslant0$
on $\R^2$, then
\begin{equation}
  b_\triangle-C_gL^{-1/5}
    \leqslant\tilde c_L\leqslant b_\triangle
  \qquad(L\geqslant L_g),
\end{equation}
with the constants of Proposition~\ref{cert:prop:free}.
In particular, $\lim_{L\to\infty}\tilde c_L=b_\triangle$.
\end{proposition}

\begin{proof}
All lattice sums below converge absolutely. Since
$\Lambda_\triangle$ has covolume one, Poisson summation and the
contact conditions give
\begin{equation}
  \begin{aligned}
    \kappa_g
      &\overset{\eqref{cert:eq:Bg}}{=}1+\hat g(0)-g(0)
      \overset{\eqref{cert:eq:poisson-general},
                 \eqref{eq:latticecontactconds}}{=}
        1+\sum_{\bm p\in\Lambda_\triangle\setminus\{0\}}g(\bm p)\\
      &\overset{\eqref{eq:latticecontactconds}}{=}
        1+\sum_{\bm p\in\Lambda_\triangle\setminus\{0\}}G(\bm p)
      \overset{\eqref{cert:eq:btri-lattice},
                 \eqref{cert:eq:G}}{=}b_\triangle.
  \end{aligned}
\end{equation}
Under the additional hypotheses, Proposition~\ref{cert:prop:free}
applies. Thus, for $L\geqslant L_g$,
\begin{equation}
  b_\triangle-C_gL^{-1/5}
    =\kappa_g-C_gL^{-1/5}
    \overset{\eqref{cert:eq:generic-rate}}{\leqslant}\tilde c_L
    \overset{\eqref{eq:certscale}}{\leqslant}b_\triangle,
\end{equation}
which also proves the limit.
\end{proof}

\begin{proof}[Proof of Theorem~\ref{cert:thm:criterion}]
Assertions~(i), (ii), and~(iii) follow directly from
Propositions~\ref{cert:prop:periodic}, \ref{cert:prop:free}, and~\ref{cert:prop:contact}, respectively. 
This establishes the criterion under the stated hypotheses on $g$. 

To deduce Theorem~\ref{thm:abrikosov} from this criterion, it remains only to
show that there exists a function satisfying both the sign and contact conditions.
The existence of such a function was established by OpenAI in \cite[Theorem~1.2]{Atomic}
(recalled here as Theorem~\ref{thm:gaussiancertificate}),
which completes the proof.
\end{proof}

\bigskip
\subsection*{Acknowledgements}
\textsc{G.D.F.} is a member of GNAMPA--INdAM. He gratefully acknowledges 
partial financial support from the GNAMPA Project CUP~E53C25002010001, 
and from the University of Naples Federico II through the FRA Project-B 
``VarMoCry'' on \emph{Variational Analysis and Modeling of Liquid Crystals}.

\subsection*{Declarations} \emph{AI Usage:} 
This work is the result of a long-term research
effort by the authors, who developed all the
mathematical strategies underlying the paper.
The identification of the effective interaction
coefficient relies on the sharp Gaussian-minorant
theorem of OpenAI \cite[Theorem~1.2]{Atomic}.
AI assistants (Claude Opus 5.5 and Gemini Pro 3.1)
were used to review intermediate calculations
and check the consistency of the authors' arguments.
The authors take full responsibility for the
content of the paper.
\emph{Data Availability:}
No data are associated with this manuscript.
\emph{Conflict of Interest:}
The authors declare no conflicts of interest.
\emph{Author Contributions:}
All authors contributed equally to this work.

\raggedbottom
\bibliographystyle{amsplain}
\bibliography{references}
\end{document}